\documentclass[11pt,A4]{article}
\usepackage{amsmath,amssymb,amsthm,epsfig}

\newcommand\nc\newcommand
\nc\calA{\mathcal{A}}\nc\calC{\mathcal{C}}
\nc\calD{\mathcal{D}}\nc\calF{\mathcal{F}}
\nc\calH{\mathcal{H}}\nc\calI{\mathcal{I}}
\nc\calJ{\mathcal{J}}\nc\calL{\mathcal{L}}
\nc\calM{\mathcal{M}}\nc\calO{\mathcal{O}}
\nc\calQ{\mathcal{Q}}\nc\calR{\mathcal{R}}
\nc\calS{\mathcal{S}}\nc\calT{\mathcal{T}}
\nc\calU{\mathcal{U}}\nc\calV{\mathcal{V}}
\renewcommand\gg{\mathfrak g}
\nc\bfL{\mathbf L}
\nc\ccc{\mathbf c}\nc\qqq{\mathbf q}
\nc\rrr{\mathbf r}\nc\sss{\mathbf s}
\nc\ttt{\mathbf t}\nc\uuu{\mathbf u}
\nc\dC{\mathbb C}\nc\dN{\mathbb N}
\nc\dR{\mathbb R}\nc\dZ{\mathbb Z}
\renewcommand\a{\alpha}\renewcommand\b{\beta}\renewcommand\d{\delta}

\renewcommand\l{\lambda}
\renewcommand\o{\omega}\nc\G{\Gamma}\nc\ga{\gamma}
\nc\D{\Delta}\nc\e{\varepsilon}
\nc\eps{\epsilon}\nc\s{\sigma}\nc\Ups{\Upsilon}
\nc\dd{{\rm d}}
\nc\rmleft{{\rm left}}
\nc\rmright{{\rm right}}
\nc\Rep{{\rm Rep}}
\nc\Lie{{\rm Lie}}
\nc\Prim{{\rm Prim}}
\nc\rmu{{\rm u}}
\nc\can{{\rm can}}
\nc\triv{{\rm triv}}
\nc\inv{{\rm inv}}
\nc\univ{{\rm univ}}\nc\gen{{\rm gen}}
\nc\Ext{{\rm Ext}}
\nc\Hom{{\rm Hom}}\nc\Lin{{\rm Lin}}
\nc\End{{\rm End}}\nc\Coend{{\rm Coend}}

\nc\im{{\rm im}}
\nc\tr{{\rm tr}}
\nc\id{{\rm id}}\nc\Id{{\rm Id}}
\nc\Mat{{\rm Mat}}
\nc\Aut{{\rm Aut}}
\nc\m{{\rm m}}
\nc\ad{{\rm ad}}\nc\Ad{{\rm Ad}}
\nc\w{\wedge}

\nc\rt{\triangleright}
\nc\lt{\triangleleft}
\nc\brt{\blacktriangleright}
\nc\blt{\blacktriangleleft}
\nc{\covdot}{\,\underline{.}\,}
\nc{\lcross}{\!>\!\!\!\!\triangleleft}
\nc{\bicross}{{\,\triangleright\!\triangleleft\,}}
\nc\tgg{\tilde{\gg}}
\nc\ts{\tilde{\sigma}}
\nc\td{\tilde{\delta}}
\nc\te{\tilde{\e}}
\nc\tket{\tilde{]}}
\nc\htheta{\hat{\theta}}
\nc\surject{\to\hspace{-0.8em}\to}
\nc\uA{\underline{A}}
\nc\uH{\underline{H}}
\nc\um{\underline{{\rm m}}}
\nc\ueta{\underline{\eta}}
\nc\ue{\underline{\e}}
\nc\uD{\underline{\D}}
\nc\uS{\underline{S}}
\nc\utr{\underline{\tr}}
\nc\uotimes{\,\underline{\otimes}\,}
\nc\1{\underline{1}}\nc\2{\underline{2}}
\nc\multi[5]{{}_{#2}^{#4}\! {#1}_{#3}^{#5}}
\nc\ceqn[2]{\begin{equation}\label{#1}\begin{array}{c}#2\end{array}
    \end{equation}}
\nc\cqfd{\hfill$\diamond$}
\nc\ie{{\it i.e.}~}
\nc\Proof{{\it Proof. }}
\nc\Remark{{\it Remark. }}
\nc\ba{\begin{array}}
\nc\ea{\end{array}}

\newtheorem{lemma}{Lemma}[section]
\newtheorem{definition}[lemma]{Definition}
\newtheorem{proposition}[lemma]{Proposition}
\newtheorem{theorem}[lemma]{Theorem}
\newtheorem{corollary}[lemma]{Corollary}

\begin{document}
\title{Braided Lie algebras and bicovariant differential calculi
    \\ over co-quasitriangular Hopf algebras.}
\author{X. Gomez\footnote{EU Marie Curie Fellow}\;
and S. Majid\footnote{Royal Society University Research Fellow}}
\date{}
\maketitle
\abstract{ We show that if $\gg_\G$ is the quantum tangent space
(or quantum Lie algebra in the sense of Woronowicz) of a
bicovariant first order differential calculus over a
co-quasitriangular Hopf algebra $(A,\rrr)$, then a certain
extension of it is a braided Lie algebra in the category of
$A$-comodules. This is used to show that the Woronowicz quantum
universal enveloping algebra $U(\gg_\G)$ is a bialgebra in the
braided category of $A$-comodules. We show that this algebra is
quadratic when the calculus is inner. Examples with this
unexpected property include finite groups and quantum groups with
their standard differential calculi. We also find a quantum Lie
functor for co-quasitriangular Hopf algebras, which has properties
analogous to the classical one. This functor gives trivial results
on standard quantum groups $\calO_q(G)$, but reasonable ones on
examples closer to the classical case, such as the cotriangular
Jordanian deformations. In addition, we show that split braided
Lie algebras define `generalised-Lie algebras' in a different
sense of deforming the adjoint representation. We construct these
and their enveloping algebras for $\calO_q(SL(n))$, recovering the
Witten algebra for $n=2$. }
\section*{Introduction.}
In a well-known article \cite{Woronowicz}, Woronowicz has given an
axiomatic treatment of so called bicovariant first order
differential calculi (FODC) over Hopf algebras. It appeared that,
given a arbitrary Hopf algebra $A$, there is no canonical way to
associate to it a bicovariant FODC. Nevertheless, each bicovariant
FODC $(\G,\dd)$ over $A$ extends to a graded differential algebra
$(\G^\w,\dd)$ - later shown to be have a Hopf superalgebra
structure \cite{Brzezinski}\cite{Schauenburg}, and has an
associated ``quantum Lie algebra'' $\gg_\G$. A quantum Lie algebra
is a vector space $\gg$ equipped with a braiding operator $\sigma$
and a ``quantum Lie bracket'' satisfying certain identities, which
coincide with the usual Lie algebra axioms when the braiding $\s$
on $\gg$ is the usual flip. This article is concerned with the
universal enveloping algebra $U(\gg)$ of a quantum Lie algebra (a
certain quotient of the tensor algebra of $\gg$). The main
question which we address is whether  $U(\gg)$ can be equipped
with a Hopf algebra or bialgebra structure. To do this we will
need additional coassociative structure $\delta:\gg\to
\gg\otimes\gg$, which we axiomatise as a `good quantum Lie
algebra' $(\gg,\sigma,[\ ,\,],\delta)$ (we then show that in the
case of $\gg_\G$ associated to an FODC there is a canonical such
structure when $A$ is coquasitriangular.)

Let us recall what are the obstructions for a bialgebra structure on
$U(\gg_\G)$.
 Before asking in which sense our required coalgebra
structure maps $\D$ and $\e$  should be algebra morphisms, the
coproduct on $U(\gg_\G)$ should be coassociative. Woronowicz
\cite{Woronowicz} has shown that bicovariant FODC over $A$ are
parametrized by (in our conventions) $\ad_L$-invariant left ideals
of $A$ contained in $\ker \e_A$. The ``cotangent space'' of the
FODC $(\G,\dd)$ corresponding to such an ideal $\calI$ is
canonically identified with $\ker\e_A/\calI$, and the quantum Lie
algebra $\gg_\G$ (or tangent space) is the dual of this. When $A$
is commutative, $\ker\e_A/\calI$ has a natural associative algebra
structure, say $\mu$ (often without unit). Therefore the dual
$\gg_\G$ has a natural comultiplication $\d=\mu^*$, and the space
$k\, 1\oplus\gg_\G$ has a coalgebra structure $(\D,\e)$ such that
$\D(1)=1\otimes 1$, $\D(x)=x\otimes 1+1\otimes x+\d(x)$,
$\e(x)=0$, for $x\in \gg_\G$. Actually, for the standard
differential calculus on Lie groups, the ideal $\calI$ in question
is $(\ker\e_A)^2$, therefore $\mu=0$, $\d=0$, and we recover the
standard coproduct on $U(\gg_\G)$. In the non commutative case,
these things do not work anymore. The key idea to solve this
problem is the observation that, when $A$ has a co-quasitriangular
$\rrr$, $\ad_L$-invariant left ideals of $A$ are two-sided ideals
of another algebra, namely $\uA$, the braided version (or
`braided-group') of the quantum group $A$ in the category of (in
our conventions, left) $A$-comodules, introduced in
\cite{Majid-braided-group}. Therefore, in this case the space
$\tgg_\G=k\, 1\oplus \gg_\G$ inherits naturally a coalgebra
structure, which is dual to that of the unital algebra
$\uA/\calI$. It remains to show that this comultiplication does
extend to a bialgebra structure on $U(\gg_\G)$, and to say in
which sense. The answer (theorem \ref{main-theo}) is that
$U(\gg_\G)$ indeed becomes a bialgebra, in the category of (right)
$A$-comodules, which is braided thanks to the existence of $\rrr$.

The reason for all this to work, and which was the starting point
of the article, is the theory of braided Lie algebras. These were
introduced in \cite{Majid-braided-Lie} as axiomatically-defined
finite-dimensional objects $\calL$ with axioms strong enough to
define a braided enveloping algebra $B(\calL)$  as a quadratic
bialgebra in a braided category (it was previously denoted
${\calU}(\calL)$ in \cite{Majid-braided-Lie}). Moreover, for
standard quantum groups $U_q(\gg)$ there is an algebra map
$B(\calL)\to U_q(\gg)$ so that $U_q(\gg)$ is essentially generated
by $\calL$ (the `Lie problem for quantum groups'). In this case
braided Lie algebra $\calL$ is of an $n^2$-dimensional matrix form
and can in fact be identified\cite{Majid-classif-bicov} with the
quantum tangent space for the FODC on $\calO_q(G)$ constructed by
the R-matrix method of Jurco \cite{Jurco}. However, a general
theorem systematically linking the Woronowicz theory and the
braided Lie algebra theory has been missing and is what we provide
now. Indeed, in our construction $\tgg_\G=k\, 1\oplus \gg_\G$ is a
braided Lie algebra in $\calM^A$ and its braided universal
enveloping bialgebra $B(\tgg_\G)$ therefore provides the
homogenization or quadratic central extension of $U(\gg_\G)$, with
$U(\gg_\G)$ a bialgebra quotient of $B(\tgg_\G)$.

Thus we prove the existence of a (braided) bialgebra structure on
$U(\gg_\G)$ for arbitrary bicovariant FODC over $A$, provided $A$
is co-quasitriangular. The existence of an antipode is more
problematic, and actually we can prove that an antipode does not
exist in many examples of interest, such as for all finite
dimensional bicovariant calculi on standard quantum groups
$\calO_q(G)$, $q$ not a root of unity, and all bicovariant calculi
on finite groups.  The reason is quite curious~: We prove in
theorem \ref{theo-quadratic} (which is general, \ie does not use
co-quasitriangularity) that when the differential $\dd$ is
implemented by a biinvariant element $\theta$, then $U(\gg_\G)$ is
a quadratic algebra, in fact a quantum symmetric algebra. In the
co-quasitriangular inner case, we find that $U(\gg_\G)$ is a quadratic
bialgebra. More precisely,
$$
U(\gg_\G)\simeq B(\calL)\qquad (\G\;\;{\rm inner})
$$
is the braided universal enveloping algebra of a braided Lie subalgebra
$\calL\subset \tgg_\G$, of the same dimension as $\gg_\G$. Since
$B(\calL)$ never has an antipode, $U(\gg_\G)$ does not have one as well.
Moreover, for the simple bicovariant FODC over standard quantum
groups $\calO_q(G)$ where the braided Lie algebra
$\calL$ in question is a matrix braided Lie algebra, its braided
universal enveloping algebra is an algebra of braided matrices
$B(R)$\cite{Majid-book}. Therefore $U(\gg_\G)$ is also an algebra
of braided matrices.

The algebra $\uA$, which is a key ingredient of our construction,
also proves useful to obtain an analogue of the Lie functor for
Lie groups -proposition \ref{Lie-functor}. Our functor goes from
the category of co-quasitriangular Hopf algebras (pairs
$(A,\rrr)$) to that of first order differential calculi (triples
$(A,\G,\dd)$). We show that it shares many of the properties of
the classical one; apparently too many because it sends to zero
most of the standard quantizations $\calO_q(G)$ (an exception is
$G=GL(n)$ for which the associated quantum Lie algebra has
dimension 1). Thus it provides another point of view on the non
triviality of these quantizations. We show in the case $G=SL(n)$
that for the softer (and less interesting) triangular deformations
of $\calO(G)$, the functor does give results which are close to
the classical ones (in particular with a reasonable dimension of
the quantum Lie algebra).

The last significant contribution of this paper is to establish a
relationship between braided Lie algebras and a third approach to
q-deformed Lie algebras defined by representation theory. Here
$\gg_q=\gg$ as a vector space but with the $q$-deformation of the
adjoint representation. For generic $q$ there remains a canonical
intertwiner $\gg_q\otimes \gg_q\to \gg_q$ which could be called
`Lie bracket'. Even though examples have been computed already for
some years\cite{Delius}, one does not know a full set of axioms
that the obtained $(\gg_q,[\ , \ ])$ should obey.  In this context
there is similarly a proposal\cite{LS} for an `enveloping algebra'
$U_{LS}(\gg)$ (say) associated to $\gg$ semisimple. An open
problem here, shown only for $\gg=sl(2)$, is to find some kind of
`homogenisation' of $U_{LS}(\gg)$ mapping onto (the locally finite
part of) $U_q(\gg)$ and thereby relating these algebras. A
corollary of the braided Lie theory is a solution of this problem
for $\gg=sl(n)$, as follows. We consider braided Lie algebras that
split as $\calL=kc\uplus\calL^+$, where $\calL^+$ is the kernel of
the counit of the braided Lie algebra and where we suppose that
$[c,\ ]$ acts as a multiple $\lambda$ of the identity. One can (in
principle) axiomatise the inherited properties of $\calL^+$ and
define its enveloping algebra as $B_{\rm
red}(\calL^+)=B(\calL)/\langle c-\lambda\rangle$. For the standard
matrix braided Lie algebra $\calL=\widetilde{sl_q(n)}$ associated
to $U_q(sl(n))$ we have $c=\utr$ (the quantum trace) and
$\calL^+=sl_q(n)$ (say) has the classical dimension $n^2-1$. The
enveloping $B_{\rm red}(sl_q(n))$ by construction has
homogenisation $B(\widetilde{sl_q(n)})$ mapping to (the locally
finite part of) $U_q(sl(n))$. It is clear already from
\cite{Majid-braided-Lie} that $B_{\rm
red}(sl_q(2))=U_{LS}(sl(2))$. We also mention that the latter is
isomorphic to the Witten algebra $W_{q^2}(sl(2))$ introduced in
\cite{Witten}, as was already noticed by L. Le Bruyn in
\cite{LeB-95}. This suggests a definition of $W_q(sl(n))$ for any
$n$ (although the physical requirements which lead to the
definition of $W_q(sl(2))$ are not considered here).
\\

Let us explain the content of the paper in more detail. In the
first preliminaries section we recall various well-known facts
about (co)-quasitriangular Hopf algebras, crossed modules and Hopf
bimodules, and quadratic bialgebras. Then in
section~\ref{sect-quantum}, we recall the definition of a quantum
Lie algebra, following theorems 5.3 and 5.4 in \cite{Woronowicz}.
(The notion of quantum Lie algebra is however not left-right
symmetric and we take the conventions opposite to
\cite{Woronowicz}, the main reason being that we want the quantum
Lie bracket to be given by the left adjoint action in the
differential calculi setting). We observe that the homogenization
of the universal enveloping algebra $U(\gg)$ of any quantum Lie
algebra $\gg$ is a quantum symmetric algebra. We then investigate
the existence of a coproduct on the universal enveloping algebra
of a quantum Lie algebra $(\gg,\s)$, not necessarily in the
context of differential calculi. For this, we suppose the
existence of an underlying braided category $(\calV,\Psi)$ in
which both $\gg$ and $U(\gg)$ live (the structure maps
$(\s,[\,,\,])$ of $\gg$ should be morphisms in $\calV$), and look
for a coproduct on $U(\gg)$ of the form $\D(x)=x\otimes 1+1\otimes
x+\d(x)$, $x\in\gg$, for some ``little coproduct'' $\d:\gg\to
\gg\otimes \gg$. Let us stress that not all quantum Lie algebras
can be equipped with such a little coproduct. Among those which
can, there is a subclass with nice properties, leading to our
notion of ``good'' quantum Lie algebras. An important feature of
these good quantum Lie algebras is that their braiding $\s$ is not
an essential datum~: it can be expressed in terms of the other
structure maps of $\gg$. Moreover, $\s$ can coincide with the
underlying braiding $\Psi_{\gg,\gg}$ only in some special cases
(which include super Lie algebras). Therefore a generic ``good
quantum Lie algebra'' is equipped with two braidings, the
categorical braiding $\Psi_{\gg,\gg}$, and the braiding $\s$,
which should not be confused.

Section \ref{sect-braided}, about braided Lie algebras, is mainly
taken from \cite{Majid-braided-Lie}, with slight improvements, in
particular on some properties of the canonical  braiding $\Ups$,
and on the connection with quantum Lie algebras. Recall that a
braided Lie algebra is already a coalgebra $(L,\D,\e)$ in a
braided category, endowed with a ``braided Lie bracket''
satisfying identities which also mimic usual Lie algebra axioms.
One of the differences is that they do not have an antisymmetry
axiom, and indeed, such an axiom is impossible to define in
general. However, one can consider the subclass of braided Lie
algebras $L$ which have a braided Lie algebra imbedding $k\to L$.
In this case, the Lie algebra-like object inside $L$ is $L^+$ (the
kernel of $\e$), and there is a natural notion of antisymmetry
axiom. We call ``good'' those braided Lie algebras which meet all
these requirements, and show that there is a 1-1 correspondence
between good braided Lie algebras and good quantum Lie algebras
(given by $L\to L^+$). Good braided Lie algebras are precisely the
ones which appear as extensions in the context of FODC over
co-quasitriangular Hopf algebras. Not all interesting braided Lie
algebras are `good' in this sense , e.g. the matrix braided Lie
algebras $\calL$ above are not (their $\calL^+$ is not a quantum
Lie algebra but a ``generalised'' one).

In section \ref{Link}, we first recall how quantum Lie algebras
arise in the work of Woronowicz \cite{Woronowicz}, and make clear
what we call the extended (co)tangent spaces of a bicovariant
FODC.  We work with right invariant 1-forms (and left invariant
vector fields), therefore most of our formulas differ from that of
\cite{Woronowicz}. We then prove the main results of this paper,
already mentionned. We give examples of non trivial calculi
arizing from the quantum Lie functor (this mainly concerns the
co-triangular case) and at the far opposite examples of
differential calculi over Hopf algebras which are ``annihilated''
by the quantum Lie functor~: finite groups and quantum groups.
These examples are well-known \cite{BDMS}
\cite{Majid-classif-bicov} \cite{KS} \cite{HS}, but they
illustrate well the fact that $U(\gg_\G)$ is quadratic when $\G$
is inner. Finally, section \ref{sect-Jurco} contains the link
between $B(\calL)$ for such calculi and generalised Lie algebras
along the lines of \cite{LS}.

\section{Preliminaries}\label{sec-Preliminaries}
Throughout, $k$ is a field, vector spaces, algebras, etc, are over
$k$. The flip is written $\tau$ ($\tau(v\otimes w)=w\otimes v)$.
We use Sweedler's notation for coproducts and coactions, omitting
the summation sign, and Einstein's convention for summation over
repeated indices.

\paragraph{Crossed modules, Hopf modules, etc.}

Let $(A,\m,\eta,\D,\e,S)$ be a Hopf algebra.  The Hopf and full
duals of $A$ are $A^\circ\subset A^*$ respectively. The pairing
between $A^*$ and $A$ is written independently $x(a)=\langle
a,x\rangle$, $x\in A^*$, $a\in A$. We let $\ad_L, \ad_R:A\to
A\otimes A$ be the left and right adjoint coaction of $A$ on
itself ($\ad_L(a)=a_{(1)}Sa_{(3)}\otimes a_{(2)}$), and $\Ad_L$,
$\Ad_R$ the left and right adjoint action of $A^\circ$ on itself
($\Ad_Lx(y):= x\rt_{\Ad}y:= x_{(1)}yS(y_{(2)})$). Then the
left coadjoint action of $A^\circ$ on $A$ is $\Ad_L^*x(a)=\langle
a^{(0)},x\rangle \, a^{(1)}$, where $a^{(0)}\otimes
a^{(1)}=\ad_R(a)$. We recall the useful lemma (if $\xi:A\otimes
A\to k$ is some linear map, we define $\xi_1, \xi_2:A\to A^*$ by
$\xi_1(a)(-)=\xi(a,-)$ and $\xi_2(a)(-)=\xi(-,a)$).

\begin{lemma}\label{intertwinner}
Let $\xi:A\otimes A\to k$ be a linear form satisfying
$\m^{op}*\xi=\xi*\m$ and $\im \xi_1\subset A^\circ$, where $*$ is
the convolution product. Then $\xi_1$ intertwines the left adjoint
and coadjoint actions of $A^\circ$ on $A^\circ$ and $A$
respectively, \ie $\Ad_L h\circ \xi_1=\xi_1 \circ \Ad_L^* h$ for
all $h\in A^\circ$. Likewise (if $\im \xi_2\subset A^\circ$),
$\xi_2$ intertwines the right ones.
\end{lemma}
We write $\multi{\calM}{A}{A}{A}{A}$, $\multi{\calC}{A}{}{A}{}$,
$\multi{\calM}{A}{}{}{}$, $\multi{\calM}{}{}{A}{}$ the categories
of Hopf bimodules, left crossed modules, left modules and left
comodules over $A$ respectively. Recall that a Hopf bimodule is a
vector space $\G$ which is both a bimodule and a bicomodule (with
coactions $\D_L:\G\to A\otimes \G$ and $\D_R:\G\to \G\otimes A$),
both coactions commuting with both actions in the natural way. A
left crossed module over $A$ is a vector space $V$ endowed with a
left $A$-action (noted $a\otimes \eta\mapsto a\rt \eta$) and a
left $A$-coaction (noted $\eta\mapsto
\d_L(\eta)=\eta^{(-1)}\otimes \eta^{(0)}$), such that $\d_L(a\rt
\eta)=a_{(1)}\eta^{(-1)}S(a_{(3)})\otimes a_{(2)}\rt \eta^{(0)}$.
The category $\multi{\calC}{A}{}{A}{}$ is a braided (monoidal)
category when the antipode of $A$ is invertible. We shall only
need the braiding on $V\otimes W$ which is given by
$$
v\otimes w\mapsto v^{(-1)}\rt w\otimes v^{(0)},
$$
with inverse $w\otimes v\mapsto v^{(0)}\otimes S^{-1}(v^{(-1)})\rt
w$. As in \cite{Woronowicz}, the categories
$\multi{\calM}{A}{A}{A}{A}$ and $\multi{\calC}{A}{}{A}{}$ are
equivalent. The left crossed module corresponding to $\G$ is
$(\G_R,\rt,\d_L)$ where $\G_R=\{\eta\in\G:\D_R(\eta)=\eta\otimes
1\}$ is the subspace of {\it right} invariants of $\G$, and $a\rt
\eta=a_{(1)}\eta S(a_{(2)})$, $\d_L(\eta)=\D_L(\eta)$. There is a
canonical projection $\pi_R:\G\to \G_R$ given by
$\pi_R(v)=v^{(0)}S(v^{(1)})$, where $v^{(0)}\otimes
v^{(1)}=\D_R(v)$. It satisfies $\pi_R(avb)=a\rt \pi_R(v)\,\e(b)$,
($a,b\in A$, $v\in\G$).
Conversely, $\G$ is recovered from $(\G_R,\rt,\d_L)$ by
$$
\G\simeq\G_R\otimes A\quad {\rm via}\quad
    \left\{ \ba{l}
    v\mapsto \pi_R(v^{(0)})\otimes v^{(1)}
    \\
    \eta\, a\leftarrow \eta\otimes a
    \ea\right.
$$
with tensor product bimodule and bicomodule structure ($\G_R$
is seen as a trivial right module and comodule,
and $A$ is the regular Hopf bimodule)~:
writing $\G=\G_R.A$ as free right $A$-module
(instead of $\G_R\otimes A$),
the extra structures are ($\eta\in \G_R$, $a\in A$)
$$
a\eta=(a_{(1)}\rt \eta). a_{(2)},\quad
\D_L(\eta. a) =\eta^{(-1)}a_{(1)}\otimes (\eta^{(0)}.a_{(2)})
,\quad
\D_R(\eta. a) = (\eta. a_{(1)})\otimes a_{(2)}
$$
A co-quasitriangular structure on a bialgebra $A$ is a linear map
$\rrr:A\otimes A\to k$ which  intertwines the multiplication of
$A$ and its opposite ($\m^{op} *\rrr=\rrr *\m$, $*$ being the
convolution product), and satisfies $\rrr(ab,c)=\rrr(a_{(2)},b)\,
\rrr(a_{(1)},c)$, and $\rrr(a,bc)=\rrr(a,c_{(1)})\,
\rrr(b,c_{(2)})$ for all $a,b,c\in A$. The maps
$\rrr_1,\rrr_2:A\to A^*$ take their values in $A^\circ$ and
satisfy~:
\begin{center}
    \begin{tabular}{l}
    $\rrr_1:A\to A^\circ$ is an algebra/anticoalgebra map,
    \\
    $\rrr_2:A\to A^\circ$ is an antialgebra/coalgebra map.
    \end{tabular}
\end{center}
When $(A,\rrr)$ is co-quasitriangular, the tensor category
$\multi{\calM}{}{}{A}{}$ is braided, the braiding on $V\otimes W$
being
$$
v\otimes w\mapsto \rrr(w^{(-1)},v^{(-1)})\, w^{(0)}\otimes v^{(0)}.
$$
If moreover $A$ has an antipode $S$, then $S^2$ is inner
(hence $S$ is bijective) and the form $\rrr$ is convolution invertible
with inverse $\bar{\rrr}$ such that
$\bar{\rrr}(a,b)= \rrr(S(a),b)$, that is
${\bar{\rrr}}_1=\rrr_1\circ S=S^{-1}\circ \rrr_1$, and
$\bar{\rrr}_2=\rrr_2\circ S^{-1}=S\circ \rrr_2$. Then the braiding on
$V\otimes W$ is invertibe with inverse
$w\otimes v\mapsto \bar{\rrr}(w^{(-1)},v^{(-1)})\, v^{(0)}\otimes w^{(0)}$.

In this article, $(\calV,\otimes)$ is a monoidal category
of the form $\multi{\calM}{A}{}{}{}$, $\multi{\calM}{}{}{A}{}$,
$\multi{\calC}{A}{}{A}{}$ or variants (switching left and right),
$A$ being a bialgebra or Hopf algebra.
Thus, its objects are in particular vector spaces,
and $k$ is the underlying vector space of the unit object.
Recall that if a braiding $\Psi$ exists (eg $A$ is co-quasitriangular
in the case of $\multi{\calM}{}{}{A}{}$), it is a collection
of {\it natural} (iso)morphisms $\Psi_{M,N}:M\otimes N\to N\otimes M$
for all pair $(M,N)$ of objects in $\calV$. The naturality
means that if $f:M\to M'$ and $g:N\to N'$ are morphisms, then
the equality
$\Psi_{M',N'}\circ (f\otimes g)=(g\otimes f)\circ \Psi_{M,N}$ holds.
Also, the structure maps of an algebra, coalgebra, etc are by assumption
{\it morphisms} in $\calV$. Finally, if $A$ and $B$ are algebras in $\calV$,
their tensor product in $\calV$ is $A\uotimes B$, with multiplication
$\m_{A\uotimes B}=(\m_A\otimes \m_B)(\id_A\otimes \Psi_{A,B}\otimes \id_B)$.
The sign $\uotimes$ is to stress the braided structure.
Likewise for coalgebras. Thus a bialgebra in $\calV$ is
both an algebra $(B,\m_B,\eta_B)$ and coalgebra $(B,\D_B,\e_B)$,
$\D_B:B\to B\uotimes B$ and $\e_B:B\to k$ being morphisms of algebras ni
this braided sense.

\paragraph{Quadratic bialgebras.}
Let $(C,\D,\e)$ be a coalgebra in $(\calV,\otimes \Psi)$ as above.
Its tensor algebra
$T(C)=\bigoplus_{n\ge 0} C^{\otimes n}$, with $C^0=k$, is naturally
a bialgebra in $\calV$. The coalgebra structure on each summand
$C^{\otimes n}$
is the (braided) tensor product one. Let $V\subset C\otimes C$
be a subobject and $\langle V\rangle\subset T(C)$ be the 2-sided ideal
generated by $V$.
Clearly, $T(C)/\langle V\rangle$ is a bialgebra in $\calV$ if and only if
$\D_{C\uotimes C}(V)\subset C^{\otimes 2}\otimes V+V\otimes C^{\otimes 2}$
and $\e_{C\uotimes C}(V)=0$.

\begin{lemma}
The bialgebra $T(V)/\langle V\rangle$ never has an antipode.
\end{lemma}

\Proof Assume that there is a map $S':T(C)\to T(C)$ such that the
induced map $S:T(C)/\langle V\rangle\to T(C)/\langle V\rangle$,
$S(a+\langle V\rangle):= S'(a)+\langle V\rangle$, is an
antipode. For all $c\in C\hookrightarrow T(C)$ one should have
$S'(c_{(1)})\otimes c_{(2)}-\e(c)1=\sum_i a_i\otimes v_i\otimes
b_i$ for some $a_i,b_i\in T(C)$, $v_i\in V\subset C\otimes C$. If
$\e(c)=1$, this would mean that $1\in \oplus_{n\ge 1} C^{\otimes
n}$, which is false. \qed
\\

One can always see $V$ as $\im(F)$ for some
morphism $F:C\otimes C\to C\otimes C$ (possibly not of coalgebras).
Then $\langle \im(F)\rangle$ is a bialgebra ideal if and only if
there exist some maps $\Phi,\Phi':C^{\otimes 4}\to C^{\otimes 4}$ such that
\begin{equation}\label{F}
\D_{C\uotimes C}\circ F=
    (\, (F\otimes \id^{\otimes 2})\circ \Phi+
    (\id^{\otimes 2}\otimes F)\circ \Phi'\,)\circ \D_{C\uotimes C},
\quad \e_{C\uotimes C}\circ F=0
\end{equation}
(we use implicitly the fact that $\D_{C\uotimes C}$ is injective
so that any map $C^{\otimes 2}\to C^{\otimes 4}$ is of the form
$\Phi\circ \D_{C\uotimes C}$ for some map $\Phi:C^{\otimes 4}\to
C^{\otimes 4}$). We shall use the following special cases.

\begin{lemma}\label{C-ts}
Let $(C,\D,\e)$ be a coalgebra in $\calV$.
If $\Ups:C\uotimes C\to C\uotimes C$ is a morphism of coalgebras,
then $T(C)/\langle \im(\id-\Ups)\rangle $ is a bialgebra in $\calV$,
without antipode.
\end{lemma}

\Proof Let $F=\id^{\otimes 2}-\Ups:C\uotimes C\to C\uotimes C$.
The hypotheses on $\Ups$ are exactly that (\ref{F}) is satisfied
with $\Phi=\id^{\otimes 4}$ and $\Phi'=\Ups\otimes \id\otimes
\id$, (or $\Phi=\id\otimes \id\otimes \Ups$ and
$\Phi'=\id^{\otimes 4}$). \qed
\\

The following is from \cite{Doi}. Let $(C,\D,\e)$ be a coalgebra
in $\multi{\calM}{k}{}{}{}$ (\ie a usual coalgebra) and
$r:C\otimes C\to k$ some linear map. Define $F_+,F_-:C\otimes C\to
C\otimes C$ by $F_+(a\otimes b)=r(a_{(1)},b_{(1)})\,
a_{(2)}\otimes b_{(2)}$ and $F_-(a\otimes b)=b_{(1)}\otimes
a_{(1)}\, r(a_{(2)}, b_{(2)})$. Then  $F=F_+-F_-$ satisfies
(\ref{F}) with $\Phi=\Phi'=\id^{\otimes 4}$, independently of the
choice of the linear map $r$. Thus $A(C,r):= T(C)/\langle
\im(F_+-F_-)\rangle$ is a bialgebra, generated $C$ with relations
\begin{equation}\label{A-C-r}
r(a_{(1)},b_{(1)})\, a_{(2)}\, b_{(2)}=
b_{(1)}\, a_{(1)}\, r(a_{(2)}, b_{(2)})
\end{equation}

Note that if $r$ is convolution invertible,
$T(C)/\langle \im(F_+-F_-)\rangle=T(C)/\langle \im(\id-\Ups)\rangle$,
where
$\Ups(a\otimes b)=\bar{r}(a_{(1)},b_{(1)})\, b_{(2)}\otimes a_{(2)}\,
    r(a_{(3)}, b_{(3)})$ is a morphism of coalgebras.

\begin{lemma} \label{lemma-FRT}
{\rm \cite{Doi}}
If $r$ is convolution invertible, the following are equivalent~:
\begin{enumerate}
\item
The linear map $r:C\otimes C\to k$ extends (uniquely) to a
co-quasitriangular structure $\rrr$ on $A(C,r)$.
\item
The identity
$r(a_{(1)},b_{(1)})\, r(a_{(2)},c_{(1)})\, r(b_{(2)},c_{(2)})
=
r(b_{(1)},c_{(1)})\, r(a_{(1)},c_{(2)})\, r(a_{(2)},b_{(2)})$
holds for all $a,b,c\in C$.
\item
The map $\Sigma~: C\otimes C\to C \otimes C$,
$\Sigma(a\otimes b)=r(b_{(1)},a_{(1)})\; b_{(2)}\otimes a_{(2)}$,
satisfies the braid relation.
\end{enumerate}
\end{lemma}


\section{Quantum Lie algebras.}
\label{sect-quantum}
\setcounter{equation}{0}

Let $(\calV,\otimes)$ be a -possibly not braided- monoidal category
as in the preliminaries. (Its objects are in
particular $k$-vector spaces, $k$ is the underlying vector
space of the unit object, and by convention $\ga$ will always stand
for a distinguished basis vector of the unit object).
\begin{definition}\label{def-q-Lie}
A left quantum Lie algebra in $\calV$ is a triple $(\gg,\s, [\, ,\,])$
where $\gg$ is an object,
$\s:\gg\otimes \gg\to \gg\otimes \gg$ and $[\, ,\,]:\gg\otimes \gg\to \gg$
are morphisms  satisfying the following axioms
\begin{enumerate}
\item
$\s$ satisfies the braid relation.
\item
Quantum Jacobi identity~:
    $[x,[y,z]]=[[x,y],z]+\sum_i [y_i,[x_i,z]]$ for all $x,y,z\in \gg$,
    where $\sum_i y_i\otimes x_i=\s(x\otimes y)$.
\item
Writing $\s_{12}=(\s\otimes \id)$, $\s_{23}=\id\otimes \s$,
and $C(x\otimes y)=[x,y]$
\begin{eqnarray}
\s\; (\id\otimes C)
    -(C\otimes \id)\; \s_{23}\; \s_{12}
    &=& 0
\\
\s\; (C\otimes \id) -(\id\otimes C) \; \s_{12}\; \s_{23}
 &=&
(C\otimes \id) \; (\id\otimes \s)
- \s\; (\id\otimes C)\;(\s\otimes \id).\label{ket-not-morphism}
\end{eqnarray}
\item
Quantum antisymmetry~: If $\sum_i x_i\otimes y_i\in \ker (\id-\s)$,
    then $\sum_i [x_i,y_i]=0$.
\end{enumerate}
The universal enveloping algebra of $(\gg,\s,[\, ,\,])$ is
\begin{equation}\label{def-U(g)}
U(\gg)=T(\gg)/\langle \im(\id^{\otimes 2}-\s-[\, ,\,])\rangle,
\end{equation}
the tensor algebra of $\gg$ divided by the two-sided ideal
generated by all elements of the form $x\otimes y-\s(x\otimes y)-[x,y]$,
$x,y\in \gg$.
\end{definition}

Woronowicz has shown that these axioms appear naturally in the
context of bicovariant differential calculi over Hopf algebras $A$
(theorems 5.3 and 5.4 in \cite{Woronowicz}). In \cite{Woronowicz},
$\calV$ is $\multi{\calM}{k}{}{}{}$, but this can be made more
precise: it is actually a quantum Lie algebra in the monoidal
category $\multi{\calM}{}{}{A}{}$ (see the comments after
proposition \ref{prop-ext-T}). In the following, ``quantum Lie
algebra'' will mean ``left quantum Lie algebra''.

\begin{lemma}\label{lemma-axiom1}
Axiom 4 of a quantum Lie algebra is the necessary and sufficient condition
for the natural map $\jmath:\gg\hookrightarrow T(\gg)\surject U(\gg)$
to be injective.
\end{lemma}

\Proof
Assume injectivity of $\jmath$
and let $v\in\ker(\id^{\otimes 2}-\s)\subset \gg\otimes \gg$.
Then, in $U(\gg)$, $0=(\id-\s-[\, ,\,])(v)=-[\, ,\,](v)\in\jmath(\gg)$.
By the injectivity of $\jmath$, we get $[\, ,\,](v)=0$.
Conversely, assume antisymmetry and let $z\in \ker \jmath$.
This means that, as an element of $T(\gg)$,
$z=\sum_i u_i\otimes(x_i\otimes y_i-\s(x_i\otimes y_i)-[x_i,y_i])\otimes v_i$
for some $u_i,v_i\in T(\gg)$, $x_i, y_i\in \gg$. On the r.h.s,
terms of degree $\ge 1$ must cancel, \ie one can take $u_i=v_i=1$.
Then, terms of degree two must cancel,
\ie $\sum_i (x_i\otimes y_i-\s(x_i\otimes y_i))=0$. By antisymmetry,
this implies $\sum_i [x_i,y_i]=0$, so $z=0$.
\qed
\\

The other three axioms of a quantum Lie algebra have the following
important interpretation.
Given an object $\gg$ in $\calV$, equipped with morphisms
$\s:\gg\otimes \gg\to \gg\otimes \gg$ and
$[\, ,\,]:\gg\otimes \gg\to \gg$, we define its extension $(\tgg,\ts)$
as follows. We set $\tgg=k\ga\oplus \gg$ and the morphism
$\ts:\tgg\otimes \tgg\to \tgg\otimes \tgg$ is defined by
\begin{equation}\label{ts-vs-s}
\ba{ll}\ts(\ga\otimes z)=z\otimes \ga, \quad \ts(z\otimes
\ga)=\ga\otimes z, &(z\in\tgg)
\\
\ts(x\otimes y)=\s(x\otimes y)+[x,y]\otimes \ga\
    &(x,y\in \gg\hookrightarrow \tgg)
\ea
\end{equation}

\begin{lemma}\label{lemma-braid-for-ts}
The following are equivalent~:
\begin{enumerate}
\item
$\tilde{\s}$ satisfies the braid relation;
\item
the triple $(\gg,\s,[\, ,\, ])$ satisfies axioms $1-3$
of a (left) quantum Lie algebra.
\end{enumerate}
\end{lemma}

\Proof Direct calculation. One would obtain the corresponding axioms
for a right quantum Lie algebra by defining
$\ts(x\otimes y)=\s(x\otimes y)+\ga\otimes [x,y]$ for $x,y\in \gg$.
\qed

\begin{lemma}\label{lemma-Ug-quad}
Let $(\gg,\s,[\,,\,])$ be a quantum Lie algebra.
Let $S_\s(\gg)=T(\gg)/\langle \im(\id-\s)\rangle$
be the {\rm quantum symmetric algebra} of $\gg$ with respect to the braiding
$\s$. Likewise, let
$S_{\ts}(\tgg)=T(\tgg)/\langle \im(\id^{\otimes 2}-\ts)\rangle$.
(i) There are  isomorphisms of algebras
\begin{equation}\label{iso-Ug-Stg}
U(\gg)\simeq S_{\ts}(\tgg)/\langle \ga-1\rangle \quad {\rm and}\quad
S_\s(\gg)\simeq S_{\ts}(\tgg)/\langle \ga\rangle.
\end{equation}
(ii) If there exists a subobject $\calL\subset \tgg$ such that
$\tgg=k \ga\oplus \calL$ and
$\ts(\calL\otimes \calL)\subset \calL\otimes \calL$, then
$S_{\ts}(\tgg)\simeq k[\ga]\otimes S_{\ts_{|\calL}}(\calL)$,
and
\begin{equation}\label{iso-Ug-SL}
U(\gg)\simeq S_{\ts_{|\calL}}(\calL)\simeq  S_\s(\gg).
\end{equation}
\end{lemma}

\Proof
$(i)$ is clear from the definition of $(\tgg,\ts)$, see (\ref{ts-vs-s}).
$(ii)$ If $\calL$ has the given properties, then
$S_{\ts}(\tgg)$ is generated by $\ga$ and $\calL$ with
relations
$x\otimes y=\ts(x\otimes y)$ and  $\ga\otimes x=x\otimes \ga$,
$(x,y\in\calL)$ and the first isomorphism follows (since
$\ts(\calL\otimes \calL)\subset \calL\otimes \calL$ by hypothesis).
Factoring out by $\langle \ga-1\rangle$, we obtain
$U(\gg)\simeq S_{\ts_{|\calL}}(\calL)$ by $(i)$.
Let $\varphi:\tgg\to \gg$
be the projection onto $\gg$ with kernel $k\ga$.
Since by hypothesis $\ga\notin \calL$, the map $\varphi$ induces a
vector space isomorphism
$$
\varphi_{|\calL}:\calL\stackrel{\simeq}{\longrightarrow}\gg,
$$
and satisfies $(\varphi\otimes \varphi)\circ \ts=
\s\circ(\varphi\otimes \varphi)$. Indeed, checking this on
$X\otimes Y$ with either $X$ or $Y$  proportional  to $\ga$ is
immediate from $\varphi(\ga)=0$ and (\ref{ts-vs-s}). Otherwise,
since $\varphi$ is the identity on $\gg$, and since $\ts(x\otimes
y)=\s(x\otimes y)+[x,y]\otimes \ga$ on $\gg\otimes \gg$, we get
$(\varphi\otimes \varphi)\circ \ts_{|\gg\otimes \gg}=\s$ =
$\s\circ(\varphi\otimes \varphi)_{|\gg\otimes \gg}$. Therefore,
$\varphi_{|\calL}$ is a vector space isomorphism which conjugates
the braidings $\ts_{|\calL}$ on $\calL$ and $\s$ on $\gg$, hence
the last isomorphism. \qed
\\

{\it Remarks.}
When $\gg$ is classical ($\s=\tau$ is the flip), the algebra
$S_{\ts}(\tgg)$ already appears in \cite{LeB-S} \cite{LeB-VdB} where
it is written $H(\gg)$ and called the homogenization of $U(\gg)$.
We keep our notation to stress that $H(\gg)=S_{\ts}(\tgg)$
is a quantum symmetric algebra.
Under the hypothesis $(ii)$ of the lemma one obtains an isomorphism
$$
T(\gg)/\langle \im(\id^{\otimes 2}-\s-[\, ,\,])\rangle
\simeq T(\gg)/\langle \im(\id^{\otimes 2}-\s)\rangle.
$$
The quantum Lie bracket has mysteriously disappeared.
However, one has to be careful with this isomorphism
since it is in general {\it not} induced by
the identity isomorphism $T(\gg)\to T(\gg)$,
(unless the quantum Lie bracket is the zero map, in which case
$\calL=\gg$ and $U(\gg)\simeq S_\s(\gg)$ is a tautology).
Obviously, if $\gg$ is classical, $\calL$ exists only in the case
described above ($\gg$ abelian), but more interesting situations
do appear in the ``quantum case'' ($\s\ne \tau$).
\cqfd
\\

We give for completeness a third characterization of quantum Lie
algebras, by a construction due to D. Bernard \cite{Bernard1}. It
shows that one can associate a co-quasitriangular bialgebra to any
quantum Lie algebra $\gg$, in which $U(\gg)$ imbeds as an algebra.
This bialgebra is not, however, what we are after (one would
expect $U(\gg)$ to be a quasitriangular bialgebra, not
co-quasitriangular).

We assume that $\calV$ is the category of vector spaces (hence braided).

Let $C$ be a matrix coalgebra, with comultiplication $\l\mapsto
\l_{(1)}\otimes \l_{(2)}$. Let $\gg$ be the (unique up to
isomorphim) simple left $C$-comodule, with coaction $x\mapsto
x^{(-1)}\otimes x^{(0)}$. It can be viewed as a $C\-k$-bicomodule
for the right coaction $x\mapsto x\otimes \ga$ ($\ga$ is the
grouplike element of the 1-dimensional coalgebra). One obtains a
coalgebra $(\hat{C},\hat{\D},\hat{\e})$, where $\hat{C}:=
C\oplus \gg\oplus k\ga$ and
$$
\hat{\D}(\l)=\l_{(1)}\otimes \l_{(2)},\quad
\hat{\D}(x)=x^{(-1)}\otimes x^{(0)}+x\otimes \ga,\quad
\hat{\D}(\ga)=\ga\otimes \ga.
$$
(Note that $\hat{C}$ is Morita equivalent to the coalgebra of upper
triangular $2\times 2$ matrices.)
Let $\hat{r}:\hat{\gg}\otimes \hat{\gg}\to k$ be a linear map satisfying
\begin{equation}\label{cond-hat-r}
\hat{r}(\ga,-)=\hat{r}(-,\ga)=\hat{\e}(-),\qquad
 \hat{r}(\gg,-)=0
\end{equation}
where ``$-$'' stands for ``anything''. Thus, $\hat{r}$ is uniquely
determined by $r:= \hat{r}_{|C\otimes C}$ and $\o:=
\hat{r}_{|C\otimes \gg}$, which can be arbitrary. We consider the
bialgebra $A(\hat{C},\hat{r})$ as in the preliminaries.
\begin{proposition}
Assume $r$ is convolution invertible.
$\hat{r}$ extends to a co-quasitriangular structure on $A(\hat{\gg},\hat{r})$
if and only if $(\gg,\s,[\,,\,])$ satisfies axioms
1-3 of a (left) quantum Lie algebra,
where
$$
\s(x\otimes y)=r(y^{(-1)},x^{(-1)})\, y^{(0)}\otimes x^{(0)},
\qquad [x,y]=\o(y^{(-1)},x)\, y^{(0)}.
$$
\end{proposition}

\Proof
If $r$ is invertible, so is $\hat{r}$ (its inverse also satisfies
(\ref{cond-hat-r}), and is given by $\l\otimes \mu\mapsto \bar{r}(\l,\mu)$,
$\l\otimes x\mapsto -r(\l_{(1)},x^{(-1)})\, \o(\l_{(2)},x^{(0)})$ for
$\l,\mu\in C$, $x\in\gg$).
By lemma \ref{lemma-FRT},
$\hat{r}$ extends to a co-quasitriangular structure on
$A(\hat{C},\hat{r})$ iff
the map $\hat{\Sigma}:\hat{C}\otimes \hat{C}\to \hat{C}\otimes \hat{C}$,
$\hat{\Sigma}(a\otimes b)=\hat{r}(b_{(1)},a_{(1)})\, b_{(2)}\otimes a_{(2)}$,
satisfies the braid relation. One easily checks that $\hat{\Sigma}$ preserves
the subspace $k\ga\oplus \gg$, where it takes the form (\ref{ts-vs-s}),
with $\s$ and $[\,,\,]$ as stated.
Therefore, if $\hat{\Sigma}$
satisfies the braid relation, by lemma \ref{lemma-braid-for-ts},
$(\gg,\s,[\,,\,])$ must satisfy axioms 1-3 of a quantum Lie algebra.
The converse is long but straightforward. We omit it.
\qed
\\

{\it Remark.}
$\ga$ is grouplike central in $A(\hat{C},\hat{r})$, therefore one can
consider the quotient $A(\hat{C},\hat{r})/\langle \ga-1\rangle$,
which is still co-quasitriangular if $A(\hat{C},\hat{r})$ is.
The other relations in $A(\hat{C},\hat{r})$ are such that the subalgebra
generated by $C$ is isomorphic is $A(C,r)$, the subalgebra generated
by $\ga$ and $\gg$ is $S_{\ts}(\tgg)$, and the crossed relations are given
by ($\l\in C$, $x\in\gg$)~:
\begin{eqnarray*}
\l\, x &=& r(x^{(-1)},\l_{(1)})\, x^{(0)}\, \l_{(2)}
\\
x\, \l&=& r(\l_{(1)},x^{(-1)})\, \l_{(2)}\, x^{(0)} +
\o(\l_{(1)},x)\, \l_{(2)}\, \ga
-x^{(-1)}\, \l_{(1)}\, \o(\l_{(2)}, x^{(0)})
\end{eqnarray*}
If $\o=0$ (\ie $[\,,\,]=0$), combining the two relations above we get
$x\, \l=r_{21}(x^{(-2)},\l_{(1)})$
$r(x^{(-1)},\l_{(2)})$
$x^{(0)}\,\l_{(3)}$.
So, if moreover $r_{21}*r=\e_C\otimes\e_C$, the bialgebra
$A(\hat{C},\hat{r})/\langle \ga-1\rangle$
is just the crossed product of $A(C,r)$ with the quantum symmetric algebra
of its simple comodule. At the far opposite, if $r_{21}*r$ is a non degenerate
bilinear form on $C$, we get $x\, \l=\l \,x=0$ for all $x\in \gg$, $\l\in C$.
When $\o\ne 0$, the terms involving $\o$ are even more unusual.
\cqfd

\paragraph*{``Good'' quantum Lie algebras.}
In this paragraph we investigate bialgebra structures on  $U(\gg$) itself.
To give a sense to this, we assume (until the end of the paper)
that $(\calV,\otimes,\Psi)$ is braided.
We stress that $\gg$ is now equipped with two braidings,
$\s$ and $\Psi_{\gg,\gg}$, which differ in general
(indeed, if $\s=\Psi_{\gg,\gg}$, one  should have
$\s\; (C\otimes \id) =(\id\otimes C) \; \s_{12}\; \s_{23}$
instead of (\ref{ket-not-morphism}) by the naturality of $\Psi$).
The algebra $U(\gg)$ has a filtration
$$
U(\gg)_{(0)}\subset U(\gg)_{(1)}\subset ...\subset U(\gg)_{(n)}\subset ...
$$
induced by the natural $\dZ_{\ge 0}$-grading of $T(\gg)$.
By lemma \ref{lemma-axiom1}
one can identify $U(\gg)_{(1)}$ with $k\, 1\oplus \gg$.
Classically (when $\s=\tau$),
$U(\gg)$ is a Hopf algebra in $\multi{\calM}{k}{}{}{}$
with coalgebra structure $(\D,\e)$ and
antipode $S$ uniquely determined by
$\D(x)=x\otimes 1+1\otimes x$, $\e(x)=0$ and $S(x)=-x$
for all $x\in \gg$.
In particular, each term of the above filtration is a subcoalgebra of
$U(\gg)$, and $\gg\subset \ker\e$.
If we require that it is so in the general case,
a hypothetical bialgebra structure $(\D,\e)$ on $U(\gg)$ is uniquely
determined by a coassociative map $\d: \gg\to \gg\otimes \gg$
(which should be a morphism in the category) such that
\begin{equation}\label{def-D-vs-delta}
\D(x)=x\otimes 1+1\otimes x+\d(x), \quad \e(x)=0,\qquad (x\in\gg).
\end{equation}
We say that $\d$ is a compatible coproduct on $(\gg,\s,[\, ,\,])$
if the above formula defines a coalgebra structure  on $U(\gg)$.
Even if there are some similarities with Lie bialgebras and their
quantization, the situation is different since $\d$ here is
coassociative, and in fact, the choice $\d=0$ is not always
possible:

\begin{lemma}\label{lemma-stand-cop}
Let $(\gg,\s,[\, ,\,])$ be a quantum Lie algebra in $\calV$.
Then $U(\gg)$ is a Hopf algebra in $\calV$ with coalgebra structure
$(\D,\e)$ and antipode $S$ such that
\begin{equation}\label{eq-standard-cop}
\D(x)=x\otimes 1+1\otimes x, \qquad \e(x)=0,\qquad S(x)=-x
\end{equation}
for all $x\in \gg$ if and only if $(1+\Psi_{\gg,\gg})(1-\s)=0$.
\end{lemma}

{\it Proof}. Direct calculation.
\qed
\\

In the general case, we shall restrict ourselves to compatible
coproducts which satisfy a ``nice'' criterion
(a sufficient but not necessary condition).
This criterion is suggested by lemma \ref{C-ts}. A posteriori
``motivations'' for this choice are
given in remark (iii) after theorem \ref{theo-POOM}.

Let $\d:\gg\to \gg\otimes \gg$ be some coassociative morphism.
The extension $\tgg=k\ga\oplus \gg$ of $\gg$
(see (\ref{ts-vs-s})) can be seen
as a coaugmented coalgebra $(\tgg,\td,\te)$ by setting
\begin{equation}\label{tgg-bis}
\td(\ga)=\ga\otimes \ga,\quad
\td(x)=x\otimes \ga+\ga\otimes x +\d(x)\qquad (x\in\gg=\ker\te).
\end{equation}
Note that in fact, $\tgg\simeq U(\gg)_{(1)}$ as a coalgebra.
Recall that $S_{\ts}(\tgg)/\langle \ga-1\rangle\simeq U(\gg)$.
Therefore, since all natural maps $\gg\hookrightarrow U(\gg)$,
$\gg\hookrightarrow \tgg\hookrightarrow S_{\ts}(\tgg)$ are injective,
$\d$ is a compatible coproduct on $\gg$ iff
$S_{\ts}(\tgg)$ is a coalgebra in $\calV$ with coproduct $\tilde{\D}$
such that $\tilde{\D}(X)=\td(X)$ for all $X\in\tgg$.
By lemma \ref{C-ts}, this is ensured if
$\ts:\tgg\uotimes \tgg\to \tgg\uotimes \tgg$ is a morphism of coalgebras,
so~:

\begin{lemma}\label{delta-compatible}
If $\ts: \tgg\uotimes \tgg\to \tgg\uotimes \tgg$
is a morphism of coalgebras, then $\d$ is a compatible coproduct on
$(\gg,\s,[\,,\,])$.
\end{lemma}

Note that what should be a condition on $\d$, the maps $\s$ and
$[\, ,\,]$ being fixed, is finally better seen as a condition on
$\ts$ (\ie $\s$ and $[\, ,\,]$) with respect to a fixed $\d$.
Solving this condition leads to the following definition. -We use
diagrammatic notations as is conventional; compositions of maps
are written from top to bottom, the braiding $\Psi_{V,W}:V\otimes
W\to W\otimes V$ and its inverse $\Psi^{-1}_{W\otimes V} :W\otimes
V\to V\otimes W$ are represented respectively by the symbols~:
\[\epsfbox{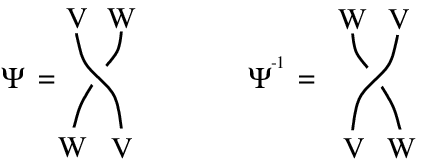}\]

\begin{definition}\label{def-POOM}
A {\rm good quantum Lie algebra} in $\calV$ is a quadruple
$(\gg,\s,[\, , \,],\d)$ where
$\gg$ is an object,
$\s:\gg\otimes \gg\to \gg\otimes\gg$,
$[\, ,\, ]:\gg\otimes \gg\to \gg$ and
$\d:\gg\to \gg\otimes \gg$ morphisms, such that $\d$ is coassociative,
and obeying the  axioms below
\[\epsfbox{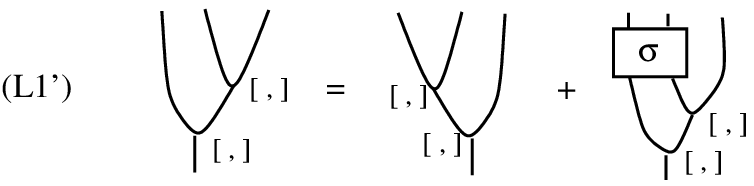}\]
\[\epsfbox{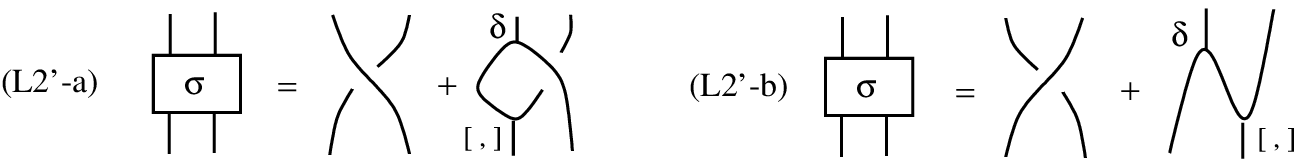}\]
\[\epsfbox{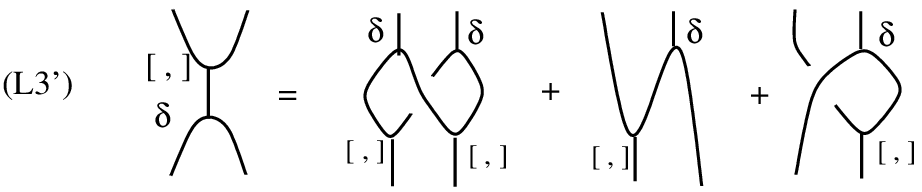}\]
and
\begin{equation}\label{antisymmetry}
\ker(\id^{\otimes 2}-\s)\subset \ker([\, ,\,]).
\end{equation}
\end{definition}

The proof of the following theorem
is given after lemma \ref{good-good}.

\begin{theorem}\label{theo-POOM}
(i) Let $(\gg,\s,[\, ,\,],\d)$ be a good quantum Lie algebra in $\calV$.
Then $(\gg,\s,[\, ,\,])$ is a quantum Lie algebra in $\calV$,
and $\ts:\tilde{\gg}\uotimes \tilde{\gg}\to
    \tilde{\gg}\uotimes \tilde{\gg}$
is a morphism of coalgebras. In particular,
$U(\gg)$ is a bialgebra in $\calV$ with coalgebra structure
given by
$$
\D(x)=x\otimes 1+1\otimes x+\d(x), \qquad \e(x)=0, \qquad (x\in \gg).
$$
(ii) Conversely,
let $(\gg,\s,[\, ,\,])$ be a quantum Lie algebra in $\calV$
equipped with a coassociative morphism
$\d:\gg\to \gg\otimes \gg$. If
$\tilde{\s}:\tilde{\gg}\uotimes \tilde{\gg}\to
    \tilde{\gg}\uotimes \tilde{\gg}$
is a morphism of coalgebras, then $(\gg,\s,[\, ,\,],\d)$
is a good quantum Lie algebra.
\end{theorem}

{\it Remarks.}
(i) The braiding $\s$ of a good quantum Lie algebra
$(\gg,\s,[\,,\,],\d)$ can be expressed
in terms of the maps $[\, ,\,]$, $\d $ and the braiding $\Psi_{\gg,\gg}$,
therefore it is not an essential datum.
Moreover, a good quantum Lie algebra $(\gg,\s,[\,,\,],\d)$
satisfies by hypothesis the axioms
2 and 4 of a quantum Lie algebra,
therefore the first claim of the theorem is that the axioms 1 and 3
are also satisfied, in particular $\s$ satisfies the braid relation.
\\
(ii) If either
$[\, ,\,]=0$ or $\d=0$, one must have
$\s=\Psi_{\gg,\gg}=(\Psi_{\gg,\gg})^{-1}$ by axiom (L2').
Said the other way round, if $\Psi_{\gg,\gg}$ is not symmetric,
neither $[\, ,\,]$ nor $\d$ can be zero
(compare with lemma \ref{lemma-stand-cop}).
\\
(iii) The axioms of a good quantum Lie algebra are satisfied by
the usual Lie algebras ($\s=\Psi_{\gg,\gg}=\tau$, $\d=0$).
Moreover, they almost characterize the standard coproduct in this
case~: if $\s=\Psi_{\gg,\gg}=\tau$, one easily finds that
$(\gg,\s,[\, ,\,],\d)$ is a good quantum Lie algebra if and  only
if $(\gg,[\,,\,])$ is a usual Lie algebra, $\d$ is coassociative,
$\d\gg'=0$ and $\d\gg\subset Z(\gg)\otimes Z(\gg)$, where
$\gg'=[\gg,\gg]$ and $Z(\gg)$ is the center of $\gg$. Therefore,
if $\gg'=\gg$ or if $Z(\gg)=0$, one must have $\d=0$.
\\
(iv) There are compatible coproducts which are not ``good''~: for
instance, take $\gg=e(1,1)$ the usual Lie algebra with basis
$e_0,e_+,e_-$ such that $[e_0,e_\pm]=\pm\, e_\pm$, $[e_+,e_-]=0$
(if $k=\dR$, $\gg$ is the Lie algebra of the pseudo-euclidean
plane). Then $U(\gg)$ is a  Hopf algebra in
$\multi{\calM}{k}{}{}{}$ for all coproducts of the form
(\ref{def-D-vs-delta}) with $\d e_\pm=0$, $\d e_0=\l\,(e_+\otimes
e_- -e_-\otimes e_+)$ for all $\l\in k$; in fact if $\l\ne 0$ it
can be rescaled, eg to $\l=1$, by rescaling $e_+$ or $e_-$. The
antipode is given by $S(x)=-x$ for $x\in\gg$, independently of
$\l$. When $\l\ne 0$, $(\gg,[\,,\,],\d)$ is not good
 since $\d\ne 0$ but $Z(\gg)=0$.
\cqfd


\section{Braided Lie algebras and Lie coalgebras}
    \label{sect-braided}
\setcounter{equation}{0}

The definition of a good quantum Lie algebra is already coming close
to that of a braided Lie algebra \cite{Majid-braided-Lie}.
In this section we recall their definition and main properties,
and discuss the connection with good quantum Lie algebras.
In particular, we shall see that the axioms of a good quantum Lie algebra
can take a much simpler form when expressed in terms of the braided Lie algebra
it corresponds to.

\subsection{$L$ and $B(L)$}\label{subsect-braided-Lie-algebras}

\begin{definition} \label{def-b-Lie}
{\rm \cite{Majid-braided-Lie}}
A (left) braided Lie algebra in $\calV$ is a coalgebra
$(L,\D,\e)$ in the category, equipped with
a morphism in $\calV$ (the braided Lie bracket)
$[\, ,\, ]:L \underline{\otimes} L\to L$ satisfying
the axioms pictured below
\[\epsfbox{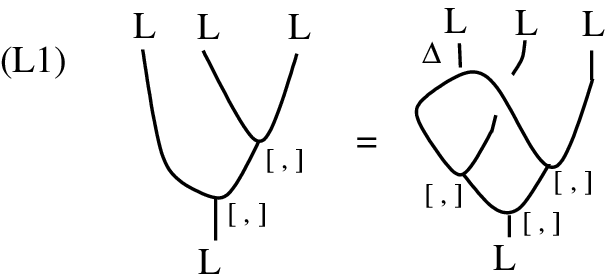}\qquad\qquad\epsfbox{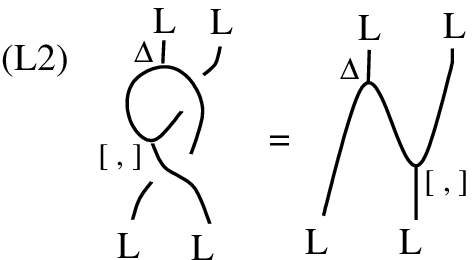}\]
\[\epsfbox{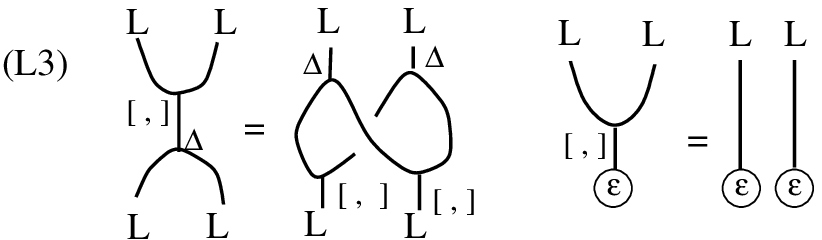}\]
Axiom {\rm (L1)} is called the left
braided Jacobi identity,  {\rm (L2)} weak braided cocommutativity,
and {\rm (L3)} states that $[\, ,\,]:L\uotimes L\to L$ is a morphism
of coalgebras.

A braided Lie subalgebra of $L$ is subcoalgebra $M$ such that
    $[M,M]\subset M$.
A morphism of braided Lie algebras in $\calV$ is a morphism
of coalgebras $\phi:L_1\to L_2$ such that
$[\, ,\,]_2\circ (\phi\otimes\phi)=\phi\circ [\, ,\,]_1$.
\end{definition}

{\it Remarks.} $(i)$ If the braiding on $L$ is symmetric and if
$(L,\D,\e)$ is cocommutative ($\Psi\circ \D=\D$), axiom (L2) is
automatically satisfied, independently of $[\, ,\,]$.

$(ii)$ By a braided Lie algebra we shall mean a left one. Axioms
for a right braided Lie algebra are obtained from that of a left
one by applying a symmetry along the medium vertical axis of each
diagram in the definition while keeping the same order of
crossing\cite{Majid-lec} (so that the diagrams of axiom (L3)
remain unchanged); see for instance \cite{Wambst}. It is observed
in \cite{Wambst} that if $(L,\D,\e,[\, ,\,])$ is a left braided
Lie algebra, then $(L,\Psi^{-1}\circ \D, \e, [\, ,\,]\circ \Psi)$
is a right one.

$(iii)$
The naturality of the braiding with respect to the morphisms
$\D$, $\e$ and
$[\, , \,]$ here means that~:
\begin{eqnarray*}
(\id\circ \D)\circ \Psi=
\Psi_{12}\circ \Psi_{23}\circ \D
    &,&
(\D\circ \id)\circ \Psi=
    \Psi_{23}\circ \Psi_{12}\circ \D,
\\
(\e\otimes  \id)\circ \Psi =
    \id\otimes \e
    &,&
(\id\otimes \e)\circ \Psi =
    \e\otimes \id,
\\
\Psi\circ ([\, ,\,]\otimes \id)=
    (\id\otimes [\, ,\, ])\circ \Psi_{12}\circ \Psi_{23}
    &,&
\Psi\circ (\id\otimes [\, ,\,])=
    ([\, ,\, ]\otimes  \id)\circ \Psi_{23}\circ \Psi_{12}.
\end{eqnarray*}
This identities should be added explicitly to the axioms if one
forgets about a background category and consider a vector space
$L$ equipped with maps $(\Psi,\D,\e,[\, ,\,])$. For instance, in
dimension $1$, there is only one isomorphism class of braided Lie
algebra. Indeed, let $L=k\, \ga$ with $\ga$ grouplike (after
scaling). Then there exists scalars $\l,q\in k$ such that
$\Psi(\ga\otimes\ga)=q\ga\otimes \ga$ and $[\ga,\ga]=\l\, \ga$.
Counit and naturality constraints then force $\l=q=1$. \cqfd

\begin{lemma}\label{lemma-dir-sum}
(i) Let $(C,\D,\e)$ be a coalgebra in $\calV$.
The map $[\, ,\,]_\triv:C\otimes C\to C$,
$[x,y]_{{\rm triv}}=\e(x)\, y$ is a braided Lie bracket on $C$
if and only if $(\Psi_{C,C})^2=\id_{C\otimes C}$.
\\
(ii) Let $(L_i,\D_i,\e_i,[\, ,\,]_i)$, $i=1,2$, be two braided Lie
algebras in $\calV$. The direct sum coalgebra $L=L_1\oplus L_2$
equipped with the map $[\, ,\,]_L:L\otimes L\to L$
$$
[x_1\oplus x_2,y_1\oplus y_2]_L
=
[x_1,y_1]_1+\e_2(x_2)\, y_1\;\oplus\;
[x_2,y_2]_2+\e_1(x_1)\, y_2
$$
is a braided Lie algebra if and only if
$\Psi_{L_i, L_j}\circ \Psi_{L_j, L_i} =\id_{L_i\otimes L_i}$ when
$i\ne j$.
In this case we call $L$ the direct sum of $L_1$ and $L_2$.
\end{lemma}

\Proof
$(i)$ Using only the counit axiom of a coalgebra and the naturality of $\Psi$,
one easily checks that this bracket always satisfies axioms (L1)
and (L3), but satisfies (L2) iff $\Psi_{C,C}=(\Psi_{C,C})^{-1}$.
Note that the trivial bracket
on a right braided Lie algebra in case $\Psi^2=\id$ would be
$[x,y]_{{\rm triv}}=x\, \e(y)$.
$(ii)$ The reasons are the same as in $(i)$.
\qed

\begin{theorem}\label{prop-Ups-vs-ket}
Let $(L,\D,\e)$ be a coalgebra in $\calV$. There is a one-to-one
correspondence between
\begin{enumerate}
\item
morphisms of coalgebras $[\, ,\,]:L\uotimes L\to L$ such
that $(L,\D,\e,[\,,\,])$ is a left braided Lie algebra.
\item
morphisms of coalgebras $\Ups:L\uotimes L\to L\uotimes L$ such that
\begin{enumerate}
    \item
    $\Ups$ satisfies the braid relation,
    \item
    $\Ups(\ker\e\otimes \ker\e)\subseteq \ker\e\otimes L$,
    \item
    the following equalities hold (in the box, $Y$ means $\Ups$)~:
\ceqn{Ups-vs-brack}{\epsfbox{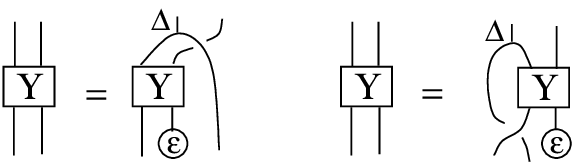}}
\end{enumerate}
\end{enumerate}
It is given by
$\Ups:=([\, ,\,]\otimes \id)\circ
    (\id\otimes \Psi)\circ (\D\otimes \id)$,
and
$[\, ,\,]:=(\id\otimes \e)\circ \Ups$.
In diagrammatic form~:
\ceqn{def-PsiD}{\epsfbox{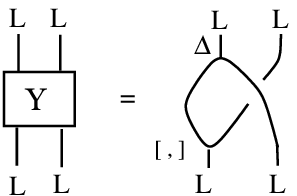}}
We call $\Ups$ the canonical braiding of $(L,\D,\e,[\,,\,])$.
\end{theorem}

\Proof
Let $(L,\D,\e,[\, ,\,])$ be a braided Lie algebra,
and define $\Ups$ as in (\ref{def-PsiD}).
The fact that it satisfies the braid relation is proved in
\cite{Majid-sol-YBE}\cite{Wambst}.
By definition of $\Ups$ and by the counit axioms, one has
\begin{equation}\label{eps-vs-PsiD}
(\e\otimes \id)\circ \Ups(x\otimes y)=\e(y)\, x,
\quad
(\id\otimes\e)\circ \Ups(x\otimes y)=[x,y].
\end{equation}
Therefore $\Ups$ satisfies $(b)$ and the counit part of the fact that
it is a morphism of coalgebras. The coproduct part, \ie the equality
$(\Ups\otimes \Ups)\circ \D_{L\uotimes L} =\D_{L\uotimes L}\circ \Ups$,
is checked as~:
\[\epsfbox{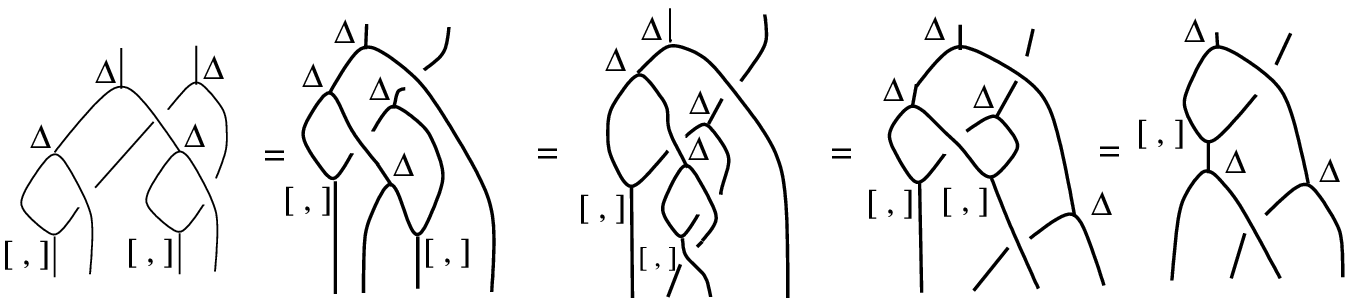}\]
The first and third equalities use only the naturality and coassociativity
axioms, the second axiom (L2) and the fourth axiom (L1).
Finally, $\Ups$ satisfies the equalities (\ref{Ups-vs-brack})
since, in view of (\ref{eps-vs-PsiD}),
the left one is nothing but its definition,
and the right one is an equivalent form of axiom
(L2) (multiplied by $\Psi^{-1}$ on the left,
\ie on the bottom).

Conversely, define $[\, ,\,]:= (\id\otimes \e)\circ
\Ups:L\uotimes L\to L$. Obviously, $[\, ,\,]$ is a morphism of
coalgebras (axiom (L3)), as composition of morphism of coalgebras.
Next, axiom (L1) is satisfied since~:
\ceqn{Ups-Jac}{\epsfbox{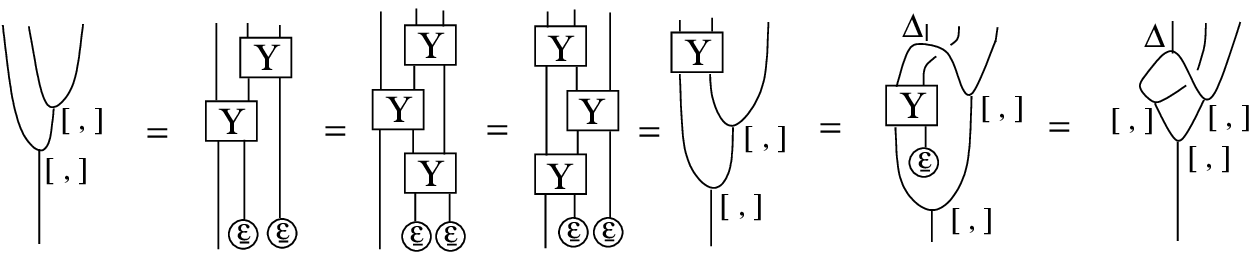}} The first equality is by
definition of $[\, ,\,]=(\id\otimes \e)\circ \Ups$, the second
uses the fact that $\Ups$ is a morphism of coalgebras, the third
is the braid relation for $\Ups$, the fourth is again the
definition of $[\, ,\,]$, the fifth uses the left equality in
(\ref{Ups-vs-brack}), and the sixth is again the definition of
$[\, ,\,]$. Finally, the left equality in (\ref{Ups-vs-brack})
means that one can reexpress $\Ups$ in terms of $[\, ,\,]:=
(\id\otimes \e)\circ \Ups$ as $\Ups=([\, ,\,]\otimes \id)\circ
(\id\otimes \Psi)\circ (\D\otimes \id)$. Therefore the second
diagram of (\ref{Ups-vs-brack}) is again nothing but axiom (L2)
for $(L,\D,\e,[\,,\,])$. \qed
\\

By definition, the braided enveloping algebra of
$(L,\D,\e,[\, ,\,])$ is the symmetric algebra of $L$ with respect
to $\Ups$~:
\begin{equation}\label{eq-def-BL}
 B(L):= S_\Ups(L)=T(L)/\langle \im(\id-\Ups)\rangle.
\end{equation}
(Its definition is motivated by the braided Jacobi identity
which can be expressed as the equality between the first and fifth
diagrams in (\ref{Ups-Jac})).
Since $\Ups$ is a morphism of coalgebras,
by lemma \ref{C-ts} the maps
$\D:L\to L\uotimes L$ and
$\e:L\to k$ extend uniquely to algebra morphisms
$B(L)\to B(L)\uotimes B(L)$ and $B(L)\to k$ respectively, \ie~:

\begin{corollary} \label{cor-BL-bialgebra}
{\rm\cite{Majid-braided-Lie}}.
$B(L)$ is a bialgebra in $\calV$.

\end{corollary}

{\it Remark.}
 $B(L)$ is quadratic, and ``$\Ups$-commutative'',
but it lives in $\calV$ where the braiding is $\Psi$. Note that
$\Ups$ need not be invertible (although we do not know any example
where it is not), and that there is no way to express the original
braiding $\Psi_{L,L}$ in terms of $\Ups$, $\D$, $\e$  and $[\,
,\,]$. \cqfd

\begin{proposition}\label{prop-BL-funct}
(i) The correspondence $L\to B(L)$ is an exact covariant functor.
\\
(ii) Let $L=L_1\oplus L_2$ be the direct sum of two braided
Lie algebras. Then $B(L)\simeq B(L_1)\uotimes B(L_2)$ is the tensor product
of $B(L_1)$ and $B(L_2)$ in the category.
\end{proposition}

\Proof $(i)$ A morphism of objects $L\to M$ induces a morphism of
algebras $T(L)\to T(M)$. By definition, a morphism of braided Lie
algebras $f:L\to M$ intertwines all structure maps, in particular
$f$ is a morphism of coalgebras and $(f\otimes f)\circ
\Ups_L=\Ups_M\circ (f\otimes f)$. Therefore $f$ induces a
bialgebra morphism $\tilde{f}:T(L)/\langle \im(\id^{\otimes
2}-\Ups_L)\rangle\to T(M)/\langle \im(\id^{\otimes
2}-\Ups_M)\rangle$. Clearly $\tilde{f}$ is injective (or
surjective) if and only if $f$ is. $(ii)$ By $(i)$, there are
bialgebras embeddings $B(L_1)\hookrightarrow B(L)\hookleftarrow
B(L_2)$. The claim follows from the observation that (by
definition of a direct sum, see lemma \ref{lemma-dir-sum}), for
$x\in L_i$, $y\in L_j$, and $i\ne j$, one has $\Ups(x\otimes
y)=\Psi(x\otimes y)$. \qed

\subsection{Good braided Lie algebras.}
We shall say that $(L,\D,\e,[\, ,\,])$ is {\it unital} if there
exists a morphism of braided Lie algebras $\eta:k\to L$ with
\begin{equation}\label{def-unital}
[\, ,\,]\circ ( \eta\otimes \id_L)= \id_L,\quad
{\rm and}\quad [\,,\,]\circ (\id_L\otimes \eta)=\eta\circ \e,
\end{equation}
or equivalently (by theorem \ref{prop-Ups-vs-ket})
\begin{equation}\label{def-unital-bis}
\Ups\circ ( \eta\otimes \id_L)= \id_L\otimes \eta,\quad
{\rm and }\quad \Ups\circ (\id_L\otimes \eta)=\eta\otimes \id_L.
\end{equation}
In terms of $\ga=\eta(1)$, $(L,\D,\e,[\, ,\,],\eta)$ is unital
if the span of $\ga$ is isomorphic to the trivial object (which
implies that it is ``bosonic''~: $\Psi(\ga\otimes z)=z\otimes
\ga$ and $\Psi(z\otimes \ga)=\ga\otimes z$ for all $z\in L$), is
grouplike ($\eta$ is a morphism of coalgebras), and $[\ga,z]=z$,
$[z,\ga]=\eps(z)\gamma$ for all $z\in L$. The last two equalities are
equivalent to $\Ups(\ga\otimes z)=z\otimes \ga$ and
$\Ups(z\otimes \ga)=\ga\otimes z$ for all $z\in L$. In particular
$\ga$ is a central grouplike in $B(L)$.
We stress that not all braided Lie
algebras are unital and the morphism $\eta:k\to L$ or the element
$\ga$, if it exists, is not unique in general.
If $(L,\D,\e,[\, ,\,],\eta)$ is unital, we define
\begin{equation}\label{g}
\gg:= \ker\e,\quad {\rm and}\quad \ga:= \eta(1),
\end{equation}
so that there is a distinguished
decomposition $L=k\ga\oplus \gg$.
By the counit axioms and by (\ref{eps-vs-PsiD}), there exists unique morphisms
$\s:\gg\otimes\gg\to \gg\otimes\gg$
and $\d:\gg\to \gg\otimes\gg$ such that
\begin{equation}\label{def-sigma}
\Ups(x\otimes y)=\s(x\otimes y)+[x,y]\otimes \ga,
\quad
\D(x)=x\otimes \ga+\ga\otimes x+\d(x).
\end{equation}
The axioms of the  braided Lie algebra $L$ and the structure of
$B(L)$, in the unital case, can be given in terms of
$(\gg,\sigma,[\ ,\ ],\d)$ as was done in  \cite[Fig
11]{Majid-braided-Lie}. Exactly as $\Ups$ is expressible in terms
of $\Psi_{L,L}$, $\D$ and $[\,,\,]$, the map $\s$ can be expressed
in terms of $\Psi_{\gg,\gg}$, $\d$ and $[\,,\,]_{|\gg\otimes
\gg}$. The braided Jacobi identity -axiom (L1)- for $L$~:
\begin{equation}
\label{Jac-L} [X,[Y,Z]]=[\ ,[\ ,\ ]](\Ups(X\otimes Y)\otimes Z),\quad
\forall X,Y,Z\in L
\end{equation} becomes when restricted to $\gg$~:
\begin{equation}\label{jacobi-sigma} [x,[y,z]]=
[\ ,[\ ,\ ]](\sigma(x\otimes y)\otimes z)+[[x,y],z],\quad\forall
x,y,z\in \gg.
\end{equation}
Similarly the braid relation for $\Ups$ implies the braid relation
for $\sigma$. Also from the first of (\ref{def-sigma}) the
relations of $B(L)$ become $\gamma$ central and
\begin{equation}\label{rel-sigma}
xy-\cdot\circ\sigma(x\otimes y)=[x,y]\ga,\quad \forall x,y\in\ga.
\end{equation}

Note now that (\ref{jacobi-sigma}) is axiom 2 of a quantum Lie
algebra for $(\gg,\sigma,[\ ,\ ])$. In fact if we identify
$L=\tgg$ then lemma~\ref{lemma-braid-for-ts} with $\Ups=\ts$
ensures that $\gg$ satisfies axioms 1 and 3 of a quantum Lie
algebra. The only item missing is an antisymmetry property for the
Lie bracket, which can be added at hand~:

\begin{definition}\label{def-good}
A {\rm good} braided Lie algebra in $\calV$ is a unital braided Lie algebra
$(L,\D,\e,[\, ,\,],\eta)$ such that
$\ker(\id_{\gg\otimes\gg}-\s)\subset \ker([\, ,\,]_{\gg\otimes\gg})$,
where $(\gg,\s,[\,,\,],\d)$ is defined as above.
\end{definition}
{}From the above discussion, if $(L,\D,\e,[\,,\,],\eta)$ is a good
braided Lie algebra, then $(\gg,\s,[\,,\,])$ is  a quantum Lie
algebra. Also it is evident from (\ref{rel-sigma}) that a
sufficient condition for a unital braided Lie algebra $L$ to be
good is that $\gamma$ is not a zero divisor in $B(L)$.
Unsurprisingly from the above discussion, and taking onto account
the coproduct one has~:

\begin{lemma}\label{good-good}
$(\gg,\s,[\,,\,], \d)$ is a good quantum Lie algebra
if and only if its extension $(\tgg,\ts,\td,\te)$ is a
good braided Lie algebra, with braided Lie bracket
${[}\,,\,\tket=(\id\otimes \te)\circ \ts$.
\end{lemma}

\Proof
Coassociativity of $\d$ and $\td$ are obviously equivalent,
and the antisymmetry axiom is postulated in both definitions.
The reader will easily check that axiom (L1), (L2) and
(L3) for $\tgg$ are equivalent to axioms (L1'), (L2'-b) and (L3')
for $\gg$, while (L2'-a) corresponds to the definition of the
canonical braiding $\Ups\equiv\ts$ for $\tgg$. We omit the details.
(Note that the braided Lie bracket on $\tgg$
is given by ${[}\ga,\ga\tket=\ga$, ${[}\ga,x\tket=x$,
${[}x,\ga\tket=0$, $x\in\gg$, and
${[}\,,\,\tket_{|\gg\otimes \gg}=[\,,\,]$).
\qed
\\

{\it Proof of theorem \ref{theo-POOM}.}
$(i)$ If $(\gg,\s,[\, ,\,],\d)$ is a good quantum Lie algebra,
then its extension
$(\tgg,\ts,\td,\te)$ is a good braided Lie algebra;
by the discussion before definition \ref{def-good}
this implies that $(\gg,\s,[\, ,\,])$ is a quantum Lie algebra.
Moreover, $\Ups_{\tgg}=\ts:\tgg\uotimes \tgg\to \tgg\uotimes \tgg$ is a
morphism of coalgebras by theorem \ref{prop-Ups-vs-ket}
so that $\d$ is a compatible
coproduct on $\gg$ by lemma  \ref{delta-compatible}.
\\
$(ii)$ Let $(\gg,\s,[\,,\,],\d)$ be a quantum Lie algebra with
$\d:\gg\to \gg\otimes \gg$ coassociative, and $(\tgg,\ts,\td,\te)$
be the extension of $\gg$. $\ts$ is a morphism of coalgebras iff
$\td_{\tgg\otimes \tgg}\circ\, \ts(X\otimes Y)= (\ts\otimes
\ts)\circ\td_{\tgg\otimes \tgg}(X\otimes Y)$, for all $X,Y\in\gg$,
where $\td_{\tgg\otimes \tgg}= (\id\otimes \Psi\otimes
\id)\circ(\td\otimes \td)$. One checks that this is trivially
satisfied if either $X$ or $Y$ is proportional to $\ga$. For
$X,Y\in \gg$, we obtain an equation in $\tgg^{\otimes 4}$. Using
the decomposition $\tgg=k\ga\oplus \gg$, this leads to $2^4$
equations, eight of which are trivially satisfied. In the
remaining eight, three are the identities corresponding to axioms
(L2') and (L3').
 Therefore, $(\gg,\s,[\,,\,],\d)$
must be a good quantum Lie algebra. The remaining five identities are
automatically satisfied by $(i)$.
\qed
\\

Any braided Lie algebra $\calL$
can be imbedded in a unital one $(L,\eta_L)$~: take the direct sum
of braided Lie algebras $L=k\ga\oplus \calL$
(in the sense of lemma \ref{lemma-dir-sum}),
with  unit $\eta_L(1)=\ga$;
here $k\ga$ is the unique one dimensional braided Lie algebra.
We call $L$ the {\it trivial extension}
of $\calL$. Not all unital braided Lie algebras are
trivial extensions.

\begin{proposition}\label{propo-triv-ext}
If $L$ is the trivial extension of some braided Lie algebra $\calL$,
then $L$ is good, and $(\gg=\ker\e_L,\s,[\,,\,],\d)$ is a good quantum
Lie algebra (see {\rm (\ref{g}) (\ref{def-sigma})}).
Moreover, $\d:\gg\to\gg\otimes \gg$ is injective,
and there are bialgebra isomorphisms $B(L)\simeq k[\ga]\otimes B(\calL)$ and
$U(\gg)\simeq B(\calL)$.
\end{proposition}

\Proof The bialgebra isomorphism $B(L)\simeq k[\ga]\otimes
B(\calL)$ (with $\ga$ grouplike, not primitive) is by proposition
\ref{prop-BL-funct}. We note that $\ga$ is not a zero divisor in
$B(L)$  so that $L$ is good. For the other statements, we apply
lemma \ref{lemma-Ug-quad} to $(\tgg,\td,\te):= (L,\D_L,\e_L)$ and
$\ts:= \Ups_L$. Clearly, $(\Ups_L)_{|\calL}=\Ups_\calL$. The
projection $\varphi:L=\tgg\to \gg$ with kernel $ k\ga$ (now given
by $\varphi(X)=X-\e(X)\ga$, $X\in L$) restricts to an isomorphism
$\varphi_{|\calL}~:\calL\stackrel{\simeq}{\longrightarrow}\gg$ and
satisfies $(\varphi\otimes \varphi)\circ \Ups_L=
\s\circ(\varphi\otimes \varphi)$ as already known,  and also
$(\varphi\otimes \varphi)\circ \D_L=\d\circ\varphi$ as is easily
checked. Therefore, when restricted to $\calL$, the previous
equation tells that $\d$ is injective (since
${\D_L}_{|\calL}=\D_\calL$ is), and that the algebra isomorphism
$U(\gg)\simeq S_{\Ups_{\calL}}(\calL)=:B(\calL)$ of lemma
\ref{lemma-Ug-quad} is also a bialgebra isomorphism in the present
case. \qed
\\

{\it Remark.} We see that {\it any} braided Lie algebra $\calL$ can help to
construct a good quantum Lie algebra $(\gg,\s,[\,,\,],\d)$
{\it of the same dimension}, with $\d$ injective~:
take $\gg=\ker\e_L$ where $L=k\ga\oplus \calL$. Clearly, not all good
quantum Lie algebras are of this form (for instance, the usual Lie algebras).
This large class of examples shows that $U(\gg)$ does not always have an
antipode, since in this case $U(\gg)\simeq B(\calL)$ cannot have an antipode.
\cqfd

\subsection{Split braided Lie algebras.}

The notion of a unital braided Lie algebra can generalized as
follows. Indeed, while we are interested in unital extensions of
braided Lie algebras $\calL$, these $\calL$ themselves are not
typically unital. Yet deformation examples should be close to
unital ones since the classical model for the entire theory in
\cite{Majid-braided-Lie} is the example $\calL=k1\oplus$ a
classical Lie algebra.

Thus, we say a braided Lie $\calL$ is {\em split} if there is a
morphism $c:k\to \calL$ in the braided category, or in concrete
terms a distinguished element $c\in \calL$ with span isomorphic to
the trivial object, such that
\begin{equation}\label{split}\e(c)=1,\quad [x,c]=\eps(x) c,\quad
\forall x\in \calL.\end{equation} Here $kc$ the trivial object
implies $\Psi(c\otimes x)=x\otimes c$ and $\Psi(x\otimes c)=c\otimes x$
for all $x\in \calL$, while the second of (\ref{split}) is
equivalent (by the counit axioms) to
\begin{equation}\label{Ups-cx} \Ups(x\otimes c)=c\otimes x,\quad
\forall x\in \calL.\end{equation} Being split is significantly
weaker than the unital case. However, (\ref{Ups-cx}) still
ensures that $c$ is still central in $B(\calL)$. We set
$\calL^+=\ker\e_\calL$. Then, by again the counit axioms, there
exist uniquely determined maps $\o:(\calL^+)^{\otimes 2}\to
(\calL^+)^{\otimes 2}$, $\rho:\calL^+\to (\calL^+)^{\otimes
2}$ and $\Theta:\calL^+\to\calL^+$ such that, for all
$x,y\in\calL^+$,
\begin{eqnarray}
\Ups(x\otimes y)&=&\o(x\otimes y) +[x,y]\otimes c,\label{gen-xy}
\\
    \Ups(c\otimes x) &=& \Theta(x) \otimes c+\rho(x);\label{gen-cx}
\\ \Theta(x)&=& [c,x].
\end{eqnarray}
As in the unital case, one may write the axioms of a braided Lie
algebra in the split case in terms of $(\calL^+,\omega,[\ ,\
],\rho,\Theta)$. For example, the braided Jacobi identity
(\ref{Jac-L}) on $\calL$ gives
\begin{equation}\label{Jac-split}[x,[y,z]]=
[\ ,[\ ,\ ]](\omega(x\otimes y)\otimes
z)+[[x,y],\Theta(z)],\quad\forall x,y,z\in\calL^+\end{equation}
generalising the unital case. Note also that, because of (\ref{gen-cx}),
lemma \ref{lemma-braid-for-ts} cannot be applied in general and therefore
$\omega$ (in the
place of $\sigma$) needs not obey the braid relation. Moreover,
$B(\calL)$ in these terms is generated by $c,\calL^+$ with $c$
central and the relations
\begin{equation}\label{rel-split} xy-\cdot\circ\omega(x\otimes
y)=[x,y]c,\quad cx-\Theta(x)c=\cdot\circ\rho x,\quad\forall
x,y\in\calL^+.\end{equation} Also note that if $c$ is not a zero
divisor in $B(\calL)$ then clearly
\begin{equation}\label{good-split}
\ker (\id-\omega)\subset \ker ([\,,\,])
\end{equation}
just as for good quantum Lie algebras above.

Let us assume now that $\calL^+$ is a simple object. Then $\Theta$
acts as a multiple $\lambda$ of the identity and we define the
reduce enveloping algebra associated to the split braided Lie
algebra to be:
\[ B_{red}(\calL^+)=B(\calL)/\langle c-\lambda\rangle.\]
It is the tensor algebra $T(\calL^+)$ modulo the ideal generated
by the relations for all $x,y\in\calL^+$:
\begin{eqnarray}\label{def-gen}
xy-\cdot\circ\o(x\otimes y)&=&\l\, [x,y] \\
\l(1-\l )\, x &=&\cdot\circ\rho(x).\label{problem}
\end{eqnarray}

Finally, we suppose that $\lambda\ne0$ and $c$ is not a zero
divisor of $B(\calL)$. If we define
\[ A=\lambda^{-1}(\id-\omega)\] then clearly  (\ref{good-split})
and (\ref{Jac-split}) appear as
\begin{equation}\label{gen-Jacobi}
\ker A\subset\ker([\ ,\ ]),\quad [\ ,[\ ,\ ]](A(x\otimes
y),z)=[[x,y],z],\quad\forall x,y,z\in\calL^+
\end{equation}
which have been proposed as the axioms of a `generalized Lie
algebra' $(\calL^+,A,[\,,\,])$ in\cite{LS}. Here
(\ref{gen-Jacobi}) is called a generalized Jacobi identity and $A$
is called a generalized antisymmetrizer (although $A$ is not
required to satisfy any further axioms in this regard). Similarly,
according to the definition of \cite{LS}, the universal enveloping
algebra of $(\calL^+,A,[\,,\,])$ is
$U_{LS}(\calL^+)=T(\calL^+)/\langle \im(A-[\,,\,])\rangle$, i.e.
generated by $\calL^+$ with the relation (\ref{def-gen}). We see
that split braided Lie algebras with simple $\calL^+$ have this
general structure but with $\omega$, $[\ ,\ ]$ and an additional
map $\rho$ obeying several more axioms inherited from the
braided Lie algebra structure. We also see that the natural
`enveloping algebra' generated by $\calL^+$ in this case, namely
$B_{red}(\calL^+)$ has potentially an additional relation
(\ref{problem}). On the other hand, $B_{red}(\calL^+)$ comes as a
quotient of a quadratic algebra $B(\calL)$ giving its
homogenisation and forming a bialgebra (in a braided category),
both of them desirable features.

\subsection{The adjoint action.}
Let $\Rep(L)$ be the category of representations of $L$~: its
objects are pairs $(V,\a)$ where $\a:L\otimes V\to V$ is a
morphism in $\calV$ satisfying axiom (R1) pictured in
(\ref{fig-R}); $\a$ is called the action of $L$ on $V$ and we
write $x\rt_\a v=\a(x\otimes v)$. A morphism of representations
(intertwiner) is a morphism $f:V\to W$ in $\calV$ satisfying
$f\circ\a_V=\a_W\circ(\id\otimes f)$. Clearly, $\Rep(L)$ is the
same as $\multi{\calM}{B(L)}{}{}{}$. It is a monoidal category
with tensor product $(V,\a_V)\otimes (W,\a_W)=(V\otimes
W,\a_{V\otimes W})$ where $\a_{V\otimes W}= (\a_V\otimes
\a_W)(\id_L\otimes \Psi_{L,V}\otimes \id_W) (\D_L\otimes
\id_{V\otimes W})$, and unit object $k$ (with action afforded by
the counit $\e$). \ceqn{fig-R}{\epsfbox{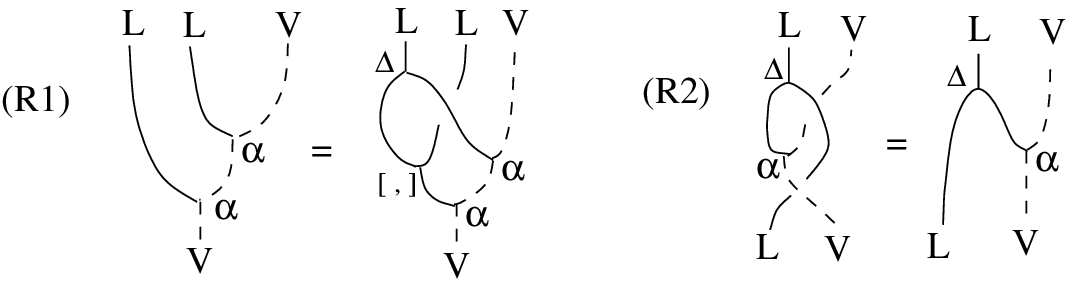}} We also let
$\Rep(L)'$ be the subcategory of representations satisfying the
property (R2) also pictured in (\ref{fig-R}). Clearly, if $0\to
U\to V\to W\to 0$ is an exact sequence in $\Rep(L)$, $(V,\a_V)$
satisfies (R2) if and only if $(U,\a_U)$ and $(W,\a_W)$ satisfy
(R2).

\begin{proposition}
$\Rep(L)'$ is a braided monoidal category with braiding
$\Psi$, the same braiding as in $\calV$.
\end{proposition}

\Proof
$\Rep(L)'$ is closed under $\otimes$~: Let $(V,\a_V)$ and $(W,\a_W)$
satisfy (R2). We check that $(V\otimes W,\a_{V\otimes W})$ also
satisfies (R2)~:
\[\epsfbox{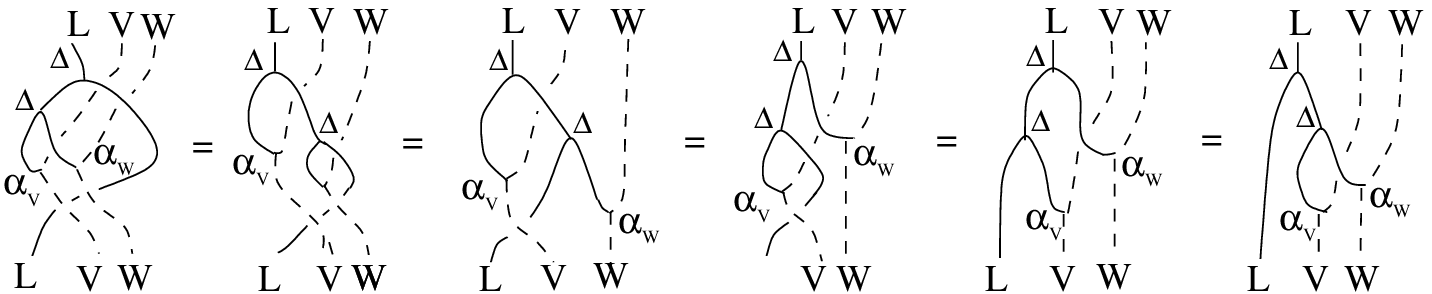}\]
The first, third and fifth equalities use the coassociativity of $\D$
and the naturality of $\Psi$, which holds since it already holds in $\calV$,
the second is (R2) for $(W,\a_W)$, the fourth is (R2) for $(V,\a_V)$.
Next, the braiding $\Psi_{V,W}$ is a morphism in $\Rep(L)'$~: we check the
equality $\Psi_{V,W}\circ \a_{V\otimes W}=
\a_{W\otimes V}\circ(\id\otimes \Psi_{V,W})$~:
\[\epsfbox{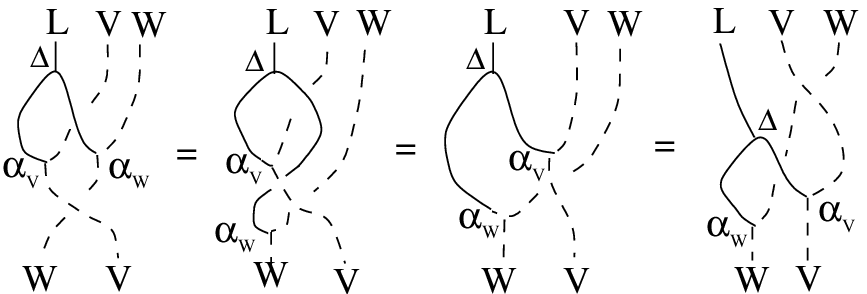}\]
The first and third equalities use the naturality of $\Psi$, the second
is property (R2) for $(V,\a_V)$. Note that (R2) is used (only for $V$),
therefore one cannot conclude anything for $\Rep(L)$ in general.
\qed
\\

\Remark This is an analogue at the Lie level of the braided
category of modules with respect to which a braided group behaves
cocommutatively (which in turn was the origin of (L2)), see
\cite{Majid-trans} for the general setting of that.\cqfd

 Obviously, $L_\Ad=(L,[\,,\,])$ is a representation of $L$
and satisfies $(R2)$ by assumption, so $L_\Ad\in \Rep(L)'$, and so
does the trivial representation. It is also clear from  the axioms
of a braided Lie algebra that the maps $\D: L_\Ad\to L_\Ad\uotimes
L_\Ad$, $[\,,\,]:L_\Ad\uotimes L_\Ad\to L_\Ad$ and $\e:L_\Ad\to
k$, are all intertwiners, and so is $\Psi_{L,L}:L_\Ad\uotimes
L_\Ad\to L_\Ad\uotimes L_\Ad$ by the previous proposition.
Therefore~:
$$
\Ups:L_\Ad\uotimes L_\Ad\to L_\Ad\uotimes L_\Ad \ is\ an\
intertwiner\
    in\ \Rep(L)'.
$$

\begin{proposition}
There exists a unique action, the adjoint action noted $x\rt_\Ad Y$,
of $L$ on $B(L)$ with the properties that
$x\rt y=[x,y]$ for $y\in L\hookrightarrow B(L)$
and that $B(L)$ is an algebra in $\Rep(L)'$ for this action.
Under this action, $B(L)=\oplus_{n\ge 0}B(L)_n$
is the direct sum of its homogeneous components.
\end{proposition}

\Proof The tensor algebra of the adjoint representation $L_\Ad$ is
an algebra in $\Rep(L)$ (as is always the case for the tensor
algebra of a representation) and, by the previous proposition,
belongs to $\Rep(L)'$ since $L_\Ad$ does. The ideal $\langle
\im(\id^{\otimes 2}-\Ups)\rangle$ is clearly graded, and a
subrepresentation by the above observation ($\Ups$ is an
intertwiner). Therefore the quotient $B(L)=T(L_\Ad)/\langle
\im(\id^{\otimes 2}-\Ups)\rangle$ is an algebra in $\Rep(L)'$ and
has the desired properties. Since $L$ generates as an algebra,
uniqueness is clear. \qed
\\

The adjoint action of $L$ on $B(L)$ defines an action of $B(L)$ on itself,
also noted $X\rt_\Ad Y$. The bracket $[\,,\,]_{B(L)}:B(L)\otimes B(L)\to B(L)$,
$[X,Y]_{B(L)}=X\rt_\Ad Y$, satisfies $(L1)$ and $(L2)$ by the above
proposition. It also satisfies (L3), since  both $[\,,\,]$
and $\Ups$ are morphisms of coalgebras. So we have~:

\begin{corollary}
The bialgebra $B(L)$ becomes a braided Lie algebra in $\Rep(L)'$
(also in $\calV$)
for the braided Lie product $[X,Y]_{B(L)}=X\rt_\Ad Y$.
Each graded summand $B(L)_n$ is a braided Lie subalgebra.
\end{corollary}

\Remark
Assume that $L$ is a good braided Lie algebra (\ie $L=\tgg$
is the extension of a good quantum Lie algebra $\gg$).
Recall that $U(\gg)\simeq B(\tgg)/\langle \ga-1\rangle$.
By taking appropriate quotients,
one easily gets that $U(\gg)$ is a left $U(\gg)$-module algebra,
and a braided Lie algebra in $\multi{\calM'}{U(\gg)}{}{}{}$ (or in $\calV$).
However, all statements concerning the grading are lost and
shoud be replaced by ``each term $U(\gg)_{(n)}$ of the natural filtration
of $U(\gg)$ is a braided Lie subalgebra of $U(\gg)$''.
But in turn, if $L=\tgg$ is the trivial extension of some braided Lie algebra
$\calL$, all the grading properties can be recovered,
since $U(\gg)\simeq B(\calL)$.
\cqfd

\subsection{The main example.}\label{main-example}
The definition of a braided Lie algebra was motivated in
\cite{Majid-braided-Lie} as follows. Let $(H,\m,\eta,\D,\e,S)$ be
a Hopf algebra in $\calV$. Its (left) braided
adjoint action is
\begin{equation}\label{def-braided-AdL}
\Ad_Lx(y)=x\rt_{\Ad} y=
    \m(\m\otimes S)(\id\otimes \Psi)(\D\otimes \id)(x\otimes y).
\end{equation}
We note $[x,y]:= x\rt_{\Ad} y$. $H$ is a left crossed module
over itself (in the braided sense) for the regular coaction $\D$
and the braided adjoint action $\Ad_L$. One easily checks that the
corresponding crossed module braiding $\Ups:= ([\, ,\,]\otimes
\id)(\id\otimes \Psi)(\D\otimes \id)$ satisfies $\m\circ\Ups=\m$
(this is the equality $xy=(x_{(1)}\,y\,Sx_{(2)})\, x_{(3)}$ when
$\Psi=\id_{H\otimes H}$). Moreover, $(\e\otimes \id)\circ
\Ups=\id\otimes \e$ and $(\id\otimes \e)\circ\Ups=[\,,\,]$.
Therefore there is a unique map $\s:\ker\e^{\otimes 2}\to
\ker\e^{\otimes 2}$ such that $\Ups(x\otimes y)=\s(x\otimes
y)+[x,y]\otimes 1$, $x,y\in \ker\e$. Multiplying this in $H$, we
get
$$
[x,y]=\m\circ (\Ups-\s)(x\otimes y)=\m\circ (\id-\s)(x\otimes y)
$$
for $x,y\in \ker\e$. This implies
$\ker(\id-\s)\subset \ker([\, ,\,]_{|(\ker\e)^{\otimes2}})$.
Therefore $(\ker\e,\s, [\,,\,])$ is always a quantum Lie algebra in $\calV$
(all maps are morphisms in $\calV$ by assumption, axioms 1-3
come from lemma \ref{lemma-braid-for-ts} applied to
$\gg=\ker\e$, $\ts=\Ups$, and axiom 4 (antisymmetry) holds by
the above equality). However, $(H,\D,\e,[\, ,\,])$ is not always
a braided Lie algebra in $\calV$.

\begin{proposition}\label{prop-H-Lie}
$H_L=(H,\D,\e,[\,,\,],\eta)$ is a good braided Lie algebra in $\calV$
if and only if axiom {\rm (L2)} is satisfied.
\end{proposition}

\Proof
As shown in \cite{Majid-braided-Lie},
axiom (L1) -braided Jacobi identity- is always satisfied
and, assuming  (L2), then (L3) is also satisfied.
Thus in this case, $H_L$ is a braided Lie algebra.
The Hopf algebra unit of $H$ is clearly a unit for $H_L$ in the sense
of (\ref{def-unital}). The good part (antisymmetry axiom)
has been checked above.
\qed
\\

Consider the case of a usual Hopf algebra $H$ (\ie a Hopf algebra in
$\multi{\calM}{k}{}{}{}$). Then axiom (L2) for $H_L$ is ensured
if $H$ is cocommutative, but fails to hold in general.
However, when the non-cocommutativity of $H$ is controlled by a
quasitriangular structure $\calR\in H\otimes H$, then axiom (L2)
remains valid, not for $H$ but for a braided version of $H$.

Let $(H,\calR)$ be a usual quasitriangular Hopf algebra
($\calR\in H\otimes H$ satisfies Drinfeld's axioms \cite{Drinfeld}) and
view $H$ as an object in $\multi{\calM}{H}{}{}{}$ by the left
adjoint action $\Ad_L$.
We will need the three braidings $\Ups$, $\Psi_\calR$ and
$\Xi_{\calR,\calR}$ on $H$ given below~:

\begin{lemma}\label{new-cross-H}
Let $\calR,\calS\in H\otimes H$ be two
co-quasitriangular structures on a (usual) Hopf algebra $H$. The coactions
$\D,\l_\calR,\d_{\calR,\calS}:H\to H\otimes H$ below define crossed module
structures on $(H,\Ad_L)$, with associated braidings $\Ups$, $\Psi_\calR$ and
$\Xi_{\calR,\calS}$ as indicated~:
\begin{equation}\label{braidings-on-H}
\begin{array}{lcl}
\D(x)=x_{(1)}\otimes x_{(2)} &,&
    \Ups(x\otimes y) = x_{(1)}\rt_{\Ad} y\otimes x_{(2)},
\\
\l_\calR(x)=\calR^{(2)}\otimes \calR^{(1)}\rt_{\Ad} x
    =\calR_{21}(1\otimes x)\calR_{21}^{-1} &,&
    \Psi_\calR(x\otimes y)=
    \calR^{(2)}\rt_{\Ad} y\otimes \calR^{(1)}\rt_{\Ad} x,
\\
\d_{\calR,\calS}(x)=\calR_{21}(1\otimes x)\calS &,&
    \Xi_{\calR,\calS}(x\otimes y)=
    ({\calR}^{(2)}\,\calS^{(1)})\rt_{\Ad} y\otimes
    \calR^{(1)}\, x\,{\calS}^{(2)}
\end{array}
\end{equation}
\end{lemma}

\Proof
The first one needs no comment. The second one is the image of
$(H,\Ad_L)$ under the monoidal functor
$\calF_\calR:\multi{\calM}{H}{}{}{}\hookrightarrow \multi{\calM}{H}{}{H}{}$
which sends arbitrary left module $(M,.)$ to $(M,.,\l_\calR^{(M)})$
where $\l^{(M)}_\calR(m)=\calR^{(2)}\otimes \calR^{(1)}. m$. The braiding
of $(M,.)$ calculated in $\multi{\calM}{H}{}{}{}$ (thanks to $\calR$) and
that of $(M,.,\l)$ in $\multi{\calM}{H}{}{H}{}$ are equal.
The reader will easily check that replacing the factor $\calR_{21}^{-1}$
in the definition of $\l_\calR$
by any other co-quasitriangular structure $\calS$ does not affect the
crossed module properties.
Note that $\l_\calR=\d_{\calR,\calR_{21}^{-1}}$, and therefore
$\Psi_\calR=\Xi_{\calR,\calR_{21}^{-1}}$.
Moreover, $(H,\Ad_L,\d_{\calR,\calS})$ is in the image of some
functor $\calF_\calT$ (if and) only if $\calS=\calR_{21}^{-1}$,
since $\l_\calT(1)=1\otimes 1$ and $\d_{\calR,\calS}(1)=\calR_{21}\calS$.
\qed
\\

Define the linear maps $\uD:H\to H\otimes H$, $\ue=\e$ and
$\uS:H\to H$ by
\begin{eqnarray}
\uD(x) & =&
    x_{(1)}\, S(\calR^{(2)})\otimes \calR^{(1)}\rt_{\Ad} x_{(2)}
    \,=:\, x_{(\underline{1})}\otimes x_{(\underline{2})}
    \label{eq-uD}
\\
\uS(x) &=& \calR^{(2)}\, S(\calR^{(1)}\rt_{\Ad} x)
\end{eqnarray}

\begin{proposition}\label{prop-H-L}
View $H$ as an object in $\multi{\calM}{H}{}{}{}$
via $\Ad_L$ (the braiding is $\Psi_\calR$).
\\
(i) $\uH=(H,\m,\eta,\uD,\ue,\uS)$
is a Hopf algebra in $\multi{\calM}{H}{}{}{}$.
It is $\Xi_{\calR,\calR}$-cocommutative in the sense that
$\Xi_{\calR,\calR}\circ \uD=\uD$. Its braided adjoint action
$\underline{\Ad}_L$ coincides
with the adjoint action $\Ad_L$ of $H$ and the crossed module braidings
$\Ups_{H}$ on $(H,\Ad_L,\D)$ and $\Ups_{\uH}$ on
$(\uH,\underline{\Ad}_L,\uD)$ also coincide.
\\
(ii) $\uH_L=(H,\uD, \ue,[\,,\,])$ is a left braided Lie algebra
in $\multi{\calM}{H}{}{}{}$ for the braided Lie bracket
$[x,y]=\Ad_Lx(y)$. Its canonical braiding is $\Ups$
{\rm (\ref{braidings-on-H})}, \ie the crossed module
braiding on $(H,\Ad_L,\D)$.
\\
(iii) Let $L\subset H$ satisfy $[H,L]\subset L$.
Then the following are equivalent~:
$$
(a)\;\; \D(L)\subset H\otimes L,\quad
(b)\;\; \uD(L)\subset H\otimes L,\quad
(c)\;\; \uD(L)\subset L\otimes H.
$$
If one of this condition holds,
then $L$ is a braided Lie subalgebra of $\uH_L$.
\end{proposition}

(We shall need the dual version which is more general,
so we omit this proof; see proposition \ref{propo-A=Lie-coalg}).

\subsection{Braided Lie coalgebras}\label{def-Lie-coalgebra}

The dual notion of left `braided Lie coalgebras' is just given in
the diagrammatic setting by turning the diagram-axioms of a right
braided Lie algebra upside-down, or by reflecting those of a left
braided Lie algebra about a horizontal axis and restoring braid
crossings, see \cite{Majid-lec}.

\begin{definition}\label{def-b-coLie}
A  (left) braided Lie coalgebra in $\calV$
is an algebra $(\calA,\mu,\eta)$
in the category endowed with a morphism (the braided Lie cobracket)
$\d:\calA\to \calA\otimes \calA$ satisfying the axioms below
\begin{eqnarray*}
\epsfbox{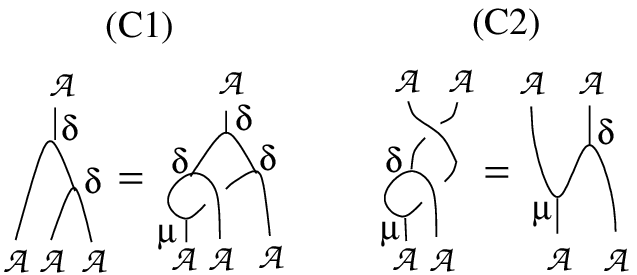} \qquad\epsfbox{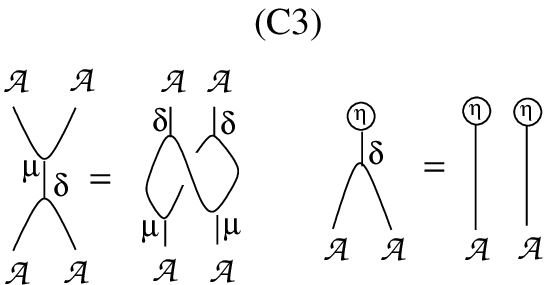}
\end{eqnarray*}
An ideal $\calI$ of $\calA$ is an algebra ideal such that
$\d(\calI)\subset \calA\otimes \calI+\calI\otimes \calA$.
A counit on $(\calA,\mu,\eta,\d)$ is a morphism of braided Lie algebras
$\e: \calA\to k$
 satisfying
\begin{equation}\label{def-counital}
(\e\otimes \id)\circ \d=\id,\qquad (\id\otimes \e)\circ \d=\eta\circ \e.
\end{equation}
\end{definition}

In the definition, $k$ is seen as a braided Lie coalgebra with
$\d(1)=1\otimes 1$. Obviously, $\calI$ is an ideal of the braided
Lie coalgebra $\calA$ iff the structure maps of $\calA$ induce a
braided Lie algebra structure on $\calA/\calI$. By turning proofs
upside-down, the following is also clear.
\\
\begin{lemma}
If $(\calA,\mu_\calA,\eta_\calA,\d_\calA)$ is a braided Lie
coalgebra, the  morphism $\Ups_\calA:\calA\uotimes \calA\to
\calA\uotimes \calA$,
\begin{equation}\label{Ups-A}
\Ups_\calA=
(\mu_\calA\otimes \id)\circ (\id\otimes \Psi)\circ (\d_\calA\otimes \id),
\end{equation}
is a morphism of algebras and satisfies the braid relation.
\end{lemma}

If $\calA$ has a right-dual $L$ in the categorical sense (in our
concrete setting it means if $\calA$ is finite-dimensional) there
are evaluation and coevaluation maps ${\rm ev}:\calA\otimes L\to
k$ and ${\rm coev}: k\to L\otimes\calA$, using which it follows by
diagrammatic methods\cite{Majid-lec} that $\calA$ is a left
braided-Lie coalgebra {\em iff} $L$ is a right braided Lie
algebra.

For our purposes we are interested in $\calA\in \calV={}^A\calM$
for some co-quasitriangular Hopf algebra $A$ but with $L$ regarded
in the braided category $\calM^A$. This is equivalent to the above
via the antipode of $A$, but in an algebraic setting it is more
natural to avoid the use of that. Thus we let $L=\calA^*$ be the
usual dual, viewed as a right $A$-comodule. We denote by
$\langle\, ,\,\rangle:\calA\otimes L\to k$ the evaluation pairing
and extend this (as for usual Hopf algebra duality) to
$(\calA\otimes \calA)\otimes (L\otimes L)\to k$ by setting
$$
\langle a\otimes b,x\otimes y\rangle := \langle a,x\rangle \;
    \langle b,y\rangle.
$$

\begin{lemma}\label{lemma-rel-Lie-coLie}
(i) Let $\calA\in {}^A\calM$ be finite-dimensional and
$L=\calA^*$. Then $(L,\D_L,\e_L,[\, ,\,]_L)$ is a braided Lie
algebra in $\calM^A$ if and only if
$(\calA,\mu_\calA,\eta_\calA,\d_\calA)$ is a braided Lie
coalgebra, with
\begin{equation}\label{L-A}
\langle \mu_\calA(a\otimes b),x\rangle :=\langle  a\otimes
b,\D(x)\rangle, \quad \langle \d_\calA(a),x\otimes y\rangle
:=\langle  a,[x,y]_L\rangle, \quad \langle
\eta_\calA(1),x\rangle :=\e_L(x)
\end{equation}
Moreover, $\eta_L:k\to L$ is a unit for $L$ if and only if
$\e_\calA:\calA\to k$ is a counit for $\calA$, where
$\e_\calA(a):=\langle a,\eta_L(1)\rangle$. (ii) Let
$L=\calA^*$ be unital. Then $L$ is the trivial extension of some
braided Lie algebra $\calL$ if and only if  $\ker\e_\calA$ is a
unital subalgebra of $\calA$ with unit $\theta$ satisfying
$\d_\calA(\theta)=1\otimes \theta$. Moreover,
$\calL=(1-\theta)^\bot$.
\end{lemma}

\Proof (i) This is a straightforward exercise from the
definitions. It is important to use compatible
conventions\cite{Majid-book} for the braidings of $\calM^A$ and
${}^A\calM$ as obtained from $\rrr:A\otimes A\to k$; one may check
that they are then adjoint. (ii) A coalgebra decomposition
$L=k\ga\oplus \calL$ is equivalent to an algebra decomposition
$$
\calA=(k\ga)^\bot\oplus \calL^\bot=\ker\e_\calA\oplus k\, \xi,
$$
for some vector $\xi$ spanning $\calL^\bot$.
Since $\calA$ is a unital algebra, so must
be $\ker\e_\calA$. Since $(\calA,\e_\calA)$ is a counital braided Lie coalgebra
(\ref{def-counital}),
one must have $\d_\calA(\xi)=1\otimes \xi+w$
for some $w\in (\ker\e_\calA)^{\otimes 2}$.
The hypothesis $[\calL,\calL]\subset \calL$ is then equivalent to
$w=0$, that is, $\d_\calA\xi=1\otimes \xi$.
Finally, if we normalize $\xi$  so that
$\langle \ga,\xi\rangle=1$, \ie $\e_\calA(\xi)=1$, then $\xi$ is an idempotent
(by definition of $\xi$ and because $\D(\ga)=\ga\otimes \ga$).
With this normalisation,  $\theta=1-\xi$ is the algebra unit of $\ker\e_\calA$
and
$\calL=\xi^\bot=(1-\theta)^\bot$.
\qed


\section{Link with differential calculi.}\label{Link}
\setcounter{equation}{0}

This section contains our main results, namely theorems connecting
the above results to bicovariant differential calculi on ordinary
Hopf algebras.

\subsection{Extended tangent spaces, inner calculi.}
Following \cite{Woronowicz},
a bicovariant first order differential calculus (bicovariant FODC)
over a Hopf algebra $A$ is a pair $(\G,\dd)$
where $\G$ is a Hopf bimodule,
with coactions $\D_L:\G\to A\otimes \G$ and $\D_R:\G\to \G\otimes A$,
and the linear map $\dd:A\to\G$ (the differential) satisfies
\begin{center}
\begin{tabular}{cl}
$\langle$1$\rangle$ &
    $\dd$ is a derivation~:
    $\dd(ab)=\dd(a)b+a \dd (b)$ for all $a,b\in A$,
\\
$\langle$2$\rangle$ &
    $\dd$ is a bicomodule map,
\\
$\langle$3$\rangle$ &
    the map $A\otimes A\to \G$,
    $a\otimes b\mapsto a\dd b$, is surjective.
\end{tabular}
\end{center}
Let $\pi_R:\G\to \G_R$ be the canonical projection on right invariants
(notations of the preliminaries).
$\G_R$ is called the (right) cotangent space of $\G$.
The differential $\dd$ and the (right handed) Maurer-Cartan map
$\o_R=\pi_R\circ \dd:A\to \G_R$ are related by
\begin{equation}\label{rel-d-omega}
\dd(a)=\o_R(a_{(1)}).a_{(2)}\; ,\quad \o_R(a)=\dd(a_{(1)})\,S(a_{(2)})
\end{equation}
and therefore are equivalent data.
Axioms $\langle 1-3\rangle$ for $(\G,\dd)$ are equivalent to
axioms $\langle 1'-3'\rangle$ below for the pair $(\G_R,\o_R)$~:
\begin{center}
\begin{tabular}{rl}
$\langle$1'$\rangle$ &
    $\o_R(ab)=a\rt \o_R(b)+\o_R(a)\, \e(b)$.
\\
$\langle$2'$\rangle$ &
    $\D_L \, \o_R=(\id\otimes \o_R)\ad_L$.
\\
$\langle$3'$\rangle$ &
    $\o_R:A\to \G_R$ is surjective.
\end{tabular}
\end{center}
($(\rt,\D_L)$ is the left crossed module structure of $\G_R$,
and $\ad_L:A\to A \otimes A$,
$\ad_L(a)=a_{(1)}\, S(a_{(3)})\otimes a_{(2)}$, is the left adjoint coaction).
A calculus $(\G,\dd)$ is called inner
if $\dd$ is an inner derivation, that is,
if there exists $\theta\in \G$ such  that for all $a\in A$
\begin{equation}\label{eq-inner}
\dd a=a\theta-\theta a \quad (\mbox{equivalently~:}\quad
\o_R(a)=a_{(1)}\,\theta\, S(a_{(2)})-\e_A(a)\, \theta).
\end{equation}

{\it Remark.} One can always assume (if necessary by replacing $\theta$
by $\pi_R(\theta)$), that $\theta\in\G_R$. Indeed,
apply $\pi_R$ to the right equality in (\ref{eq-inner}),
we get, using the properties of $\pi_R$,
$\o_R(a)=\pi_R\, \o_R(a)=
\pi_R(a_{(1)}\,\theta\, S(a_{(2)}))-\e_A(a)\pi_R(\theta)=
(a-\e_A(a))\rt \pi_R(\theta)$. Thus
$\theta'=\pi_R(\theta)$ has the same
property as $\theta$ and is right invariant.
\cqfd
\\

Any bicovariant FODC $(\G,\dd)$ can be extended to a pair
$(\widetilde{\G},\dd)$
which satisfies all axioms except $\langle 3\rangle$,
with the property that it contains
$\G$ as a Hopf sub-bimodule, and that the derivation
$A\stackrel{\dd}{\longrightarrow} \G\hookrightarrow \widetilde{\G}$
is inner; one takes
$\widetilde{\G}=\G\oplus \Theta. A$ as a right $A$-module
($\Theta$ a free variable),
with missing structures fixed by~:
$\Theta$  biinvariant
and left action $a \Theta=\dd a+\Theta. a$ in $\widetilde{\G}$.
$(\widetilde{\G},\tilde{d})$ is called the extended bimodule of $\G$
\cite{Woronowicz},
and the crossed module $(\widetilde{\G}_R,\rt,\tilde{\D}_L)$
of right invariants of $\widetilde{\G}$ is called
the (right) extended cotangent space of $\G$.
Let  $\tilde{\pi}_R:\widetilde{\G}\to \widetilde{\G}_R$
be the canonical projection and
define
$\tilde{\o}_R:A\to \widetilde{\G}_R$ by
$$
\tilde{\o}_R(a)=a\rt \Theta=\o_R(a)+\e_A(a)\, \Theta
$$

Let $(A,\m,\ad_L)$ be the left crossed $A$-module,
where the left action is the left regular one  ($a\rt b=\m(a\otimes b)=ab$).
Axioms $\langle 1'-3'\rangle$ can again be restated
as the fact that
$\tilde{\o}_R:\widetilde{A}\to \widetilde{\G}_R$ is a surjective
crossed module homomorphism; equivalently (since $\o_R(1)=0$), that
the restriction $\o_R:\ker\e_A\to \G_R$ is a surjective
crossed module homomorphism.
Therefore, $\tilde{\o}_R$ and $\o_R$ respectively
induce crossed module isomorphisms
\begin{equation}
(A,\m,\ad_L)/\calI_\G\stackrel{\simeq}{\longrightarrow}\widetilde{\G}_R,
\qquad
\ker\e_A/\calI_\G\stackrel{\simeq}{\longrightarrow}\G_R,
\end{equation}
where
$\calI_\G:= \ker\tilde{\o}_R=\ker\o_R\cap \ker\e_A$
is called ``the left ideal associated to $\G$''.
The case $\calI_\G=0$ corresponds to the
universal extended cotangent space
$\widetilde{\G}_{R,\univ}=A$, with $\tilde{\o}_{R,\univ}(a)=a$
(and therefore $\G_{R,\univ}=\ker\e_A$, ${\o}_{R,\univ}(a)=a-\e_A(a)1$).
By definition, the (right) {\it tangent space} $\gg_\G$ and
{\it extended tangent space} $\tgg_\G$ of $\G$ are
\begin{eqnarray}
\gg_\G &=& \{x\in A^*:  x(\calI_\G)=0\;\; {\rm and}\;\; x(1)=0\},\label{def-T}
\\
\tgg_\G &=& \{X\in A^*: X(\calI_\G)=0\},\label{def-ext-T}
\end{eqnarray}
As subspaces of $A^*$, one has
$\tgg_\G=k\, 1_{A^\circ}\oplus \gg_\G$ and
$\gg_\G=\tgg_\G\cap \ker\e_{A^\circ}$.
Let $[\,,\,]:A^\circ\otimes A^\circ\to A^\circ$ be defined by
$$
[X,Y]:= X_{(1)}\, Y\, S(X_{(2)}).
$$
Recall from \cite{Woronowicz} that there is a unique bilinear form
$(\, ,\,):\G\times \gg_\G\to k$ such that
$(\o_R(a).b,x)=\e(b)\langle a,x\rangle$ and,
if $\dim_k \G_R<\infty$, it allows to identify $\gg_\G$ with
$(\G_R)^*$ as right crossed modules over $A$.
In this case, by axioms $\langle 1',2'\rangle$ of a bicovariant FODC,
the right action $\lt$ -such that $(a\rt\o_R(b),x)=(\o_R(b),x\lt a)$,
coaction $x\mapsto x^{(0)}\otimes x^{(1)}$ on $\gg_\G$, and
the corresponding crossed module braiding $\s$  are defined by~:
\begin{eqnarray}
&&x\lt a=\langle a,x_{(1)}\rangle \, x_{(2)}-\langle a,x\rangle 1_{A^\circ},
\quad
x^{(0)}\, \langle x^{(1)},h\rangle = [h,x]
    \label{def-lt}
\\
&&
\s(x\otimes y)=[x_{(1)},y]\otimes x_{(2)}-[x,y]\otimes 1_{A^\circ}.
    \label{def-s}
\end{eqnarray}
for all $h\in A^\circ$. The direct analogue of this is
(see for instance {\rm \cite{KS}})~:

\begin{proposition}\label{prop-ext-T}
(i) There is a unique bilinear form
$(\, ,\,\tilde{)}:\widetilde{\G}\times \tgg_\G\to k$ such that
$(\tilde{\o}_R(a).b,X\tilde{)}=\langle a,X\rangle\, \e(b)$.
(ii) Assume that $\dim A/\calI_\G < \infty$.
Then $\tgg_\G\subset A^\circ$ and
$\tgg_\G$ has the following properties~:
(a) $1_{A^\circ}\in \tgg_\G$,
(b) $\D(\tgg_\G)\subset A^\circ \otimes \tgg_\G$,
(c) $[A^\circ, \tgg_\G]\subset \tgg_\G$.
Conversely, if $A^\circ$ separates the elements of $A$,
a subspace $\tgg_\G$ of $A^\circ$  satisfying
the properties (a), (b) and (c) is the extended
tangent space of a unique (up to isomorphism)
bicovariant FODC over $A$, with associated ideal
$\calI_\G=\{a\in A| \forall X\in \tgg_\G, X(a)=0\}$.
Moreover,
$\tgg_\G\simeq(\widetilde{\G}_R)^*$ as right crossed modules over $A$,
with right action $\leftharpoonup$,
coaction $X\mapsto X^{(0)}\otimes X^{(1)}$ and crossed module
braiding $\ts$ below~:
\begin{equation}\label{br-ext-T}
X\leftharpoonup a=\langle a,X_{(1)}\rangle \, X_{(2)},
\quad
X^{(0)}\, \langle X^{(1)},h\rangle = [h,X],\quad
\tilde{\s}(X\otimes Y)=[X_{(1)},Y]\otimes X_{(2)},
\end{equation}
for all $h\in A^\circ$.
\end{proposition}

Note that $\gg_\G$ is not a crossed submodule of $\tgg_\G$;
it is rather isomorphic to the quotient $\tgg_\G/k1_{A^\circ}$,
via the projection $\varphi:\tgg_\G\to \gg_\G$, $X\mapsto X-\e(X)1_{A^\circ}$.
Thus, one has
$x\lt a=\varphi(x\leftharpoonup a)$,
$\s(x\otimes y) =(\varphi\otimes \varphi)\ts(x\otimes y)$.
Woronowicz (\cite{Woronowicz}, Th. 5.3 and Th. 5.4)
has shown that the triple $(\gg_\G, \s,[\, ,\,])$
satisfies the axioms of a (left) quantum Lie algebra.
This can be recovered by lemma \ref{lemma-braid-for-ts}
for axioms 1-3, and antisymmetry (axiom 4) comes from (\ref{def-s})
- after multiplication in $A^\circ$, it gives
\begin{equation}\label{bra-vs-m}
[x,y]=xy-\m_{A^\circ}\circ \s(x\otimes y),\qquad (x,y\in\gg_\G).
\end{equation}
One can be more precise~:
$(\gg_\G,\s,[\,,\,])$ is a quantum Lie algebra in
$\multi{\calM}{}{}{}{A}$, but in general not in
$\multi{\calC}{}{A}{}{A}$.
The braiding $\s$ makes no problem since it is the crossed module
braiding on $\gg_\G$ in $\multi{\calC}{}{A}{}{A}$.
The quantum Lie bracket is also a morphism in $\multi{\calM}{}{}{}{A}$~:
this follows from the identity
$[h,[x,y]]=[[h_{(1)},x],[h_{(2)},y]]$ for all $h,x,y\in A^\circ$.
However, the quantum Lie bracket is not $A$-linear since
$\D([x,y])\ne [x_{(1)},y_{(1)}]\otimes [x_{(2)},y_{(2)}]$ in general.

Finally, we note from (\ref{bra-vs-m})
that if $\jmath:\gg_\G\hookrightarrow U(\gg_\G)$
is the natural imbedding, then there is a unique algebra homomorphism
$U(\gg_\G)\to A^\circ$ such that $\jmath(x)\mapsto x$ for all
$x\in \gg_\G$. This algebra homomorphism needs not be injective
nor surjective in general.

By lemma \ref{lemma-Ug-quad}, there is an algebra isomorphism
$S_{\tilde{\s}}(\tgg_\G)/\langle 1_{A^\circ}-1\rangle
\simeq U(\gg_\G)$. (The role of $\ga$ is played by $1_{A^\circ}\in\tgg_\G$,
which is an
element of degree 1 in the tensor algebra of $\tgg_\G$ and should not be
confused with the unit element $1\in k\subset T(\tgg_\G)$).
Also by lemma \ref{lemma-Ug-quad},
$U(\gg_\G)$ can sometimes be itself a quantum symmetric algebra.
We show that
this happens when the calculus is inner.

Let $\s^t:\G_R^{\otimes 2}\to \G_R^{\otimes 2}$
be the crossed module braiding on $\G_R\simeq (\tgg_\G)^*$.
By definition, the quadratic extension of $\G_R$ is
$$
\G_R^{\w, quad}:= T(\G_R)/
    \langle \ker(\id^{\otimes 2}+\s^t)\rangle.
$$
It is known that the crossed product
$\G^{\w,quad}=\G_R^{\w, quad} \lcross\, A$ has a structure
of a graded differential Hopf algebra \cite{Brzezinski} and
maps onto Woronowicz' external algebra $\G^\w$ \cite{Woronowicz}.
In some cases, it coincides with it. This is the case for instance
for the standard $n^2$-dimensional bicovariant FODC on $GL_q(n)$
and $SL_q(n)$ \cite{Schueler}.

Below we let $\G_\inv=\G_L\cap\G_R=$
$\{\eta\in \G| \D_L(\eta)=1\otimes \eta,
\D_R(\eta)=\eta\otimes 1\}$.

\begin{theorem}\label{theo-quadratic}
Let $(\G,\dd)$ be a finite dimensional
bicovariant FODC over some Hopf algebra $A$.
If there exists $\theta\in \G_\inv$ such that
$\dd(a)=a\theta-\theta a$ for all $a\in A$, then
\begin{equation}\label{Ug-quad-dual}
U(\gg_\G)\simeq (\G_R^{\w, quad})^!
\end{equation}
is isomorphic to the quadratic dual of $\G_R^{\w, quad}$.
\end{theorem}

First, the hypothesis of the theorem
can be interpreted as follows.

\begin{lemma}\label{lemma-inner}
There is a 1-1 correspondence between
\begin{enumerate}
\item
$\theta\in\G_\inv$
    such that $\dd(a)=a\theta-\theta a$ for all $a\in A$,
\item
elements $\htheta\in\ker\e_A\,(\bmod \calI_\G)$ satisfying
    $a\htheta\equiv a\bmod \calI_\G$ for all $a\in \ker\e_A$ and
    $\ad_L(\htheta)\equiv 1\otimes \htheta\bmod A\otimes \calI_\G$,
\item
subspaces $\calL$ of $\tgg_\G$ satisfying
    $\tgg_\G=k 1_{A^\circ}\oplus \calL$,
    $[A^\circ,\calL]\subset \calL$ and
    $\D(\calL)\subset A^\circ \otimes \calL$.
\end{enumerate}
It is given by
\begin{equation}\label{rel-theta-L}
\o_R(\htheta)=\theta,\qquad
\calL=\{x\in \tgg_\G|\langle 1-\htheta, x\rangle =0\}=(\Theta-\theta)^\bot.
\end{equation}
(The orthogonality is with respect to the
bilinear form $(\,,\,\tilde{)}$ on $\widetilde{\G}\times \tgg_\G$).
\end{lemma}

\Proof
The equivalence
$\theta\Leftrightarrow \htheta$ follows directly
from the equivalence of sets $\langle 1-3\rangle$ and
$\langle 1'-3'\rangle$ of axioms of a bicovariant FODC.
Equivalence $\htheta\Leftrightarrow \calL$~:
Let $\htheta$ have the given properties,
and set $\calI'=\calI_\G\oplus k(1-\htheta)$.
The sum is direct since $\e_A(\calI_\G)=0$ and $\e_A(1-\htheta)=1$.
This implies that $\calL$ has codimension 1 in $\tgg_\G$, with possible
complement $k1_{A^\circ}$. $\calI'$ is obviously a left ideal
of $A$ closed under $\ad_L$, therefore $\calL$ is a left co-ideal
of $A^\circ$, invariant under the left adjoint action of $A^\circ$.
Conversely, let $\calL$ have the given properties.
Then $\calI'=\{a\in A| \langle a,\calL\rangle =0\}\supset \calI_\G$
is a left ideal of $A$ closed under $\ad_L$, and
$\calI_\G=\calI'\cap \ker\e_A$.
Therefore $\calI_\G$ has codimension 1 in $\calI'$, \ie there exists
$\xi\in A$, $\xi\notin \ker\e_A$, such that $\calI'=\calI_\G\oplus k\xi$.
We normalize $\xi$ such that $\e_A(\xi)=1$ and set
$\htheta=1-\xi\in\ker\e_A$.
Since $\calI_\G\subset \calI'$ are left ideals,
one must have for all $a\in \ker\e_A$,
$a\xi\in \calI'\cap \ker\e_A=\calI_\G$ (\ie $a\htheta\equiv a\bmod \calI_\G$).
Since $\calI_\G\subset \calI'$ are closed under $\ad_L$, one must have
$\ad_L(\xi)=a\otimes \xi+w$ for some $a\in A$ and $w\in A\otimes \calI_\G$.
By the counit axioms and from $\e_A(\xi)=1$, one has
$a=(\id\otimes \e_A)\ad_L(\xi)=\e_A(\xi)=1$.
\qed
\\

{\it Proof of the theorem.}
Define $\htheta$ and $\calL$ as in the lemma. $\calL$ satisfies
$\tgg_\G=k 1_{A^\circ}\oplus \calL$ and $\ts(\calL\otimes \calL)\subset
\calL\otimes \calL$ (since $[A^\circ,\calL]\subset \calL$ and
$\D(\calL)\subset A^\circ\otimes \calL$). Therefore,
by lemma \ref{lemma-Ug-quad},
$U(\gg_\G)\simeq S_\s(\gg_\G)=T(\gg)/\langle \im (\id^{\otimes 2}-\s)\rangle$,
whose quadratic dual is by definition
$T(\gg^*)/\langle \ker(\id^{\otimes 2}+\s^t)\rangle$,
\ie $\G_R^{\w,quad}$.
\qed
\\

Note that the simplest way to construct a bicovariant FODC
is to pick some $\ad_L$-invariant element $a$ and
a left ideal $J$ of $A$, such that either $a$ or $J$ belongs to $\ker \e_A$,
and set $\calI_\G=J a$.
In particular, we achieve the hypothesis of the theorem if we take
$\calI_\G=\ker\e_A(1-\htheta)$ with $\htheta\in\ker\e_A$  $\ad_L$-invariant.
Another known construction of
bicovariant FODC is by picking some central element $c$ of $A^\circ$.
We identify when this calculus is inner.

Let $h\leftharpoonup a=\langle a,h_{(1)}\rangle\, h_{(2)}$
be the right (co)regular action of $a\in A$ on $h\in A^\circ$.
One has~:

\begin{lemma}\label{lemma-central}
Let $c$ be central in $A^\circ$ and define $\calL_c := c\leftharpoonup A$.
\\
(i) {\rm \cite{Majid-classif-bicov}\cite{KS}}
$\tgg(c)=k1_{A^\circ}+\calL_c$ is the extended tangent space
of a bicovariant FODC $\G(c)$ over $A$.
(ii) If $1_{A^\circ}\notin \calL_c$, then $\G(c)$ is inner, with differential
implemented by a biinvariant element.
\end{lemma}

\Proof $(i)$ is well-known  but we need its proof for $(ii)$. We
apply lemma to $\xi(a,b):= \langle a \otimes b,\D(c)\rangle$.
It obviously satisfies $\m^{op}*\xi=\xi *\m$ and
$\im(\xi_1)\subset A^\circ$. So, for $h\in A^\circ$ and  $a\in A$,
one has $\Ad_L h(c\leftharpoonup a)= \Ad_L
h(\xi_1(a))=\xi_1(\Ad_L^*h(a))=c\leftharpoonup \Ad_L^*h(a)$, \ie
$\calL_c$ is a submodule of $A^\circ$ for $\Ad_L$. Since $\calL_c$
is a submodule for the right (co)regular action of $A$ by
hypothesis, it is a left co-ideal of $A^\circ$.  Therefore by
lemma \ref{prop-ext-T}, $\tgg(c)=k 1_{A^\circ}+\calL_c$ is an
extended tangent space for $A$. Under the hypothesis of $(ii)$,
one $\tgg(c)=k1\oplus \calL_c$, therefore by lemma
\ref{lemma-inner}, the element $\theta_c\in \G_R(c)$ uniquely
determined by
$k(\Theta-\theta_c)=\calL_c^\bot\cap\widetilde{\G}_R$ implements
$\dd$ and is biinvariant. \qed

\begin{lemma}
The exists $\theta\in\G_R$ (resp. $\theta\in\G_\inv$) such that
$\dd(a)=a\theta-\theta a$ for all $a\in A$  if and only if the imbedding
$\G_R\hookrightarrow \widetilde{\G}_R$ splits in
$\multi{\calM}{A}{}{}{}$ (resp in $\multi{\calC}{A}{}{A}{}$).
\end{lemma}

\Proof
Since $\widetilde{\G}_R/\G_R\simeq k$ is  trivial both as module
and comodule, the imbedding
$\G_R\hookrightarrow \widetilde{\G}_R$ splits in
$\multi{\calM}{A}{}{}{}$ iff
$\widetilde{\G}_R=\G_R\oplus k\xi$ for some $\xi$ spanning the trivial
$A$-module. Given such a $\xi$, normalized so that $\Theta-\xi\in \G_R$,
then $\theta=\Theta-\xi$ satisfies
$\o_R(a)=(a-\e(a))\rt\Theta=(a-\e(a))\rt\theta$ for all $a$.
Conversely, define $\xi=\Theta-\theta$.
Finally, the imbedding
$\G_R\hookrightarrow \widetilde{\G}_R$ splits in $\multi{\calC}{A}{}{A}{}$
iff it splits in both $\multi{\calM}{A}{}{}{}$ and $\multi{\calM}{}{}{A}{}$
(since $k$ is simple is both), so that the above $\theta$ and $\xi$ must
be also left invariant.
\qed
\\

{}From the lemma, we see that $\G_\univ$ (and all its quotients)
is inner if and only if $A$ is semi-simple, hence finite
dimensional. For infinite dimensional Hopf algebras, $\G_\univ$
cannot be inner, but all finite dimensional bicovariant FODC over
$A$ are inner when the category $\multi{\calM}{A}{}{}{(f)}$ of
finite dimensional left $A$-modules, or
$\multi{\calC}{A}{}{A}{(f)}$ of finite dimensional left crossed
modules, is semi-simple. For $A=\calO(G)$, the algebra of
polynomial functions on some matrix group $G$, neither
$\multi{\calM}{A}{}{}{(f)}$ nor $\multi{\calC}{A}{}{A}{(f)}$ are
semi-simple, and there are non inner differential calculi. For the
quantizations $\calO_q(G)$, $q$ not a root of unity, it is known
that, at least for $G=SL(n)$ or $G=Sp(n)$ \cite{HS}, all finite
dimensional bicovariant FODC are semi-simple and inner. This is an
indication that $\multi{\calC}{A}{}{A}{(f)}$ is semi-simple in
this case, although there is apparently  still no proof of this.


\subsection{The co-quasitriangular case.}\label{sect-main-results}
{}From now on, $(A,\rrr)$ is a co-quasitriangular Hopf algebra. We
work in the category $\multi{\calM}{}{}{A\, }{}$ of left
$A$-comodules. $A$ itself will always be a left $A$-comodule via
the left adjoint coaction $\ad_L$ ($a^{(-1)}\otimes a^{(0)}:=
\ad_L(a)$ in the following). The analogues of the facts given in
section \ref{main-example} for a quasitriangular Hopf algebra
$(H,\calR)$ are as follows. The (exact, monoidal) functor
$\calF_\rrr:\multi{\calM}{}{}{A\, }{}\to
\multi{\calC}{A\,}{}{A\,}{}$ now sends a left $A$-comodule
$(M,\d_L)$ to the left crossed module $(M,\brt_\rrr,\d_L)$, where
$a\brt_\rrr m:= \langle m^{(-1)}, \rrr_2(a)\rangle\, m^{(0)}$.
For $(A,\ad_L)$, this gives~:
\begin{eqnarray}
a\brt_\rrr b &=& \langle b^{(-1)},\rrr_2(a)\rangle\; b^{(0)}
    \,=\, \rrr(b_{(1)},a_{(1)})\, b_{(2)}\, \bar{\rrr}(b_{(3)},a_{(2)}).
    \label{eq-PsiR-for-A}
\\
\Psi_\rrr(a\otimes b) &=& a^{(-1)}\brt_\rrr b\otimes a^{(0)}
    = \langle b^{(-1)},\rrr_2(a^{(-1)})\rangle\; b^{(0)} \otimes a^{(0)}
    \nonumber
\\
    &=&
    b_{(3)}\otimes a_{(3)}\,
    \rrr(b_{(2)},a_{(1)})\, \rrr(b_{(1)},Sa_{(5)})\,
    \bar{\rrr}(b_{(4)},a_{(2)})\, \rrr(b_{(5)},a_{(4)})\label{Psi-r}
\end{eqnarray}
The analogue of lemma \ref{new-cross-H} is~:

\begin{lemma}\label{lemma-new-cross}
Let $\rrr$ and $\sss$ be two co-quasitriangular structures on $A$.
\\
The map $A\to A\otimes A^\circ$,
    $a\mapsto a^{[0]}\otimes a^{[1]}=
    a_{(2)}\otimes \rrr_1(a_{(1)})\,\sss_2(a_{(3)})$,
is a right coaction of $A^\circ$ on $A$.
Let $a\rightharpoonup_{\rrr,\sss} b= b^{[0]}\, \langle a,b^{[1]}\rangle$
be the corresponding left action of $A$ on itself.
Then $(A,\rightharpoonup_{\rrr,\sss},\ad_L)$ is a left crossed $A$-module
with braiding $\Xi_{\rrr,\sss}$ below~:
\begin{eqnarray}
\Xi_{\rrr,\sss}(a\otimes b) &=& a^{(-1)}\rightharpoonup_{\rrr,\sss}
        b\otimes a^{(0)}=
         b^{[0]}\otimes a^{(0)}\,
            \langle a^{(-1)}, b^{[1]}\rangle
\label{eq-Xi}
\end{eqnarray}
Moreover, $(A,\rightharpoonup_{\rrr,\sss},\ad_L)$ is in the image of
the functor $\calF_\ttt$, for some co-quasitriangular structure $\ttt$ on $A$,
(if and) only if $\sss=\bar{\rrr}_{21}$.
In particular, $\Psi_\rrr=\Xi_{\rrr,\bar{\rrr}_{21}}$.
\end{lemma}

Define the linear maps
$\um:A\, \otimes\, A\to A$,
$\underline{\eta}=\eta:k\to A$,
$\underline{S}:A\to A$ by
\begin{eqnarray}
\um(a\otimes b)
    \, =\,
    a\covdot b &:=& a_{(1)}(S(a_{(2)})\brt_\rrr b)
    \, =\,
    a_{(1)}\, b_{(2)}\,
    \langle b_{(1)}\, S(b_{(3)}), \rrr_2 S(a_{(2)})\rangle
    \label{def-covdot}
\\
\underline{S}(a) &=& a_{(1)}\brt_\rrr S(a_{(2)})
\end{eqnarray}

\begin{proposition}\label{propo-A=Lie-coalg}
(i) $\underline{A}\equiv
    (A,\um,\underline{\eta},\D,\e,
    \underline{S})$ is a Hopf algebra in
    $\multi{\calM}{}{}{A\, }{}$, where the braiding is $\Psi_\rrr$.
It is $\Xi_{\rrr,\rrr}$-commutative,
in the sense that $\um\circ \Xi_{\rrr,\rrr} =\um$.
Its braided adjoint
coaction $\underline{\ad}_L$ coincides with the adjoint coaction $\ad_L$,
and the crossed module braiding on $(\uA,\um,\underline{\ad}_L)$ coincides
with crossed module braiding on $(A,\m,\ad_L)$.
\\
(ii) $(\underline{A},\um,\ueta,\d)$ is a braided Lie coalgebra
in $\multi{\calM}{}{}{A}{}$
for the braided Lie coproduct $\d=\ad_L$.
It is counital with counit $\e_A$.
Its canonical  braiding $\Ups_{\uA}$ (see {\rm (\ref{Ups-A})}) coincides
with the crossed module braiding on $(A,\m,\ad_L)$.
\\
(iii) Let $\calI$ be a subcomodule of $(A,\ad_L)$. Then $(\calI,\ad_L,\brt_\rrr)$
is a crossed submodule of $\uA=(A,\ad_L,\brt_\rrr)$,
\ie $A\brt_\rrr \calI\subset \calI$,
and the following are equivalent~:
$$
(a)\;\; A\calI\subset \calI,
\qquad
(b)\;\; A\covdot\calI\subset \calI,
\qquad
(c)\;\; \calI\covdot A\subset \calI
$$
\end{proposition}

\Proof $(i)$ is well-known \cite{Majid-book} although the version
of $\uA$ given there is in the category of right $A$-comodules.
The $\Xi_{\rrr,\rrr}$-commutativity of $\uA$ is the equivalent of
the quasi-commutativity of $A$ ($\m^{op}*\rrr= \rrr * \m$). The
new observation is that $\Xi_{\rrr,\rrr}$ is a crossed module
braiding. Note that if $\bar{\rrr}_{21}=\rrr$, then
$\Xi_{\rrr,\rrr}=\Psi_{\rrr}$, \ie $\uA$ is a commutative algebra
in $\multi{\calM}{}{}{A}{}$. Finally, the braided Hopf algebra
$\uA$ is defined in \cite{Majid-book} by the requirements that
$\underline{\ad}_L=\ad_L$ and that the crossed module braidings on
$(\uA,\um,\underline{\ad}_L)$ and $(A,\m,\ad_L)$ coincide, so
these statements are just reminders.
\\
$(ii)$ Clearly, if $\uA$ is a braided Lie coalgebra,
its canonical braiding (see {\rm (\ref{Ups-A})}) is the crossed module
braiding on $(\uA,\um,\underline{\ad}_L)$, \ie on $(A,\m,\ad_L)$ by $(i)$,
that is ~: $\Ups_{\uA}(a\otimes b)=a^{(-1)}b\otimes a^{(0)}$.
According to (a dual version of) proposition \ref{prop-H-Lie},
we only need to check that axiom (C2) holds. It can be expressed as the
equality $\Psi_\rrr\circ \Ups_{\uA}=(\um\otimes \id)\circ (\id\otimes\ad_L)$,
Thus, we need to check that
$\rrr(b^{(-2)},a^{(-1)})\, b^{(-1)}a^{(0)}\otimes b^{(0)}=
    a\covdot b^{(-1)}\otimes b^{(0)}$ for all $a,b\in A$.
This is implied by the equality
$\rrr(b_{(1)},a^{(-1)})\, b_{(2)}a^{(0)}=a\covdot b$, which holds since
\begin{eqnarray*}
a\covdot b &=& a_{(1)}\, b_{(2)}
    \langle b_{(1)}\, S(b_{(3)}),\rrr_2S(a_{(2)})\rangle
\\
 &=& \underline{ a_{(1)}\, b_{(2)}}\,
    \langle b_{(1)},\rrr_2S(a_{(3)})\rangle\,
    \underline{\langle S(b_{(3)}),\rrr_2S(a_{(2)})\rangle}
\\
 &=& b_{(3)}a_{(2)}\, \langle b_{(1)},\rrr_2S(a_{(3)})\rangle
    \langle b_{(2)},\rrr_2(a_{(1)})\rangle
\\
 &=& b_{(2)}a_{(2)}\, \langle b_{(1)},\rrr_2(a_{(1)}\,Sa_{(3)})\rangle
    =\rrr(b_{(1)},a^{(-1)})\, b_{(2)}a^{(0)}.
\end{eqnarray*}
The underlined terms are changed using $S\rrr_2S=\rrr_2$ and
$\m^{op}*\rrr=\rrr*\m$.
\\
 $(iii)$ The first claim follows from the covariance of the functor
$\calF_\rrr:\multi{\calM}{}{}{A\, }{}\to
    \multi{\calC}{A\,}{}{A\,}{}$.
Then the  equivalence of $(a)$ and $(b)$ follows directly
from $\ad_L(\calI)\subset A\otimes \calI$ and the relation between
$\m$ and $\um$ (the reverse being
$ab=a_{(1)}\covdot (a_{(2)}\brt_\rrr b)$).
The equivalence of $(b)$ and $(c)$ follows for the
$\Xi_{\rrr,\rrr}$-commutativity property of $\uA$~:
If $\ad_L(\calI)\subset A\otimes \calI$, then
$\Xi_{\rrr,\rrr}(\calI\otimes A)\subset A\otimes \calI$;
this is in fact
an equality since $\Xi_{\rrr,\rrr}$ is a crossed module braiding and
the antipode of $A$ is invertible.
Therefore
$\Xi_{\rrr,\rrr}(\calI\otimes A)= A\otimes \calI$,
and  $A\covdot \calI=\calI\covdot A$.
\qed
\\

{\it Remarks.}
$(i)$ If $\rrr=\calR\in H\otimes H$, with $H=A^\circ$,
then the braided Hopf algebra $\uA$ is dual to the braided Hopf algebra
$\uH$ of proposition \ref{prop-H-L} in the sense
that $\langle \Psi_\rrr(a\otimes b),x\otimes y\rangle=
\langle a\otimes b,\Psi_\calR(x\otimes y)\rangle$,
$\langle a\covdot b,x\rangle=\langle a\otimes b,\uD(x)\rangle$
and
$\langle \uS(a),x\rangle = \langle a,\uS(x)\rangle$ for all $a,b\in A$,
$x,y\in H$, where the pairing between $\uA\otimes \uA$ and
$\uH\otimes \uH$ is given by
$\langle a\otimes b, x\otimes y\rangle=
\langle a,x\rangle\; \langle b,y\rangle$ (this is the opposite convention of
\cite{Majid-book}).
\\
$(ii)$ Following \cite{Schmudgen}, a central bicharacter on $A$ is
a convolution invertible map $\ccc :A\otimes A\to k$ such that
$\ccc(ab,c)=\ccc(a,c_{(1)})\ccc(b,c_{(2)})$,
$\ccc(a,bc)=\ccc(a_{(2)},b)\ccc(a_{(1)},c)$; $\ccc(a,b_{(1)})\,
b_{(2)}=\ccc(a,b_{(2)})\, b_{(1)}$ and $\ccc(a_{(1)},b) a_{(2)}=
a_{(1)}\ccc(a_{(2)},b)$ for all $a,b,c\in A$. If $\ccc$ is a
central bicharacter, then $\ccc*\rrr$ is also a co-quasitriangular
structure on $A$. But one easily checks that changing $\rrr$ in
$\ccc*\rrr$ does not affect $\um$.
\\
$(iii)$ If $C$ is a subcoalgebra of $A$,
then $C$ generates $\underline{A}$
as an algebra if and only if $C$ generates $A$ as an algebra
(this follows from the reciprocal relations between $m$ and $\um$,
and the fact that subcoalgebras are subcomodules of $(A,\ad_L)$).
\\
$(iv)$ The subspace $A^{\ad_L}=\{a\in A:\, \ad_L(a)=1\otimes a\}$
is a commutative subalgebra of $A$, and belongs to the center of $\uA$.
Indeed, for $b\in A^{\ad_L}$ and all $a\in A$ one has
$b\covdot a=a\covdot b=ab$ (the first equality comes from
$\Xi_{\rrr,\rrr}(b\otimes a)=a\otimes b$) and also
$\ad_L(ab)=\ad_L(a)\, \ad_L(b)$,
$\ad_L(a\covdot b)=\ad_L(a)\covdot \ad_L(b)$ (the last equality actually
holds for all $a,b$ by axiom (C3)).
\cqfd

\begin{theorem}\label{main-theo}
Let $(A,\rrr)$ be co-quasitriangular,
$(\G,d)$ be a  bicovariant FODC over $A$.
\\
(i) $\widetilde{\G}_R\simeq A/\calI_\G$ is a counital
braided Lie coalgebra in $\multi{\calM}{}{}{A}{}$,
when regarded as $\uA/\calI_\G$.
\\
(ii) If $(\G,d)$ is finite dimensional,
$\tgg_\G$ is a unital braided Lie algebra in $\multi{\calM}{}{}{}{A}$.
Moreover,
$U(\gg_\G)\simeq B(\tgg_\G)/\langle 1_{A^\circ}-1\rangle$
is a bialgebra in $\multi{\calM}{}{}{}{A}$.
\\
(iii) If moreover there exists a biinvariant $\theta\in\G$ such
$\dd a=a\theta-\theta a$ for all $a\in A$, then $U(\gg_\G)\simeq B(\calL)$
as bialgebras,
where $\calL=(\Theta-\theta)^\bot$ is a braided Lie subalgebra of $\tgg_\G$.
\end{theorem}

\Proof $(i)$ By hypothesis, $\calI_\G\subset \ker\e_A$ satisfies
$\ad_L(\calI_\G)\subset A\otimes \calI_\G$ and $A\calI_\G\subset
\calI_\G$. By the above proposition, this implies that $\calI_\G$
is a 2-sided ideal of $\uA$, and a fortiori that
$\ad_L(\calI_\G)\subset A\otimes \calI_\G+\calI_\G\otimes A$.
Therefore $\calI_\G$ is an ideal of $(\uA,\um,\ueta,\ad_L)$ in the
sense of definition \ref{def-Lie-coalgebra}, and
$\widetilde{\G}_R\simeq A/\calI_\G$ is a counital braided Lie
coalgebra in $\multi{\calM}{}{}{A}{}$.
$(ii)$ According to $(i)$
and lemma \ref{lemma-rel-Lie-coLie}, $\tgg_\G$ is a braided Lie
algebra in $\calM^A$ with braided Lie bracket
$[\ ,\ ]=\Ad_L$ and coproduct $\uD$ adjoint
(in the conventions of lemma \ref{lemma-rel-Lie-coLie})
to the multiplication in $\uA/\calI_\G$.
$(iii)$ This is proposition
\ref{propo-triv-ext} and theorem \ref{theo-quadratic} put
together. \qed

\subsection{A ``quantum Lie functor''.}
Let $\calD$ be the category of bicovariant first order differential
calculi~: its objects are triples $(A,\G,\dd)$ where $A$ is a Hopf
algebra and $(\G,\dd)$ a bicovariant first order differential calculus
over $A$. Morphisms are pairs
$(\varphi^0,\varphi^1):(A,\G_A,\dd_A)\to (B,\G_B,\dd_B)$
such that $\varphi^0:A\to B$ is a Hopf algebra homomorphism
and $\varphi^1:(\G_A,\dd_A)\to (\G_B,\dd_B)$ is a morphism of Hopf bimodules
(over $A$) such that $\varphi^1\circ \dd_A=\dd_B\circ \varphi^0$~;
equivalently, such that $\varphi^1\circ \o_{R,A}=\o_{R,B}\circ \varphi^0$.
Because of the surjectivity axiom $\langle 3\rangle$ of a bicovariant FODC,
$\varphi^1$, if it exists, is uniquely determined by $\varphi^0$.
The condition of existence is easily seen to be that
$\varphi^0(\calI_{\G_A})\subset \calI_{\G_B}$.
Let $\calC\calQ\calT$ be the category
of co-quasitriangular Hopf algebras~: it consists of pair $(A,\rrr)$
where $A$ is a Hopf algebra and $\rrr$ is a co-quasitriangular structure
on $A$. Morphisms are Hopf algebra morphisms $\varphi:A\to B$ satisfying
$\rrr_B\circ \varphi =\rrr_A$.

\begin{proposition}\label{Lie-functor}
There is an exact functor $\bfL:\calC\calQ\calT\to \calD$, which
sends $(A,\rrr)$ to $(A,\G(\rrr),\dd)$ where $(\G(\rrr),\dd)$ is
the bicovariant FODC over $A$ whose associated left ideal is
$\calI(\rrr):= \ker\e_A\covdot \ker\e_A$ (the product in
$\uA=\uA(\rrr)$).
\end{proposition}

\Proof
$\uA$ is in particular a coalgebra in $\multi{\calM}{}{}{A}{}$,
\ie $\ad_L(a\covdot b)=a^{-1}b^{(-1)}\otimes a^{(0)}\covdot b^{(0)}$ for
all $a,b\in A$. Therefore,
 if $I$ and $J$ are any $\ad_L$-invariant left ideals of $\uA$,
their covariantized product $I\covdot J$ is also an
$\ad_L$-invariant ideal of $\uA$. But by proposition \ref{propo-A=Lie-coalg},
$I$, $J$ and $I\covdot J$ are also left ideals of $A$.
This holds in particular
for $I=J=\ker\e_A$.
Next one easily checks that if
$\varphi:(A,\rrr_A)\to (\varphi,\rrr_B)$ is a morphism in $\calC\calQ\calT$,
then the same map $\varphi$ is a morphism of $k$-algebras
$\uA\to \underline{B}$, and therefore satisfies
$\varphi(\ker\e_A\covdot\ker\e_A)\subset \ker\e_B\covdot \ker\e_B$.
This gives the functoriality property of $\bfL$. Finally, if $\varphi$ is
surjective ({\it resp.} injective), its restriction
$\ker\e_A\covdot\ker\e_A\to \ker\e_B\covdot \ker\e_B$ is also
surjective ({\it resp.} injective), proving exactness. Note that for
the associated quantum Lie algebras, if $\varphi:A\to B$ is surjective
({\it resp.} injective), it means that $\gg_B$ imbeds in
({\it resp.} maps onto)  $\gg_A$ as quantum Lie algebras in $\calM^A$.
\qed
\\

{\it Remarks~:}
If $(\calA,\e)$ is any augmented algebra, the space $(\ker\e/(\ker\e)^2)^*$
can be seen in either of the following ways~:
\begin{eqnarray*}
(\ker\e/(\ker\e)^2)^* &=&
    \{\chi\in\calA^\circ| \chi(1)=1, \chi((\ker\e)^2)=0\}
\\
    &=& \Prim_\e(\calA^\circ)
\\
    &\simeq& \Ext^1_\calA(k_\e,k_\e)
\end{eqnarray*}
where $\Prim_\e(\calA^\circ)=
    \{\chi\in \calA^\circ| \chi(ab)=\chi(a)\e(b)+\e(a)\chi(b)\}$
is the space of $\e$-primitive elements of the coalgebra
$\calA^\circ$, and $\Ext^1_\calA(k_\e,k_\e)$ parameterizes the
exact sequences $0\to k_\e\to M\to k_\e\to 0$ of $\calA$-modules
($k_\e$ is the 1-dimensional $\calA$-module afforded by
$\e:\calA\to k$, and the 2-dimensional module $M_\chi$ associated
to $\chi\in\Prim_\e(\calA^\circ)$ has basis $v_0$, $v_1$ such that
$av_0=\e(a)v_0+\chi(a)v_1$, $av_1=\e(a)v_1$. Note that if
$\chi\in\Prim_\e(\calA^\circ)$ and $a,b\in \calA$ satisfy
$\e(a)=1$, $\e(b)=0$ and $ab=qba$ for some $q\in k$, then one must
have $(1-q)\chi(b)=0$. Therefore, if $q\ne 1$, $\chi(b)=0$.
Alternatively, one has $(1-q)b=-qb(a-1)+(a-1)b\in (\ker\e)^2$,
therefore if $q\ne 1$,
 $b\in (\ker\e)^2$.
\cqfd
\\

Until the end of this paragraph, $\rrr$ is fixed and we write
$\G:= \G(\rrr)$ as above. Its tangent space $\gg_\G=
\Prim_\e((\uA)^\circ)$ is the space of primitive elements of
$(\uA)^\circ$. Note that if $A$ is commutative, $\uA=A$ (recall
that if $A=\calO(G)$ is the algebra of polynomial functions on
some algebraic group $G$, $\Prim(A^\circ)=\Lie(G)$ is the Lie
algebra of $G$ -if $G$ is finite, this is zero). In the general
case, it is not obvious to determine $\Prim_\e((\uA)^\circ)$, but
one can still describe some nice properties of $(\G,\dd)$. First,
$(\G,\dd)$ is clearly never inner (unless it is zero). Moreover,
$U(\gg_\G)$ is a Hopf algebra in $\calM^A$ with Hopf structure
$(\uD,\ue,\uS)$ uniquely determined by
\begin{equation}
\uD(x)=x\otimes 1+1\otimes x,\quad \ue(x),\quad \uS(x)=-x
\end{equation}
Indeed, the little coproduct $\d$ on $\gg_\G$
(see  (\ref{def-D-vs-delta})) is dual to the multiplication on
$\underline{\ker\e}_A/(\ker\e_A\covdot \ker\e_A)$ which is zero.
Thus, $\G$ has many properties of the standard differential calculus on
Lie groups. To complete the analogy, we show below that
the action $\rt$ of $A$ on $\G_R$ can be nicely ``linearized'',
and that the braiding $\s$
on  $\gg_\G$ in $\calC^A_A$ coincides with
the braiding in $\calM^A$, \ie the category in which $\gg_\G$
lives as a quantum Lie algebra (recall that they don't  in general).

Let $\qqq=\rrr_{21}*\rrr:A\otimes A\to k$. The maps
$\qqq_1,\qqq_2:A\to A^\circ$ are given by $\qqq_1=\rrr_2*\rrr_1$ and
$\qqq_2=\rrr_1*\rrr_2$. Recall from \cite{RS} that $(A,\rrr)$ is
called co-triangular if $\qqq=\e_A\otimes \e_A$ and co-factorizable
if $\qqq_1$ (equivalently $\qqq_2$) is injective.

\begin{lemma}\label{lemma-g-Lie}
(i) The factor crossed modules
$(\ker\e_A,\m,\ad_L)/\calI(\rrr)$ and
$(\ker\e_A,\brt_\rrr,\ad_L)/\calI(\rrr)$ are isomorphic.
Therefore, the left action of
$A$ on $(\G_R,\rt,\D_L)$ and its crossed module braiding $\s^t$ are
given by
\begin{eqnarray}&&a\rt \o_R(b) = \o_R(a\brt_\rrr b) = \rrr(b_{(1)},a_{(1)})\,
    \o_R(b_{(2)})\, \bar{\rrr}(b_{(3)},a_{(2)}),
    \label{commut-triangular}
\\
&&
\s^t\circ (\o_R\otimes \o_R)=(\o_R\otimes \o_R)\circ \Psi_{\rrr}
    \label{s=Psi-r}
\end{eqnarray}
where $\Psi_\rrr$ is defined in {\rm (\ref{Psi-r})}.
Moreover, for all $a\in A$, one has
\begin{equation}\label{ad-L}
\ad_L(a)\equiv -1\otimes a+\D(a)+ (\uS\otimes \id)\Psi_\rrr \, \D(a)
\bmod \calI(\rrr)\otimes A.
\end{equation}
In particular, for all $x,y\in \gg_\G$, one has
\begin{equation}
\langle a,[x,y]\rangle =
    \langle a_{(1)},x\rangle\,\langle a_{2)},y\rangle
    -\langle \Psi_\rrr(a_{(1)}\otimes a_{(2)}),x\otimes y\rangle
\end{equation}
(ii)
The $\Xi_{\rrr,\rrr}$-commutativity of $\uA$ implies that for all
$\chi\in\Prim_\e((\uA)^\circ)=\gg_\G$ and $a\in A$ one must have
\begin{equation}\label{cond-chi}
\qqq_1(a^{(-1)})\, \chi(a^{(0)})=\chi(a)\, 1_{A^\circ}
\end{equation}
(equivalently~:
$\ad_L(a)-1\otimes a\in \ker\qqq_1\otimes A+\ker\e_A\otimes \ker\chi$.)
\end{lemma}

\Proof First one obviously has, for all $a,b\in A$, $a\covdot
b\equiv \e(a)\, b+(a-\e(a)1)\e(b) \bmod(\ker\e_A\covdot\ker\e_A)$.
{}From this and the relation $ab=a_{(1)}\covdot (a_{(2)}\brt_\rrr
b)$, we get
\begin{equation}\label{ab-acovdotb}
a b\equiv a\brt_\rrr b+\e(b)\, (a-\e(a)1)\bmod (\ker\e_A\covdot\ker\e_A).
\end{equation}
In particular, if $\e(b)=0$,
$ab\equiv a\brt_\rrr b\bmod (\ker\e_A\covdot\ker\e_A)$,
hence the crossed module isomorphism as stated.
By the axiom $\langle 1'\rangle$ of a bicovariant FODC, one has
$a\rt\o_R(b)=\o_R(a(b-\e(b)1))=\o_R(a\brt_\rrr(b-\e(b)1))=\o_R(a\brt_\rrr b)$,
where we have used (\ref{ab-acovdotb}) for the second equality and
 $\o_R(a\brt_\rrr 1)=\e(a)\o_R(1)=0$ for the third.
(\ref{s=Psi-r}) follows. Finally, for $a\in A$, one has by
(\ref{ab-acovdotb}), $\ad_L(a)=a_{(1)}Sa_{(3)}\otimes
a_{(2)}\equiv (a_{(1)}-\e(a_{(1)}))\otimes
a_{(2)}+a_{(1)}\brt_\rrr Sa_{(3)}\otimes a_{(2)}$. Using
$\uS(a)=S(a_{(1)})\brt_\rrr \uS(a_{(2)})=\uS(Sa_{(1)}\brt_\rrr
a_{(2)})$, we get $a_{(1)}\brt_\rrr Sa_{(3)}\otimes a_{(2)}$ $=
\uS((a_{(1)} Sa_{(3)})\brt_\rrr a_{(4)}) \otimes a_{(2)}$ $=
(\uS\otimes \id)\Psi_\rrr\D(a)$. (We have used that the braided
antipode intertwines the coaction $\ad_L$, and therefore also the
action $\brt_\rrr$). For the next formula, we use that if
$x\in\gg_\G=\Prim_\e((\uA)^\circ)$, then $\langle 1,x\rangle=0$
and $\langle \uS(a),x\rangle=-\langle a,x\rangle$.
\\
(ii) By hypothesis, $\chi(a\covdot b)=\chi(a)\e(b)+\e(a)\chi(b)$
for all $a,b\in A$.
On the other hand, by the $\Xi_{\rrr,\rrr}$-commutativity of $\uA$,
one also has (with the notations of lemma \ref{lemma-new-cross})
$\chi(a\covdot b)=
\chi(b^{[0]}\covdot a^{(0)})\, \langle a^{(-1)},b^{[1]}\rangle$.
One then uses the identities $a^{(-1)}\, \e(a^{(0)})=\e(a)1$ and
$\e(b^{[0]})\, b^{[1]}=\qqq_2(b)$ to get
$\chi(a\covdot b)=\chi(b)\e(a)+\qqq(a^{(-1)},b)\chi(a^{(0)})$
for all $a,b\in A$. Comparing with the previous expression gives the claim.
\qed
\\

Note that in terms of $\dd(a)=\o_R(a_{(1)})\, a_{(2)}$,
(\ref{commut-triangular}) can be rewritten
\begin{equation}\label{d-triangular}
a\, \dd(b)=
    \rrr(b_{(1)},a_{(1)})\, \dd(b_{(2)})\,
     a_{(2)}\, \bar{\rrr}(b_{(3)},a_{(3)}).
\end{equation}

We apply the lemma above to describe $\G(\rrr)$ when $(A,\rrr)$
is of $GL(n)$ or $SL(n)$ type. In this case,
the condition $\m^{op}*\rrr=\rrr*\m$ for $A$
(ie $\um=\um\circ\Xi_{\rrr,\rrr}$ for $\uA$)
is almost the only set of relations.

Let $R=(R^i{}_j{}^k{}_l)$ be a  bi-invertible
solution of the Yang-Baxter equation on $k^n\otimes k^n$.
Recall that this means that the inverse $R^{-1}$, such
that $R^i{}_a{}^k{}_b (R^{-1})^a{}_j{}^b{}_l=\d^i_j\d^k_l$, exists,
as well as the second inverse $\tilde{R}$, such that
$\tilde{R}^i{}_b{}^a{}_j R^k{}_a{}^b{}_l=\d^i_l\d^k_j$.
Let $A(R)=A(C,r)$ be the FRT bialgebra \cite{FRT}
defined by $R$ (see lemma \ref{lemma-FRT}),
where  $C$ is the $n\times n$ matrix coalgebra
with basis ${t^i}_j$ ($\D{t^i}_j={t^i}_k\otimes {t^k}_j$)
and $r:C\otimes C\to k$ is the bilinear form such that
$r({t^i}_j,{t^k}_l)=R^i{}_j{}^k{}_l$. We assume that there exists a
central grouplike
$\det$ of degree $n$ (the quantum determinant) such that
(1) $\det$ is  not a zero divisor and there exists a (bi)algebra
automomorphism $f$ of $A(R)$ such that
$\det \,a=f(a)\, \det$ for all $a\in A(R)$,
(2) there exists a matrix of elements $({\tilde{t}^i}_j)$ of $A(R)$
such that
${t^i}_a{\tilde{t}^a}_j=\d^i_j\, \det$ for all $i,j$.
Then the localization of $A(R)$ at $\det$,
which we note $\calO_R(GL(n))$, has an antipode $S$ such that
$S({t^i}_j)={\tilde{t}^i}_j\, {\det}^{-1}={\det}^{-1}f({\tilde{t}^i}_j)$.
It is co-quasitriangular with $\rrr=z\, r$ when restricted to $C\otimes C$,
$z\in k^\times$ is arbitrary (in our setting $z$ plays no role since in general
the covariantized multiplication is invariant under the change
$\rrr\mapsto \ccc*\rrr$ for any central bicharacter $\ccc$ on $A$).
If $\det$ is central in $A(R)$ and if $(\det-1)$ belongs to the left and
right radicals of $\rrr$, we set
$\calO_q(SL(n))=A(R)/\langle \det-1\rangle$.
When $R=R_q$ is the  Drinfel'd-Jimbo $R$-matrix of type $A_n$, over $\dC$,
we write $\calO_q(G)$ instead of $\calO_R(G)$, $G=GL(n)$ or $G=SL(n)$,
and $\calO(G)=\calO_{q=1}(G)$.

The left quantum trace of $\calO_R(G)$, $G=GL(n)$ or $G=SL(n)$, is
the element $\utr={(\tilde{R})^{aj}}_{ia}\, {t^i}_j$. It is
$ad_L$-invariant~: $\ad_L(\utr)=1\otimes \utr$. We shall assume
that $\e(\utr)\ne 0$ (recall that $\e(\utr)= 0$ can happen; for
instance, for $\calO_q(G)$, $\utr=\sum_{i=1}^n q^{-2i}\, {t^i}_i$
has counit zero when $q^{2n}=1$; likewise, for $\calO(G)$,
$\utr=\tr$ has counit zero when $char(k)$ divides $n$).

Take first $A=\calO_R(G)$, $G=GL(n)$.
$\chi\in\gg_\G=\Prim_\e((\uA)^\circ)$ is uniquely determined by
its values on algebra generators of $\uA$, \ie on $C$ and
${\det}^{-1}$. {}From $\det\covdot {\det}^{-1}=1$, we get
$\chi({\det}^{-1})=-\chi(\det)$. Therefore $\Prim_\e((\uA)^\circ)$
can be identified with the space of linear functionals on $C$ such
that (\ref{cond-chi}) holds for all $a\in A$ (since $\um=\um\circ
\Xi_{\rrr,\rrr}$ are the only left relations to be checked).
{}From the assumption $\e(\utr)\ne 0$, $C$ decomposes as a direct
sum of sub-comodules (for $\ad_L$), $C=k\, \utr\oplus C^+$, where
$C^+=C\cap\ker\e$. So we have a vector space decomposition
$$
\gg_\G=k\, z\oplus \gg_\G^+
$$
where $\gg_\G^+=\{\chi\in\gg_\G| \chi(\utr)=0\}$ and
$k \,z=\{\chi\in\gg_\G:\chi(C^+)=0\}$.
Clearly, whatever $R$ is, $k z\ne 0$,
\ie there exists a non zero functional $z$ on $C$ such that
$z(C^+)=0$, $z(\utr)=1$ satisfying (\ref{cond-chi}).
One also easily checks that $\gg_\G^+$ is a quantum Lie subalgebra of
$\gg_\G$.
We would like to prove  (when it makes sense, \ie when
$\calO_R(GL(n))$ has a Hopf algebra quotient $\calO_R(SL(n))$),
that $\gg_\G^+$ is the quantum Lie algebra of $\calO_R(SL(n))$,
which we know imbeds into $\gg_\G$ by the exactness of the functor $\bfL$.
This follows from the lemma below~:

\begin{lemma}\label{lemma-det-tr}
Let $A=\calO_R(GL(n))$ as above.
Assume that $\e(\utr)\ne 0$, and that moreover
(1) $\utr$ is not a zero divisor in $\uA$, (2)
$C^+$ contains no $\ad_L$-invariant elements.
Then for all $\chi\in \Prim_\e((\uA)^\circ)$, one has
$$
\chi(\det)=\frac{n}{\e(\utr)}\; \chi(\utr).
$$
(In particular, $\chi(\det)=0$ if and only if $\chi(\utr)=0$).
\end{lemma}

\Proof
To prove this properly, one should take the general formula for $\det$
(see eg \cite{DMMZ})
 and reexpress $\det$ in terms of the covariantized product $\um$
of $\uA$. This is quite complicated and we use a trick that
requires the listed assumptions, which are probably not necessary for
a good proof.

For $a\in A$, let $a^{\underline{n}}$ be the $n$-th power of $a$
calculated in $\uA$. By iteration one checks that, for $\chi\in
\Prim_\e((\uA)^\circ)$,
$\chi(a^{\underline{n}})=n\,\chi(a)\,\e(a)^{n-1}$ holds. {}From
the decomposition $C=k\utr\oplus C^+$, from the fact that $\det\in
C^n=C^{\underline{n}}$, $\e(\det)=1$, and that $\utr$ is central
in $\uA$, there must exists elements $\a_i\in
(C^+)^{\underline{i}}$ such that
$$
\e(\utr)^n\, \det=(\utr)^{\underline{n}}+
\sum_{i=1}^n (\utr)^{\underline{n-i}}\covdot \a_i,
$$
\ie $\det$ is a polynomial in $\utr$ with coefficient $k$
for the leading one and in $C^+$ for the others.
Since $\chi$ vanishes on $C^+\covdot C^+\subset \ker\e\covdot \ker\e$,
we get
$\e(\utr)\chi(\det)=n\, \chi(\utr)+\chi(\a_1)$. We need to show that $\a_1=0$.
Since $\det$ is grouplike,
it is $\ad_L$-invariant. It is easy to see, using the property
of the covariantized product and the $\ad_L$-invariance of all powers
of $\utr$, that each of the terms of the above
decomposition of $\det$ must be $\ad_L$-invariant.
Let $\tau=(\utr)^{\underline{n-1}}$.
One must have $\ad_L(\tau\covdot\a_1)=1\otimes (\tau\covdot \a_1)=
(1\otimes \tau)\covdot (1\otimes \a_1)$. On the other hand,
$\ad_L(\tau\covdot \a_1)=\ad_L(\tau)\covdot \ad_L(\a_1)=
(1\otimes \tau)\covdot \ad_L(\a_1)$.
Since $1\otimes \tau$ is not a zero divisor by hypothesis,
this implies that $\a_1\in C^+$ is $\ad_L$-invariant, therefore zero by the
second hypothesis.
\qed
\\

The next proposition sums up the results for $\calO_R(G)$,
$G=SL(n)$ or $G=GL(n)$. We assume that $R$ has all the good
properties listed above~: $\calO_R(G)$ are well defined Hopf
algebras, $\e(\utr)\ne 0$, and  the statement of lemma
\ref{lemma-det-tr} hold, whenever its hypothesis are necessary or
not.

\begin{proposition}
(i) Assume $R_{21}{R}=1$.
For $G=GL(n)$, one has $\dim\gg_\G=n^2$, with basis
$\{{\chi_i}^j:i,j=1,...,n\}$ such that
$\langle {t^a}_b,{\chi_i}^j\rangle=\d_i^a\,\d_b^j$.
Let ${\o^a}_b=\o_R({t^a}_b)$ and $\o=({\o^a}_b)$.
The left crossed module structure of $\G_R$ and the corresponding $\s^t$
are given by
$$
\ttt_1\rt \mathbf{\o}_2=R_{21}\mathbf{\o}_{2}R,
\quad
\D_L({\o^i}_j)={t^i}_a\, S{t^b}_j\otimes {\o^a}_b,
\quad
\s^t(\o_1\otimes R_{21} \o_2 R)=R_{21}\o_2 R\otimes \o_1
$$
Equivalently, the first relation is
$R\dd(\ttt_1) \ttt_2=\ttt_2\, \dd(\ttt_1) R$.
The quantum Lie algebra structure of $\gg_\G$ is given by
$$
\s({\chi_i}^j\otimes {\chi_k}^l)=
    (\s_i{}^j{}_k{}^l)_b{}^a{}_d{}^c {\chi_a }^b \otimes {\chi_c }^d,
\quad
[{\chi_i}^j, {\chi_k}^l]= \d_k^j\, {\chi_i}^l-
    (\s_i{}^j{}_k{}^l)_r{}^a{}_b{}^r
    {\chi_a }^b
$$
where $(\s_i{}^j{}_k{}^l)_b{}^a{}_d{}^c=
    ({R_{21}})^j{}_\b{}^\a{}_b
    R^\ga{}_i{}^l{}_\a
    R^\b{}_d{}^a{}_\d
    ({R_{21}})^c{}_\ga{}^\d{}_k$.
The quantum Lie algebra for $SL(n)$ is the quantum Lie subalgebra
$\gg_\G^+=\{\chi\in\gg_\G|\langle \utr,\chi\rangle=0\}$, of dimension $n^2-1$.
\\
(ii) If $q$ is not a root of unity, the quantum Lie algebra
of $\calO_q(SL(n))$ -obtained by the functor $\bfL$- is zero.
\end{proposition}

\Proof
(i) If $R_{21}R=1$, one has $\rrr_{21}*\rrr=\e_A\otimes \e_A$,
therefore $\qqq_1(a)=\e(a) 1_{A^\circ}$ and
(\ref{cond-chi}) is trivially satisfied for all $a\in A$.
So, the fact that $\gg_\G\simeq C^*$
follows from the previous discussion and all formulas
follow from corresponding ones in lemma \ref{lemma-g-Lie}
(for the formula of $\s$, one can use that $R^{-1}=\tilde{R}=R_{21}$).
\\
(ii) For $A=\calO_q(SL(n))$, one has $\chi(\det)=\chi(1)=0$ if
$\chi\in\Prim_\e((\uA)^\circ)$, which is equivalent to
$\chi(\utr)=0$, \ie its quantum Lie algebra can be identified with
the space of functionals $\chi$ on $C$ satisfying $\chi(\utr)=0$
and (\ref{cond-chi}). Recall that when $q$ is not a root of unity,
$\calO_q(SL(n)$ is factorizable \cite{HS}, therefore $\qqq_1$ is
injective, and that $(C^+,\ad_L)$ is a simple comodule. The first
property and (\ref{cond-chi}) tell that, if $\chi\ne 0$,
$\ker\chi\cap C^+$ is a proper submodule of $(C^+,\ad_L)$, which
is simple. This is impossible therefore $\chi=0$. Note that these
arguments are also valid for $\calO_q(G)$, $G=SO(n)$ or $G=Sp(n)$,
$q$ not a root of unity, \ie their quantum Lie algebra is zero for
these as well. For $q$ a root of unity, the comodule $(C^+,\ad_L)$
remains simple, therefore the arguments are also valid
 provided $\qqq_1$ is injective on the coefficient subalgebra
of $(C^+,\ad_L)$ -the smallest subcoalgebra $T$ of $\calO_q(G)$ such that
$\ad_L(C^+)\subset T\otimes C^+$. We do not know when this is true.
\qed
\\

So to conclude,
for the standard deformation of $\calO(G)$, the quantum Lie functor
$\bfL$ gives uninteresting results. However, $\calO(G)$ has
some ``softer'' deformations, defined through
triangular $R$-matrices, which behaves better
(but which are less interesting in many other aspects).
Such examples do exist~: see eg \cite{JC}  and the references
cited there for $n=2$. (According to the classification of \cite{JC},
the differential calculus that we would
obtain for the ``Jordanian'' quantum group
$GL_{h,g}(2)$ via the functor $\bfL$ belongs to a 1-parameter family
a calculi that we cannot predict.)
Triangular $R$-matrices are known in higher dimension \cite{EH}
but, up to our knowledge, associated Hopf algebras have not been studied
yet.

\subsection{Example~: Finite groups.}
We illustrate the results of theorem \ref{main-theo}
with  the example of finite groups. The more interesting case of
quantum groups is considered in the next section.

Let $G$ be a finite group with unit element $e$, $H=kG$ and $A=k(G)$
  its dual with basis $\{f_g:g\in G\}$ such that $f_g(g')=\d_{g,g'}$.
$A$ is co-quasitriangular with $\rrr=\e_A\otimes \e_A$ so that the
braiding in $\multi{\calM}{}{}{}{A}\simeq \multi{\calM}{H}{}{}{}$
is the usual flip, and $\uA=A$, $\uH=H$. It is well-known that
bicovariant FODC over $A$ are in 1-1 correspondence with
$\Ad$-invariant subsets of $G$ not containing $e$, irreducible
calculi corresponding to conjugacy classes. Since $A$ is
semi-simple, they are all inner. For an $\Ad$-invariant subset
$\calC\subset G$, let $\theta_\calC=\sum_{g\in \calC} f_g$ and
$c_\calC=\sum_{g\in \calC} g$ ($\theta_\calC\in A$ is
$\ad_L$-invariant and $c_\calC\in kG$ is central). The calculus
$\G_\calC$ corresponding to $\calC$ has associated ideal
$\calI_\calC=\ker\e_A(1-\theta_\calC)$ with basis $\{f_g: g\ne e,
g\notin \calC\}$, and extended tangent space $\tgg_\calC=k\,
e\oplus c_\calC\leftharpoonup A$ with basis $\{X_g=g:g\in
\calC\cup \{e\}\}$. The braided Lie algebra structure of
$\tgg_\calC$ is given by and
$$
\Psi_\rrr(X_g\otimes X_h)=X_h\otimes X_g,\quad
\uD (X_g)=\D (X_g)=
    X_g\otimes X_g,\quad \e(X_g)=1, \quad [X_g,X_h]=X_{ghg^{-1}}.
$$
Therefore the canonical braiding on $\tgg_\calC$ is given by
$$
\Ups(X_g\otimes X_h)=X_{ghg^{-1}}\otimes X_g
$$
and $B(\tgg_\calC)$ is the (usual) bialgebra generated by
$\{X_g:g\in\calC\cup\{e\}\}$ with the above coproduct and relations
$X_g X_h=X_{ghg^{-1}} X_g$.
Note that  $\tgg_\calC=kX_e\oplus \calL_\calC$ is indeed the trivial extension
of the braided Lie subalgebra
$\calL_\calC:=  \{X\in\tgg_\calC|\langle 1-\chi_\calC,X\rangle=0\}$,
with basis $\{X_g:g\in \calC\}$.
The quantum Lie algebra $(\gg_\calC,\s,[\,,\,]\}$ of the differential calculus
has basis $\{x_g=g-e| g\in\calC\}$; its structure maps and the coalgebra
structure on $U(\gg_\calC)\simeq B(\tgg_\calC)/\langle X_e-1\rangle$
are given by
\begin{eqnarray*}
\quad
[x_g,x_h]=x_{ghg^{-1}}-x_h &,&
\D (x_g)=x_g\otimes 1+1\otimes x_g+x_g\otimes x_g\\
\s(x_g\otimes x_h)=x_{ghg^{-1}}\otimes x_g
&,& \e(x_g)=0.
\end{eqnarray*}
Thus, $U(\gg_\calC)\simeq B(\calL_\calC)$
is generated by $1$ and $x_g$, $g\in\calC$, with relations
$x_g\,x_h-x_{ghg^{-1}}x_g=x_{ghg^{-1}}-x_h$. It is not quadratic with respect
to the set of generators $\{x_g\}$, but it is with respect to the set
$\{X_g=1+x_g\}$.
\\

{\it Remark.} $U(\gg_\calC)\simeq B(\calL_\calC)$ has no antipode.
One could think that by localizing at
the multiplicative set generated by $\{X_g:g\in \calC\}$ one would
get a Hopf algebra with antipode $S(X_g)=(X_g)^{-1}$.
This turns out to be wrong because the elements $X_g$ can be zero divisors.
Example~: let $G=S_3$ with Coxeter generators $s_1,s_2$ and relations
$s_i^2=e$, $s_1s_2s_1=s_2s_1s_2$. Let $\calC$ be the class of transpositions
and let $X_1=X_{s_1}$, $X_2=X_{s_2}$, $X_3=X_{s_1s_2s_1}$ be the generators
of $B(\calL_\calC)$. The relations are~:
$$
X_iX_j=X_kX_i \qquad ((i,j,k) \;{\rm any\;permutation \; of}\; (1,2,3))
$$
Playing with these relations, one gets $X_iX_j^2=X_j^2X_i$ and
$X_iX_j^2=X_k^2X_i$, ($(i,j,k)$ all distinct). Therefore $X_i^2$
is central in $B(\calL_\calC)$ and $X_i(X_j^2-X_k^2)=0$. \cqfd
\section{Differential calculi and matrix braided Lie algebras on
$\calO_q(G)$}\label{sect-Jurco}
\setcounter{equation}{0}

In this section we apply the above general results for
co-quasitriangular Hopf algebras to the standard $q$-deformations
$\calO_q(G)$ and their variants. These are characterised by the
`quantum Killing form' $\qqq=\rrr_{21}*\rrr$ being nontrivial and
in this case there is a standard construction
\cite{Jurco}\cite{KS}\cite{Majid-classif-bicov} for their
bicovariant differential calculi going back to B.~Jurco in an
R-matrix setting. The corresponding braided Lie algebras in this
case are the matrix ones in \cite{Majid-braided-Lie}. For the
general treatment we allow $\qqq$ to be built in fact from pairs
of coquasitriangular structures.

\subsection{Construction of the calculi}

Let $(A,\rrr)$ be quasitriangular, $\rrr$ fixed. $A$ also has a
braided version in $\calM^A$, which is the one usually appearing
in the literature \cite{Majid-book}. We note it
$\uA(\rrr)^\rmright$ to distinguish it from
$\uA=\uA(\rrr)^\rmleft$ previously given.
 We note $\ad_R(a)=a^{(0)}\otimes a^{(1)}$.
Exactly as in \S\ref{sect-main-results},
given arbitrary other co-quasitriangular structure $\sss$ on $A$,
we have a left $A^\circ$-coaction $\l_{\rrr,\sss}:A\to A^\circ \otimes A$
given by $a\mapsto a^{[-1]}\otimes a^{[0]}=
\sss_2(a_{(1)})\rrr_1(a_{(3)})\otimes a_{(2)}$.
This gives a right crossed $A$-module $(A,\leftharpoonup_{\rrr,\sss},\ad_R)$
with braiding $\Xi_{\rrr,\sss}^\rmright$ where
\begin{eqnarray*}
a\leftharpoonup_{\rrr,\sss}b &=&
    }\langle b, \sss_2(a_{(1)})\rrr_1(a_{(3)})\rangle\,a_{(2),
\\
\Xi_{\rrr,\sss}^\rmright (a\otimes b)&=&
    b_{(2)}\otimes a_{(2)}\, \langle Sb_{(1)}\, b_{(3)},
    \sss_2(a_{(1)})\rrr_1(a_{(3)})\rangle.
\end{eqnarray*}
Note that $\Xi_{\rrr,\sss}^\rmright$ can also be written
\begin{equation}\label{Xi-commute}
\sss_{21}(a_{(1)},b_{(1)})\,\Xi_{\rrr,\sss}( X_{a_{(2)}}\otimes X_{b_{(3)}})\
     \rrr(a_{(3)},b_{(2)})
=
\sss_{21}(a_{(1)},b_{(2)})\, X_{b_{(1)}}\otimes X_{a_{(2)}}\,
     \rrr(a_{(3)},b_{(3)}).
\end{equation}
When $\sss=\bar{\rrr}_{21}$, we write
$a\blt_\rrr b=a\leftharpoonup_{\rrr,\bar{\sss}_{21}}b=
b_{(2)}\,\langle Sa_{(1)}a_{(3)},\rrr_2(b)\rangle$ and
$\Psi_\rrr^{\rmright} =\Xi_{\rrr,\bar{\rrr}_{21}}^{\rmright}$,
the braiding on $(A,\ad_R)$ in $\calM^A$ thanks to $\rrr$.
The multiplication in $\uA(\rrr)^\rmright$ is
$a\covdot b=(a\blt_\rrr Sb_{(1)})b_{(2)}$, and $\uA(\rrr)^\rmright$ is
$\Xi_{\rrr,\rrr}^{\rmright}$-commutative.
We let
$$
\qqq=\sss_{21}*\rrr.
$$
It still satisfies $\qqq*\m=\m*\qqq$. One has
$\qqq_1=\sss_2 * \rrr_1$, $\qqq_2=\sss_1 * \rrr_2$,
and by straightforward applications of the properties of $\rrr$
and $\sss$, one obtains  for all $a,b\in A$,
\begin{eqnarray}
\qqq_1(ab)=
    \sss_2(b_{(1)})\, \qqq_1(a)\, \rrr_1(b_{(2)}),
    &&
\D (\qqq_1(a))=
    \sss_2(a_{(1)})\, \rrr_1(a_{(3)})\otimes \qqq_1(a_{(2)}).
    \label{eq-Q1}
\\
\qqq_2(ab)=
    \sss_1(a_{(1)})\, \qqq_2(b)\, \rrr_2(a_{(2)}),
    &&
\D (\qqq_2(a))=
    \qqq_2(a_{(2)})\otimes\sss_1(a_{(1)})\, \rrr_2(a_{(3)}).
    \label{eq-Q2}
\end{eqnarray}
In addition, $\e(\qqq_i(a))=\e(a)$, $\qqq_i(1)=1$, ($i=1,2$).
The following is well-known when $\sss=\rrr$.

\begin{lemma}\label{propofq}
{\rm (Intertwining properties of $\qqq$).}
\\
(i) $\qqq_1$ intertwines the left adjoint and coadjoint actions of
$A^\circ$, and $\qqq_2$ the right ones.
\\
(ii) $\qqq_1:\uA(\sss)^\rmright\to A^\circ$ and
$\qqq_2:\uA(\rrr)^\rmleft\to A^\circ$ are homomorphisms of algebras.
In particular, $\im \qqq_1$ is a subalgebra of $A^\circ$,
and a left coideal by {\rm (\ref{eq-Q1})}.
Moreover, $\qqq_1(A^{\ad_R})$ and $\qqq_2(A^{\ad_L})$
belong to the center of $A^\circ$.
\end{lemma}

\Proof
(i) This is lemma \ref{intertwinner} applied to $\xi=\qqq$.
(ii) We prove it for $\qqq_2$ and $\uA(\rrr)^\rmleft$.
For $a,b\in A$, using repeatedly
that $\rrr_2$ is an antialgebra map
and a coalgebra map, we get
\begin{eqnarray*}
\qqq_2(a\covdot b) &=&
    \qqq_2(a_{(1)}\, b_{(2)})\,
    \langle b_{(1)}S(b_{(3)}), \rrr_2 S(a_{(2)})\rangle
\\
    & \stackrel{(\ref{eq-Q2})}{=} &
    \langle b_{(1)}, \rrr_2 S(a_{(4)})\rangle\,
    \sss_1(a_{(1)})\,
    \underline{\qqq_2(b_{(2)})\, \rrr_2(a_{(2)})\,
    \langle S(b_{(3)}), \rrr_2 S(a_{(3)})\rangle}
\\
    &=&
    \langle b_{(1)}, \rrr_2 S(a_{(4)})\rangle\,
    \sss_1(a_{(1)})\,
    \rrr_2(a_{(2)})\,
    \langle b_{(2)}, \rrr_2 (a_{(3)})\rangle\,
    \qqq_2(b_{(3)})
\\
    &=&
    \langle b_{(1)}, \rrr_2 S(a_{(3)})\, \rrr_2(a_{(2)})\rangle\,
    \qqq_2(a_{(1)})\, \qqq_2(b_{(2)})=\qqq_2(a)\, \qqq_2(b)
\end{eqnarray*}
(The underlined term is transformed using $S\rrr_2S=\rrr_2$
and $\m*\qqq=\qqq*\m$).
Finally,  if $\ad_L(a)=1\otimes a$,
then for all $h\in A^\circ$,
$Ad_R^* h(a)=\langle 1,h\rangle \, a$, therefore
$\Ad_R h(\qqq_2(a))=\e(h)\, \qqq_2(a)$ by (i), \ie $\qqq_2(a)$ is central
in $A^\circ$. For $\qqq_1$, the proof is analogous.
\qed

\begin{proposition}\label{Jurco}
Let $C_1$ be a subcoalgebra of $A$ containing $1_A$. Then
$\tgg(C_1,\qqq):=\qqq_1(C_1)$ is the extended
tangent space of a bicovariant FODC over $(A,\rrr)$,
with associated left ideal
\begin{equation}\label{I-C-Q}
\calI(C_1,\qqq)=
    \{a\in A|\forall c\in C_1, \qqq(c,a)=0\}\; \supset \ker\qqq_2.
\end{equation}
(ii) Conversely, let $(\G,d)$ be a bicovariant FODC over $A$
with associated left ideal $\calI_\G$.
If $\qqq_1$ is injective, then
$$
C_1=\{a\in A| \; \qqq(a,\calI_\G)=0\}
$$
is a subcoalgebra of $A$ containing $1_A$. If moreover
$\tgg_\G\subset \im\qqq_1$, then $\tgg_\G=\qqq_1(C_1)$.
\\
(iii) Assume $C_1=k\, 1_A\oplus C$ for some subcoalgebra $C$
and that $\qqq_1$ is injective on $C$.
Then $\tgg_\G=1_{A^\circ}\oplus \calL$
where $\calL=\qqq_1(C)$ is a braided Lie subalgebra. The associated
bicovariant FODC is inner,
with tangent space $\gg_\G=\{x-\e(x) 1: x\in \calL\}$ and
$U(\gg_\G)\simeq B(\calL)$ as (quadratic) bialgebras.
\end{proposition}

Note that if we take $C_1=A$, we get $\calI(C_1,\qqq)=\ker\qqq_2$
which is indeed a left crossed submodule of
$(A,\m,\ad_L)$. It is the
smallest ideal we can divide by through this construction,
thus we should assume that $\qqq\ne \e_A\otimes \e_A$.
\\

\Proof $(i)$ We check the properties $(a)$, $(b)$ and $(c)$ of
lemma \ref{prop-ext-T} for $\tgg_\G=\qqq_1(C_1)$. Since $1_A\in
C_1$, one has $1_{A^\circ}\in \qqq_1(C_1)$. Since $C_1$ is a
subcomodule for $\ad_L$, it is also a left submodule the left
adjoint coaction of $A^\circ$ on $A$, therefore by the
intertwining property of $\qqq_1$~: $\Ad_L
h(\qqq_1(c))=\qqq_1(\Ad_L^*h (c))\in \qqq_1(C_1)$ for any $h\in
A^\circ$ and $c\in C_1$. Finally, $\qqq_1(C_1)$ is a left coideal
by the left equality in (\ref{eq-Q1}).
\\
$(ii)$
For $a\in A$, $b\in \calI_\G$ and $c\in C_1$, one has
since $\calI_\G$ is a left sided ideal of $\uA(\rrr)^\rmleft$,
$0=\langle c,\qqq_2(a\covdot b)\rangle =
\langle c,\qqq_2(a)\, \qqq_2(b)\rangle=
\langle \D(c),\qqq_2(a)\otimes \qqq_2(b)\rangle$
(where $\covdot$ is the multiplication in $\uA(\rrr)^\rmleft$).
If $\qqq_2$ is injective, $\qqq_2(A)$ separates the elements of $A$,
therefore $\D(C_1)\subset A\otimes C_1$. Since $\calI_\G$ is also a right ideal
of $\uA(\rrr)^\rmleft$ (proposition \ref{propo-A=Lie-coalg}),
one obtains likewise $\D(C_1)\subset C_1\otimes A$,
so $\D(C_1)\subset C_1\otimes C_1$. By definition,
$\tgg_\G=\{x\in A^\circ| \langle \calI_\G,x\rangle =0\}$ therefore
if $\tgg_\G\subset \im\qqq_1$, we immediately get $\tgg_\G=\qqq_1(C_1)$
by definition of $C_1$.
\\
$(iii)$ $\G$ is inner by lemma \ref{lemma-inner} and we apply theorem
\ref{main-theo}.
\qed
\\

\Remark
The above proposition is essentially well-known,
but is slightly more general than analogous results in
\cite{KS} \cite{Majid-classif-bicov} because it can describe differential
calculi which are not inner (the coalgebra imbedding $k1_A\hookrightarrow C_1$
might be non split in the non semi-simple case, for instance for
quantum groups at roots of unity).
It is shown however in \cite{HS} that
for the standard quantum groups $\calO_q(G)$, $G=SL(n)$ or $G=Sp(n)$,
$q$ not a root of unity, any
bicovariant FODC arises in this way for a uniquely determined
subcoalgebra $C_1$ and for some pair of co-quasitriangular structures
$(\rrr,\sss)$, where $\sss=\ccc*\rrr$, $\ccc$ a central bicharacter.
\cqfd
\\

The next proposition  describes the braided and quantum Lie algebras
$\tgg_\G$ and $\gg_\G$ associated to $\qqq$ and $C_1$.
We assume that $\qqq_1$ is injective on $C_1$, so that $\tgg_\G$ can be
identified with $C_1$. (The given formulas are actually true in the non
injective case, but should be considered with care).

\begin{proposition}
Let $C_1$ be a subcoalgebra of $A$ containing $1_A$ as before,
and $\tgg_\G=\qqq_1(C_1)$.
\\
$(i)$ For $c\in C_1$, let $X_c=\qqq_1(c)\in\tgg_\G$. The right
crossed module structure of $\tgg_\G$ over $A$ is given by (we
write $c^{(0)}\otimes c^{(1)}:= \ad_R(c)$)~:
\begin{eqnarray*}
X_c\leftharpoonup a &=&
    X_{c_{(2)}}\,
    \langle a, \sss_2(c_{(1)})\, \rrr_1(c_{(3)})\rangle
    =: X_{c\leftharpoonup_{\rrr,\sss} a}
\\
\d_R(X_c) &=& X_{c_{(2)}}\otimes S(c_{(1)})\, c_{(3)}=
    X_{c^{(0)}}\otimes c^{(1)}.
\end{eqnarray*}
The corresponding braidings
$\Psi=\Psi_\rrr$ (in $\multi{\calM}{}{}{}{A}$)
and
$\Ups=\tilde{\s}$ (in $\multi{\calC}{}{A}{}{A}$)
on $\tgg_\G$ are~:
\begin{eqnarray*}
\Psi(X_a\otimes X_b) &=&
         X_{b_{(2)}}\otimes X_{a_{(2)}}\,
        \langle S(b_{(1)})\, b_{(3)},
            \rrr_1S(a_{(1)})\,\rrr_1(a_{(3)})\rangle
\\
\Ups(X_a\otimes X_b) &=&
         X_{b_{(2)}}\otimes X_{a_{(2)}}\,
        \langle S(b_{(1)})\, b_{(3)},
            \sss_2(a_{(1)})\,\rrr_1(a_{(3)})\rangle
\end{eqnarray*}
The braided Lie algebra structure of $\tgg_\G$ in
$\multi{\calM}{}{}{}{A}$ is given by~:
$$
\uD(X_c) = X_{c_{(1)}}\otimes X_{c_{(2)}},\qquad
\ue(X_c) = \e(c),\qquad
{[}X_a,X_b{]} = X_{b^{(0)}}\,\qqq(a,b^{(1)}).
$$
Let $Y_c$ be the image of $X_c\in \tgg_\G\subset A^\circ$ in $B(\tgg_\G)$.
The defining relations of $B(\tgg_\G)$ can be written~:
\begin{equation}\label{rel-B-tgg}
\sss_{21}(a_{(1)},b_{(1)})\, Y_{a_{(2)}}\, Y_{b_{(3)}} \rrr(a_{(3)},b_{(2)})
=
\sss_{21}(a_{(1)},b_{(2)})\, Y_{b_{(1)}}\, Y_{a_{(2)}} \rrr(a_{(3)},b_{(3)}).
\end{equation}
(ii)
For $c\in C$, let $x_c=\qqq_1(c)-\e(c)1_{A^\circ}\in \gg_\G$.
The braiding $\Psi_\rrr$ (in $\calM^A$) and $\s$ (in $\calC^A_A$)
on $\gg_\G$ are given by the same formulas as above,
with $X$ replaced by $x$;
the quantum Lie bracket on $\gg_\G$ and the
braided coproduct on $U(\gg_\G)$ are given by
$$
\uD(x_c) = x_c\otimes 1+1\otimes x_c+x_{c_{(\1)}}\otimes x_{c_{(\2)}}
,\qquad
{[}x_a,x_b{]} =\qqq(a,b^{(1)})\, x_{b^{(0)}}-\e(a)\; x_b.
$$
Moreover, $U(\gg_\G)$ can also be seen as the algebra generated
by the elements $Y_c$ with relations {\rm (\ref{rel-B-tgg})} and $Y_1=1$.
\end{proposition}

\Proof This is a direct application of the definitions and the
intertwining properties of $\qqq$. The right coaction of $A$ on
$\tgg_\G$ is given by (proposition \ref{prop-ext-T})
$X_c\leftharpoonup a=$ $\langle (\qqq_1(c))_{(1)},a\rangle\,
(\qqq_1(c))_{(2)}$ $\stackrel{(\ref{eq-Q2})}{=}$ $\langle
a,\sss_2(c_{(1)})\rrr_1(c_{(3)})\rangle \, \qqq_1(c_{(2)})$.
Likewise, the coproduct $\uD$ on $\tgg_\G$ is such that $\langle
a\otimes b,\uD \qqq_1(c)\rangle=$ $\langle a\covdot
b,\qqq_1(c)\rangle=$ $\langle \qqq_2(a\covdot b), c\rangle=$
$\langle \qqq_2(a)\, \qqq_2(b), c\rangle=$ $\langle a\otimes b,
(\qqq_1\otimes \qqq_1)\D(c)\rangle$, where $\covdot$ is the
product in $\uA(\rrr)^\rmleft$. Finally, the braided Lie bracket
is $[X_a,X_b]=\Ad_L \qqq_1(a)(\qqq_1(b))$
$=\qqq_1(\Ad_L^*\qqq_1(a)(b))=\qqq_1(b^{(0)})\, \langle
b^{(1)},\qqq_1(a)\rangle$. (It can also be obtained from
$[X_a,X_b]=(\id\otimes \e)\ts(X_a\otimes X_b)$). The braidings are
obtained directly from their definition~: $\Psi_\rrr(X_a\otimes
X_b)=X_{b^{(0)}}\otimes X_{b^{(0)}} \, \rrr(a^{(1)},b^{(1)})$ and
$\ts(X_a\otimes X_b)=X_{b^{(0)}}\otimes X_a\leftharpoonup
b^{(1)}$. Note that $\Psi_\rrr$ and $\Ups=\tilde{\s}$ almost
coincide on $\tgg_\G=\qqq_1(C)$. They are equal when
$\sss=\bar{\rrr}_{21}$, but in this case, $\qqq_1(C_1)=k\,
1_{A^\circ}$, \ie the calculus is trivial. The formulas for
$\Psi_\rrr$, $\ts$ and the braided Lie bracket could be further
developed. For instance, one has
\begin{eqnarray}
{[}X_a,X_b{]} &=& X_{b_{(2)}}\, \langle a,\qqq_2(Sb_{(1)}\,b_{(3)})\rangle\,
    =\, X_{b_{(3)}} \, \sss(Sb_{(2)}, a_{(1)})\, \qqq(a_{(2)}, b_{(4)})\,
    \rrr(a_{(3)}, Sb_{(1)}) \nonumber
    \\
&=&  X_{b_{(3)}} \, \bar{\sss}(b_{(2)}, a_{(1)})\,
    \sss(b_{(4)},a_{(2)})\, \rrr(a_{(3)},b_{(5)})\,
    \rrr(a_{(4)}, Sb_{(1)})\label{struc-constants}
\end{eqnarray}
$(ii)$ is clear from $(i)$.
\qed
\\

\Remark
$(i)$ Let $\ccc$ be a central bicharacter on $A$ \cite{Schmudgen}.
If we take $\sss=\ccc*\bar{\rrr}_{21}$, we get $\qqq=\ccc_{21}$
which can hardly be non degenerate.
For $\sss=\sss*\rrr$, the right action of $A$ on $\tgg_\G=\qqq_1(C_1)$,
and therefore also the left action on $\widetilde{\G}_R$,
depends on $\ccc$, but all the remaining defining structure maps
of $\tgg_\G$, \ie $\Psi$, $\Ups$, $[\,,\,]$ and $\uD$, do not
as can be easily checked. Therefore $\ccc$ controls how
$\tgg_\G$ sits inside $A^\circ$, and distinct $\ccc$'s
give non isomorphic calculi, but the corresponding extended tangent spaces
$\tgg_\G=\qqq_1(C_1)$ are isomorphic as abstract braided Lie algebras.
In the following, we focus only on these, therefore
we assume that $\sss=\rrr$ and write $\uA^\rmright=\uA(\rrr)^\rmright$.
\\
$(ii)$
When $1\hookrightarrow C_1$ splits, which is the case we are interested in,
one has $\tgg_\G=k1_{A^\circ}\oplus \calL$, $\calL=\qqq_1(C)$.
Then the three spaces $\gg_\G$, $\calL$ and $C$ and be identified,
and $U(\gg_\G)\simeq B(\calL)$ is generated by $Y_c$, $c\in C$,
with relations (\ref{rel-B-tgg}).
According to the remarks following (\ref{bra-vs-m}),
there is an algebra homomorphism $U(\gg_\G)\simeq B(\calL)\to A^\circ$ such
that $Y_c\mapsto X_c$. If $\sss=\rrr$ (or more generally if
$\sss=\ccc*\rrr$), we see comparing (\ref{rel-B-tgg})
and the $\Xi_{\rrr,\rrr}^\rmright$-commutativity of $\uA^\rmright$
(\ref{Xi-commute})
that this homomorphism
factors through a homomorphism
$B(\calL)\to \uA^\rmright$,
that is, it is the composition~:
\begin{equation}\label{can-hom}
\begin{array}{rcccc}
U(\gg_\G)\simeq B(\calL) & \longrightarrow & \uA^\rmright &
    \stackrel{\qqq_1}{\longrightarrow} & A^\circ
\\
Y_c &\mapsto & c &\mapsto &\qqq_1(c)=X_c
\end{array}
\end{equation}
$(iii)$ If moreover $C$ is simple, let $\{{\l^i}_j:i,j=1,...,n\}$
be a basis such that $\D{\l^i}_j={\l^i}_k\otimes {\l^k}_j$, let
$\{{X^i}_j=X_{{\l^i}_j}=\qqq_1({\l^i}_j):i,j=1,...,n\}$ be the
corresponding basis of $\calL$, and ${Y^i}_j=Y_{{\l^i}_j}$ the
image of ${X^i}_j$ in $B(\calL)$. Note that one has $\uD{X^i}_j
={X^i}_k\otimes {X^k}_j$ (coproduct in the braided Lie algebra
$\calL$) and $\D {X^i}_j=\rrr_2({t^i}_a)\, \rrr_1({t^b}_j)\otimes
{X^a}_b$ (coproduct in $A^\circ$). One defines the tensors
$R^i{}_j{}^k{}_l=\rrr({\l^i}_j,{\l^k}_l)$,
$(R^{-1})^i{}_j{}^k{}_l=\rrr(S{\l^i}_j,{\l^k}_l)$,
$\tilde{R}^i{}_j{}^k{}_l=\rrr({\l^i}_j,S{\l^k}_l)$, and $Q=R_{21}R$.
Then in the numerical suffix notation\cite{FRT}, the structure
maps $\Psi$, $\Ups$ and $[\,,\,]$ of $\calL$ are given by
\begin{eqnarray*}
\Psi(R_{21}X_1 R_{21}^{-1}\otimes  X_2) &=& X_2\otimes  R_{21}X_1R_{21}^{-1}
\\
\Ups(R_{21}X_1 R\otimes  X_2) &=& X_2\otimes  R_{21}X_1R,
\\
{[}R_{21}X_1 R, X_2] &=& X_2\, Q.
\end{eqnarray*}
as in \cite{Majid-braided-Lie}\cite{Majid-sol-YBE}, i.e. we obtain
the `matrix braided Lie algebras' introduced there. Here the
algebra $B(\calL)$ is abstractly generated by the $n^2$ elements
${Y^i}_j$ with relations
$$
(R_{21}Y_1 R)\, Y_2=Y_2\, (R_{21}Y_1R)
$$
and coincides with the bialgebra of braided matrices $B(R)$ (a
bialgebra in the category of right $A(R)$-comodules
\cite{Majid-book}). Finally, the quantum Lie bracket on $\gg_\G$
is ${[}R_{21}x_1 R, x_2] = x_2\, Q - Q x_2$. \cqfd

\subsection{Example: $\calO_q(SL(n))$}\label{SLn}
Let $A=\calO_q(SL(n))$, $C$ its fundamental subcoalgebra,
with basis $\{{t^i}_j:i,j=1,...,n\}$, $R$ its standard $R$-matrix,
$\rrr$ the unique co-quasitriangular structure on $A$ such
$\rrr({t^i}_j,{t^k}_l)=R^i{}_j{}^k{}_l$,
$\det$ the quantum determinant.
Recall that $A$ (resp. $\uA^\rmright$)
is the quotient of $A(R)$ (resp. $B(R)$) by the two-sided ideal
generated by the central element $\det-1$.
Here $A(R)$  and $B(R)$ are the algebras generated by the matrix $\ttt$
of elements ${t^i}_j$ and relations
$$
A(R)~:\;\;R\ttt_1\ttt_2=\ttt_2\ttt_1R,\qquad
B(R)~:\;\;(R_{21}\ttt_1 R)\covdot \ttt_2=\ttt_2\covdot (R_{21}\ttt_1R).
$$
Consider the standard $n^2$-dimensional bicovariant FODC $\G$
over $A$ corresponding to the subcoalgebra $C$.
One has $U(\gg_\G)\simeq B(R)$ by the previous section, and therefore
there exists a central grouplike element (also written $\det$)
inside $U(\gg_\G)$ such that
$U(\gg_\G)/\langle \det-1\rangle\simeq \uA^\rmright$.
This gives the kernel of the first map in (\ref{can-hom}),
for all values of $q$.
If $q$ is not root of unity, the second is injective \cite{RS}
\cite{HS} \cite{BS}.
Its image is described in \cite{BS} proposition 5.
With the definition of $\calU_q(sl(n))$ given in \cite{BS},
$\qqq_1(A)=\calF_\ell(\calU_q(sl(n))$ is the locally finite part of
the left adjoint $\calU_q(sl(n))$-module.
Therefore we have~:
\begin{proposition}
Let $A=\calO_q(SL(n))=A(R)/\langle \det-1\rangle$, $q$ not a root
of unity, and $\G$ the standard $n^2$-dimensional bicovariant FODC
over $A$. Then $U(\gg_\G)\simeq B(R)$, and
$B_q(SL(n))=B(R)/\langle \det-1\rangle$ $\simeq
\uA^{\rmright}\simeq \calF_\ell(\calU_q(sl(n))$ is a Hopf algebra
in $\calM^A$.
\end{proposition}

This makes more precise the sense in which braided Lie algebras
solve the `Lie algebra problem' for quantum groups in
\cite{Majid-braided-Lie}. It can also be viewed as the
self-duality of the braided versions of quantum groups, i.e. the
above quotient $B_q(SL(n))$ is isomorphic to a (braided) version
of $U_q(sl(n))$ (with algebra the locally finite part) via the
quantum Killing form\cite{Majid-book}.

The $n^2$-dimensional braided Lie algebra in this example can be
denoted $\widetilde{sl_{q}(n)}$ and the $n=2$ case is computed
explicitly in \cite[Ex. 5.5]{Majid-braided-Lie}. The enveloping
algebra $B(R)=BM_q(2)$ is the standard $2\times 2$ braided matrix
algebra \cite{Majid-book}. The structural form of the general
$B(R)=BM_q(n)$ and their homological properties appear in
\cite{LeB-94}. In particular, for generic $q$ it is known that
they have the same Hilbert series as polynomials in $n^2$
variables.

We can give the the relations of $B(R)$ more explicitly as
follows, in fact for the full multiparameter $SL(n)$-type family.
In our conventions (which are slightly different from
\cite{LeB-94}) the R-matrix is
\[
R^i{}_k{}^j{}_l=\delta^i{}_k\delta^j{}_lM_{ij}+\delta^i{}_l\delta^j{}_k
L_{ij},\quad
M_{ij}=q\delta_{ij}+\theta_{ji}\frac{q}{r_{ij}}+\theta_{ij}\frac{r_{ji}}{q};\quad
L_{ij}=\theta_{ji}(q-q^{-1}),\] where $\theta_{ij}$ denotes the
function which is $1$ {\em iff} $i>j$ and otherwise zero, and
$r_{ij}\ne 0$ are multiparameters defined for $i<j$ and
constrained by $\prod_{i<j}{r_{ij}\over q}=\prod_{i>j}{r_{ji}\over q}$ for all $j$ as explained in \cite{LeB-94}. The standard $\calO_q(SL(n))$ case is $r_{ij}=q$. We let $\mu=q-q^{-1}$. Let us also introduce the `cocycle' defined
for $i,j,k$ all distinct by
\[ \sigma_{ijk}=\left({q\, r_{\sigma(k),\sigma(i)}\over
r_{\sigma(j),\sigma(i)}\,
r_{\sigma(k),\sigma(j)}}\right)^{(-1)^{l(\sigma)}}\] where
$\sigma\in S_3$ is the unique permutation of $i,j,k$ such that
$\sigma(i)>\sigma(j)>\sigma(k)$ and $l(\sigma)$ is its length. By
convention, $\sigma_{ijk}=0$ if the $i,j,k$ are not distinct.
Finally, in order to make computations we need the matrices for
$R^{-1}$ and $\tilde{R}$. Using that $R$ is q-Hecke, one can show
that
\[ R^{-1}(q,\{r_{ij}\})=R(q^{-1},\{ r^{-1}_{ij}\}),\quad
\tilde{R}^i{}_k{}^j{}_l=R^{-1}{}^i{}_k{}^j{}_lq^{2(l-j)}\] which
means that they have the same form as the above with $M_{ij}^{-1}$
in place of $M_{ij}$ and $-L_{ij}$ or $-L_{ij}q^{2(i-j)}$
respectively in place of $L_{ij}$. Let us denote by
$\bar\sigma_{ijk}$ the same expression as $\sigma_{ijk}$ but with
$q,r_{ij}$ inverted.

\begin{lemma}
\[
M^{-1}_{ki}M_{ji}M^{-1}_{jk}
=\sigma_{ijk}+q\delta_{ij}+q^{-1}(\delta_{ik}+\delta_{jk})
-(q+q^{-1})\delta_{ij}\delta_{jk}\]
\[
M^{-1}_{ki}M^{-1}_{jk}M_{li}M_{jl}=\sigma_{ijk}\bar\sigma_{ijl}
+q^{-1}\bar\sigma_{ijl}(\delta_{ki}+\delta_{kj})
+q\sigma_{ijk}(\delta_{li}+\delta_{lj})\quad\quad\]
\[ \quad\quad\quad\quad\quad\quad\quad\quad\quad\quad\quad+\delta_{ij}(1+\mu
q\delta_{il}-q^{-1}\mu\delta_{jk})+\delta_{ik}\delta_{jl}+\delta_{kj}\delta_{il}\]
if $k\ne l$ and $1$ if $k=l$.
\end{lemma}

Finally, we obtain $\Ups$ from the explicit R-matrix formula for
it in \cite{Majid-braided-Lie}, namely
\[ \Ups(X^i{}_j\otimes
X^k{}_l)=R^{-1}{}^a{}_m{}^i{}_bR^n{}_c{}^b{}_o R^p{}_d{}^c{}_l
\tilde R^d {}_j{}^k{}_a X^m{}_n\otimes X^o{}_p.\] This determines
the relations of $B(R)$ as $\Ups$-commutative (one can also work
from the `reflection' form of the $B(R)$ relations as in
\cite{LeB-94} if one does not need $\Ups$ explicitly). Since we
are dealing with an abstract braided Lie algebra, we do not
distinguish between the $n^2$ basis elements ${X^i}_j$ of $\calL$
and their images in $B(\calL)$ as was done for instance in
(\ref{can-hom}).

\begin{proposition} The $\left({n^2\atop 2}\right)$
relations of the multiparameter $B(\widetilde{sl_{n,q}})=BM_{q}(n)$
may be listed for distinct $i,j,k,l$ as follows.

(i) For $i<k$:\quad\quad $X^i{}_i X^k{}_k=X^k{}_k X^i{}_i$.

(ii.a) For $k,l>i$:\quad\quad $X^i{}_i X^k{}_l=X^k{}_lX^i{}_i$.

(ii.b) For $i>k,l$:\quad\quad $X^i{}_i X^k{}_l=X^k{}_lX^i{}_i$.

(ii.c) For $l>i>k$:\quad\quad $X^i{}_i
X^k{}_l-X^k{}_lX^i{}_i=-\mu\bar\sigma_{kil}X^i{}_lX^k{}_i $.

(ii.d) For $k>i>l$:\quad\quad $X^i{}_i
X^k{}_l-X^k{}_lX^i{}_i=\mu\bar\sigma_{ilk}X^i{}_lX^k{}_i$.

(ii.e) For $i<l$:\quad\quad $X^i{}_i X^i{}_l-q^{-2}X^i{}_lX^i{}_i=
-q^{-1}\mu \sum_{a<i}X^i{}_aX^a{}_l$.

(ii.f) For $i>l$:\quad\quad
$X^i{}_iX^i{}_l-X^i{}_lX^i{}_i=-q^{-1}\mu
\sum_{a<i}X^i{}_aX^a{}_l$.

(ii.g) For $i>k$:\quad\quad $X^i{}_iX^k{}_i-X^k{}_iX^i{}_i=\mu
q^{-1}\sum_{a<i} X^k{}_a X^a{}_i $.

(ii.h) For $i<k$:\quad\quad $X^i{}_iX^k{}_i-q^2 X^k{}_iX^i{}_i=\mu
q \sum_{a<i}X^k{}_aX^a{}_i$.

(iii.a) For $i<k$:\quad\quad
$X^i{}_jX^k{}_l-\sigma_{ijk}\bar\sigma_{ijl}X^k{}_lX^i{}_j=\mu\theta_{jl}
X^k{}_jX^i{}_l \sigma_{ijk}$.

(iii.b) For $j<l$:\quad\quad
$qX^i{}_jX^i{}_l=\bar\sigma_{ijl}X^i{}_lX^i{}_j$.

(iii.c) For $i<k$:\quad\quad $X^i{}_j
X^k{}_j=q\sigma_{ijk}X^k{}_jX^i{}_j$.

(iii.d) For $i<j$:\quad\quad
$qX^i{}_jX^j{}_l-\bar\sigma_{jli}X^j{}_lX^i{}_j=-\mu
\sum_{a<j}X^i{}_aX^a{}_l+\mu\theta_{jl} X^j{}_jX^i{}_l$.

(iii.e) For $i<k$:\quad\quad
$\bar\sigma_{ijk}X^i{}_jX^k{}_i-qX^k{}_iX^i{}_j
=\mu\sum_{a<i}X^k{}_aX^a{}_j+\mu\theta_{ji}X^k{}_jX^i{}_i$.

(iii.f) For $i<j$: \begin{eqnarray*} X^i{}_jX^j{}_i-X^j{}_iX^i{}_j
&=& q\mu\sum_{a<i}X^a{}_jX^j{}_aq^{2(a-i)}+q^{-1}\mu
X^j{}_jX^i{}_i
-q^{-1}\mu\sum_{a<j}X^i{}_aX^a{}_i\\
&&+\mu^2\sum_{b<j;a<i}X^a{}_bX^b{}_aq^{2(a-i)} -\mu^2
X^j{}_j\sum_{a<i}X^a{}_aq^{2(a-i)}.\end{eqnarray*}
\end{proposition}
\Proof We write each $R$ as the sum of an $L$ and an $M$ term as
explained above, giving 16 terms for $\Ups$. We use a standard
graphical notation to follow the values forced for the summed
indices by the $\delta$-functions in $L$ and $M$. We then use the
above lemma to break down the results further, to obtain:
\begin{eqnarray*}&&\kern-20pt \Ups(X^i{}_j\otimes X^k{}_l)
= X^k{}_l\otimes X^i{}_j(1-\delta_{kl})\left(\sigma_{ijk}\bar\sigma_{ijl}
+q^{-1}(\delta_{ki}+\delta_{kj})\bar\sigma_{ijl}
+q(\delta_{li}+\delta_{lj})\sigma_{ijk}\right)\\
&&-\mu \sum_{a<k}X^a{}_l\otimes X^i{}_a \delta_{kj} q^{2(a-k)}\bar\sigma_{ial}
+\mu\sum_{a<i}X^k{}_a\otimes X^a{}_j\delta_{il}\sigma_{ijk}\\
&& +\mu X^k{}_j\otimes X^i{}_l\theta_{jl}\sigma_{ijk}
-\mu X^i{}_l\otimes X^k{}_j\theta_{ik}\bar\sigma_{kjl}\\
&&+X^k{}_l\otimes X^i{}_j\left(\delta_{ij}+\delta_{kl}-\delta_{ij}\delta_{kl}
+\mu\delta_{ij}(q \delta_{jl}-q^{-1}\delta_{jk})+\delta_{ki}\delta_{jl}
+\delta_{kj}\delta_{il}-(2+\mu^2)\delta_{ij}\delta_{jk}\delta_{kl}\right)\\
&&-q\mu \sum_{a<k}X^a{}_l\otimes X^i{}_aq^{2(a-k)}\delta_{kj}\delta_{il}
-q\mu X^l{}_l\otimes X^i{}_l q^{2(l-k)}\theta_{kl}\delta_{kj}
-q^{-1}\mu X^i{}_l\otimes X^i{}_i q^{2(i-k)}\theta_{ki}\delta_{kj}\\
&&+\mu(q+q^{-1})X^i{}_i\otimes X^i{}_i q^{2(i-k)}\theta_{ki}\delta_{il}\delta_{kj}\\
&&+\mu\left(X^k{}_j\otimes X^i{}_l\theta_{jl}
+\sum_{a<i}X^k{}_a\otimes X^a{}_j\delta_{il}\right)\left(q^{-1}\delta_{ik}+q^{-1}\delta_{jk}
+q\delta_{ij}-(q+q^{-1})\delta_{ij}\delta_{jk}\right)\\
&&-\mu X^i{}_l\otimes
X^k{}_j\theta_{ik}\left(q^{-1}\delta_{jk}+q\delta_{jl}
+q\delta_{kl}-(q+q^{-1})\delta_{jk}\delta_{kl}\right)\\
&&+\mu^2\sum_{a<i}X^k{}_a\otimes X^a{}_l\theta_{jl}\delta_{ij}
-\mu^2\sum_{a<k}X^i{}_a\otimes X^a{}_j\theta_{ik}\delta_{kl}
+\mu^2 X^i{}_l\otimes\sum_{a<i,k}X^a{}_a q^{2(a-k)}\delta_{jk}\\
&&-\mu^2\sum_{k>a>l}X^a{}_a\otimes X^i{}_l q^{2(a-k)}\delta_{jk}
-\mu^2\sum_{a<k;\ b<i} X^a{}_b\otimes X^b{}_a
q^{2(a-k)}\delta_{il}\delta_{jk}-\mu^2 X^i{}_j\otimes
X^k{}_l\theta_{ik}\theta_{jl}.
\end{eqnarray*}
Next we break up the $\left({n^2\atop 2}\right)$ relations into
convenient special cases as stated. We then compute the relation
$X^i{}_j X^k{}_l=\cdot\Ups(X^i{}_j \otimes X^k{}_l)$ for $i,j,k,l$
according to the leading term on the left hand side in each of
the cases stated. We then simplify the resulting set of
equations. In some cases the simpler version arises from
$X^k{}_lX^i{}_j=\cdot\Ups(X^k{}_l\otimes X^i{}_j)$ instead.  \qed

Among other things, one may verify what is known from general
R-matrix methods for braided matrices\cite{Majid-book} that the
right-invariant $q$-trace element
\[\utr=\sum_i X^i{}_i q^{2i}\]
is central.

The braided Lie bracket of the multiparameter
$\widetilde{sl_q(n)}$ may be computed as $(\id\otimes\eps)\Ups$ or
directly from the R-matrix relations in \cite{Majid-braided-Lie}.
In the $q$-Hecke case these reduce to
\[ [X^i{}_j,X^k{}_l]=\delta^i{}_jX^k{}_l-\mu q q^{-2j} \delta^k{}_j
X^i{}_l+\mu R^{-1}{}^a{}_m{}^i{}_b R^n{}_c{}^b{}_l \tilde
R^c{}_j{}^k{}_a X^m{}_n\] and the matrix coalgebra structure. Note
also that $[\ ,\ ]$ necessarily closes on $\ker\eps$ which should
be thought of as the infinitesimal elements of the braided Lie
algebra (the classical model of a braided Lie algebra in
\cite{Majid-braided-Lie} is $\calL=k\oplus\gg$ with $\gg=\ker\eps$
a classical Lie algebra). We write $sl_q(n)=\ker\eps$ inside
$\widetilde{sl_q(n)}$.

\begin{proposition} Suppose that $\eps(\utr)\ne 0$, i.e. $q^{2n}\ne
1$. Then $\widetilde{sl_q(n)}=k\utr\oplus sl_q(n)$ in the braided
category, where
\[ sl_q(n)={\rm span}\{X^i{}_j,\ h_k\ |\ i\ne
j;\  k=1,\cdots,n-1\};\quad h_i:= X^i{}_i-X^{i+1}{}_{i+1}.\]
The q-Lie brackets are as follows. We  let
\[ H_i:= \sum_{a<i}[a]_{q^2}h_a;\quad
[a]_{q^2}:= {1-q^{2a}\over 1-q^2}.\]

(i) $ [x,\utr]=0,\quad [\utr,\utr]=\eps(\utr)\utr,\quad [\utr,
x]=\eps(\utr)\lambda x,\quad \forall x\in sl_q(n), \quad
\lambda=1+\mu^2.$

(ii.a) `Cartan' relations
\begin{eqnarray*}[h_i,h_i] &=&-\mu^2[2]_{q^2}q^{-2i}H_i-\mu^2
q^{-2i}[i+1]_{q^2}h_i, \quad [h_i,h_{i+1}]=\mu^2 q^{-2i}H_{i+1}\\
{}[h_i,h_{i-1}] &=&\mu^2 q^{-2(i-1)}H_i,\qquad\qquad\quad
[h_i,h_j]=0,\quad\forall |i-j|>1.\end{eqnarray*}

(ii.b) `Weight' relations for $k\ne l$:
\[ [h_i,X^k{}_l]=\mu
X^k{}_l\left(q^{-2i}(q^{-1}\delta_{k,i+1}-q\delta_{ki})+q^{-1}
\delta_{li}-q\delta_{l,i+1}\right)\]
\[ [X^k{}_l,h_i]=\mu
X^k{}_l\left(q^{-1}\delta_{ki}-q\delta_{k,i+1}+q^{-2i}(q^{-1}\delta_{l,i+1}-q\delta_{li})\right)\]

(ii.c) `Root' relations for $i\ne j$, $k\ne l$:
\[ [X^i{}_j,X^k{}_l]=-q\mu
q^{-2j}\delta_{kj}X^i{}_l+\mu\delta_{il}\sigma_{ijk}X^k{}_l
+q\mu\delta_{il}\delta_{jk}q^{-2k}X^k{}_k-\mu^2\delta_{il}\delta_{jk}q^{-2(k-1)}H_k.\]
\end{proposition}
\Proof Working from either $\Ups$ in the proof of the the
preceding proposition or from the R-matrix formula, we obtain
\[ [X^i{}_j,X^k{}_l]=\delta_{ij}X^k{}_l-q\mu q^{-2j}
\delta_{kj}X^i{}_l+\mu\delta_{il}(\sigma_{ijk} +
q\delta_{ij}+q^{-1}(\delta_{ik}+\delta_{jk})
-(q+q^{-1})\delta_{ij}\delta_{jk}\]
\[\quad\quad-\mu^2\delta_{il}\delta_{jk}\sum_{a<k}X^a{}_a
q^{2(a-k)}-\mu^2\delta_{kl}\theta_{ik}X^i{}_j+\mu^2
\delta_{ij}\theta_{jl} X^k{}_l.\] One finds in particular that for
all $k\ne l$,
\[ [X^i{}_i,X^k{}_k]=X^k{}_k(1-q\mu\delta_{ki}
q^{-2i}+q^{-1}\mu\delta_{ik}+\mu^2\theta_{ik}) -\mu^2\delta_{ik}
\sum_{a<k}q^{2(a-k)}X^a{}_a-\mu^2\theta_{ik}X^i{}_i\] which gives
the `Cartan' relations after further computation. Similarly
\[
[X^i{}_i,X^k{}_l]=X^k{}_l(1+q\mu\delta_{il}-q\mu
q^{-2i}\delta_{ik}+\mu^2\theta_{il}),\] etc give the other
relations. The $H_i$ arise from the summed $X^a{}_a$ terms
written in terms of the $h_j$. \qed

It is also possible to present the (ii.b) and (ii.c) relations
above in terms of  generators $X_i:= X^{i+1}{}_i$, $Y_i:=
X^i{}_{i+1}$ with other `root vectors' $X^i{}_j$ generated by
repeated Lie brackets of these. Among these, we have (from the
above):
\[\begin{array}{rll}{}[h_i,X_{i-1}]&=-q\mu q^{-2i}X_{i-1}&=-q^{-2(i-1)}\ [X_{i-1},h_i] \\
{}[h_i,X_i]&=q^{-1}\mu(1+q^{-2i})X_i&=-q^{-2}\ [X_i,h_i] \\
{}[h_i,X_{i+1}]&=-q\mu X_{i+1}&=-q^{2(i+1)}\ [X_{i+1},h_i]\\
{}[h_i,Y_{i-1}]&=q^{-1}\mu Y_{i-1}&=-q^{2(i-1)}\ [Y_{i-1},h_i] \\
{}[h_i,Y_i]&=-q\mu (1+q^{-2i})Y_i&=-q^2\ [Y_i,h_i] \\
{}[h_i,Y_{i+1}]&=q^{-1}\mu q^{-2i}Y_{i+1}&=-q^{-2(i+1)}\ [Y_{i+1},h_i] \\
{}[X_i,Y_j]&=\mu \delta_{ij}q^{-2i+1}(h_i-q\mu H_i)& \\
{}[Y_j,X_i]&=-\mu \delta_{ij}q^{-2i-1}(q^{2i}h_i+q\mu
H_{i}).&\end{array}\] The above results reduce for $n=2$ to the
computations in \cite[Ex. 5.5]{Majid-braided-Lie}, where $h=a-d$,
$X=c$, $Y=b$ in the notation there. The classical $q\to 1$ limit
should of course be taken after rescaling all the generators by
$\mu^{-1}$. Similarly in $BM_q(n)$ we would obtain a commutative
algebra (the coordinate algebra of the space of $n\times n$
matrices) without rescaling.

\subsection{Generalised Lie algebras $sl_{q}(n)$}

In this concluding section observe that a different quotient of
the braided enveloping algebra or braided matrices
$B(\widetilde{sl_q(n)})=BM_q(n)$ gives what could reasonably be
called the enveloping algebra of `generalised Lie algebras'
$sl_q(n)$ of the type suggested by representation
theory\cite{LS}\cite{Delius}. Indeed, we have already seen above
that $\widetilde{sl_q(n)}=k\oplus sl_q(n)$ where the $k$ is
spanned by $\utr$ and $sl_q(n)=\ker\eps$ cf.
\cite{Majid-braided-Lie}. Since the braided Lie bracket restricted
to $sl_q(n)$ is covariant it must also coincide with the
`generalised Lie bracket' defined via the q-deformed adjoint
representation in the representation theory approach.

In fact we are in the situation of Section~3.3, i.e. the braided
Lie algebra is split. As we explain now, this is a general feature
of the setting of Section~5.1 with $C_1=k\oplus C$ and $C$ simple
(this includes in principle all simple FODC over standard quantum
groups, although clearly the case $\calO_q(SL(n))$ with its
$n^2$-dimensional calculus is the most relevant). Let
$\utr=\tilde{R}^j{}_a{}^a{}_i\,{\l^i}_j$ be the right quantum
trace of $C$. We have to assume that $\e(\utr)\ne 0$. It is
$\ad_R$-invariant ($\ad_R(\utr)=\utr\otimes 1$), therefore for all
$a\in A$ it satisfies $a\otimes \utr\mapsto \utr\otimes a$ for any
braiding associated to a right crossed module $(A,?,\ad_R)$. By
the intertwining properties of $\qqq_1$, the element
$c=\qqq_1(\utr)/\e(\utr)\in\calL$ is central in $A^\circ$,
satisfies $\e(c)=1$, $\Psi(-\otimes c)=c\otimes -$, $\Psi(c\otimes
-)=-\otimes c$ and $\Ups(-\otimes c)=c\otimes -$. Finally, since
$c$ is central in $A^\circ$ and since $\calL^+:= \ker\e_\calL$ is
simple for the left adjoint action (since $C=k\utr\oplus C^+$ is a
semi-simple $A$-comodule for the adjoint coaction), one must have
$[c,x]=\l\, x$ for all $x\in \calL^+$, where $\lambda$ is a
constant, which we assume $\ne 0$. In this case the braided Lie
algebra $\calL$ has a distinguished decomposition
$$ \calL=k\,c\oplus\calL^+$$ and we have a
decomposition of the canonical braiding $\Ups$ of $\calL$ as in
Section~3.3, \ie $\Ups(z\otimes c) = c\otimes z$ for all $z\in\calL$, and
\begin{equation}
\Ups(x\otimes y)=\o(x\otimes y) +[x,y]\otimes c,
\quad
\Ups(c\otimes x) = \l\,x\otimes c+\rho(x)
\end{equation}
for all $x,y\in \calL^+$, for  uniquely determined maps
$\o:(\calL^+)^{\otimes 2}\to (\calL^+)^{\otimes 2}$ and
$\rho:\calL^+\to (\calL^+)^{\otimes 2}$. Moreover, if $c$ is not a
zero divisor then $(\calL^+,[\,,\,])$ is among other things a
generalized Lie algebra in the sense of \cite{LS} with generalized
antisymmetrizer $\l^{-1}(\id-\o)$. The zero divisor condition
holds in the multiparameter case by \cite{LeB-94}.

Next, we have for any split braided Lie algebra with $\calL^+$
simple its reduced enveloping algebra $B_{red}(\calL^+)$ as
explained in Section~3.3. In our case of interest, it means
\[ B_{red}(sl_q(n))=B(\widetilde{sl_q(n)})/\langle
\utr-\eps(\utr)\lambda\rangle.\]  As explained in Section~3.3 the
relations of $B_{red}(sl_q(n))$ contain the defining
`antisymmetrizer' relations of the `eneveloping algebra' of the
algebra $U_{LS}(sl(n))$ (say) proposed in \cite{LS} but in
principle could contain further relations. One may check that at
least for $n=2$ the two constructions do coincide. This is the
algebra
\[ q^{-2}hX-Xh=\lambda (1-q^{-4})X,\quad q^2hY-Yh
=-\lambda q^2(1-q^{-4})Y,\quad [X,Y]= {q^2-1\over
q^2+1}h^2+\lambda (1-q^{-2})h\] isomorphic after rescaling the
generators to the Witten algebra $W_{q^2}(sl(2))$ as noted in
\cite{LeB-95}. Hence we propose (multiparameter)
$B_{red}(sl_q(n))$ as a generalisation of the Witten algebra for
$n\ge 2$ and the (multiparameter) braided matrices $BM_q(n)$ in
Section~5.2 as its homogenisation.

\bigskip
\bigskip
{\scriptsize SCHOOL OF MATHEMATICAL SCIENCES, QUEEN MARY,
UNIVERSITY OF LONDON, MILE END RD, LONDON E1 4NS, UK}
\bigskip
%
%
%
%
%
%

\end{document}